\documentclass[final]{siamltex}
\usepackage{amsmath,amssymb,graphicx}
 \newtheorem{remark}[theorem]{Remark}

\title{Numerical study of shock formation in the dispersionless 
Kadomtsev-Petviashvili equation and dispersive 
regularizations\thanks{We thank B.~Dubrovin, T.~Grava, N.~Tzvetkov 
and A.~Weideman for helpful discussions and hints. 
This work has been supported by 
the project FroM-PDE funded by the European
Research Council through the Advanced Investigator Grant Scheme,
the ANR via the program ANR-09-BLAN-0117-01,  the Austrian Science Foundation FWF, project SFB F41
("VICOM") and project I830-N13 ("LODIQUAS"). We are grateful for 
access to the HPC resources from 
GENCI-CINES/IDRIS (Grant 2013-106628) on which part of the 
computations in this paper has been done, and to the Vienna Scientific Cluster (VSC).}}

\author{C.~Klein\thanks{Institut de Math\'ematiques de Bourgogne,
		Universit\'e de Bourgogne, 9 avenue Alain Savary, 21078 Dijon
		Cedex, France
    ({\tt christian.klein@u-bourgogne.fr})}
\and
K.~Roidot\thanks{SISSA, Via Bonomea 265, I-34136 Trieste, Italy
    ({\tt kristelle.roidot@sisssa.it})}
}

\begin{document}
\maketitle

\begin{abstract}
    The formation of singularities in solutions to the dispersionless 
    Kadomtsev-Petviashvili (dKP)  equation is studied numerically for 
    different classes of initial data. The asymptotic behavior of the Fourier 
    coefficients is used  to 
    quantitatively identify the critical time and location and the 
    type of the singularity. The approach is first tested in detail 
    in $1+1$ dimensions for the known  case of 
    the Hopf equation, where it is shown that the break-up of the 
    solution can be identified 
    with prescribed accuracy. For 
    dissipative regularizations of this shock formation as the 
    Burgers' equation and for dispersive 
    regularizations as the Korteweg-de Vries equation, the Fourier 
    coefficients indicate as expected global regularity of the 
    solutions. The Kadomtsev-Petviashvili (KP) equation can be seen as a 
    dispersive regularization of the dKP equation. The behavior of KP 
    solutions for small dispersion parameter $\epsilon\ll 1$ 
 near a break-up of corresponding dKP solutions is 
    studied. It is found that the difference between KP and dKP 
    solutions for the same initial data at the critical point scales roughly as $\epsilon^{2/7}$ 
    as for the Korteweg-de Vries equation. 

\end{abstract}

\begin{keywords}
    Asymptotic Fourier analysis, Kadomtsev-Petviashvili equation, Burgers' equation, Korteweg-de Vries 
    equation, 
    dispersive shocks
\end{keywords}

\begin{AMS}
    Primary, 65M70; Secondary, 65L05, 65M20
\end{AMS}

\section{Introduction}
Nonlinear evolution equations without dispersion and dissipation 
generically have solutions which show the \emph{wave breaking} 
phenomenon, 
i.e., the formation of a shock, a gradient catastrophe in finite 
time. A standard example in this context is the Hopf equation
\begin{equation}
    \label{Hopf}
u_t+6uu_x=0.
\end{equation}
It is well known that the solution of (\ref{Hopf}) for an initial value problem $u(x,0)=u_0(x)$
can be obtained via the method of characteristics in the implicit form
\begin{equation}
\label{Hopfsol}
u(x,t)=u_0(\xi),\quad x=6tu_0(\xi)+\xi.
\end{equation}
If the initial data are such that
$t_c=\dfrac{1}{\max_{\xi\in\mathbb{R}}[-6u'_0(\xi)]}$ is positive,
the solution reaches a point of gradient catastrophe $x_{c}$ at $t_{c}$ where the 
derivative of the Hopf solution blows up, but where the solution 
stays finite. 

Regularizations of this equation with small dissipation $\epsilon$, the Burgers' 
equation,
\begin{equation}
    u_{t}+6uu_{x}=\epsilon u_{xx}
    \label{burgers},
\end{equation}
or small dispersion $\epsilon$, the Korteweg-de Vries (KdV) equation,
\begin{equation}
    u_{t}+6uu_{x}+\epsilon^{2}u_{xxx}=0
    \label{KdV},
\end{equation}
will 
have solutions which stay regular at the shock of the Hopf solution 
for the same initial data, but show some critical behavior in the 
vicinity of the point $(x_{c},t_{c})$. At the critical point, the difference 
between Hopf solution and the solution of the regularized equation shows 
a characteristic scaling in $\epsilon$, for KdV $\epsilon^{2/7}$. 
Dubrovin \cite{Dub06} conjectured a universal behavior of solutions to Hamiltonian 
regularizations of the Hopf equation (among which KdV is the most 
prominent example) in the vicinity of a shock.  
In \cite{GK,GK08,DGK11,GK12} strong 
numerical evidence for this conjecture was given, which was proven 
for KdV in \cite{CG} via Riemann-Hilbert techniques. For $t\ll t_{c}$ 
this difference  scales as $\epsilon^{2}$. An asymptotic description 
of dissipative regularizations was presented in \cite{DE13}.

Dubrovin's conjecture is based on a double scaling limit where $\epsilon\to 0$ 
and simultaneously $x\to x_{c}$, $t\to t_{c}$ in such a way that the limits
\[
\lim_{\begin{matrix}
\epsilon\rightarrow 0\\
x-6u_ct\rightarrow x_c-6u_ct_c
\end{matrix}}\left[ \dfrac{x-x_c-6u_c(t-t_c)}{\epsilon^{\frac{6}{7}}}\right] \,\,\mbox{and}\,\,    \\    \quad 
\lim_{\begin{matrix}
\epsilon\rightarrow 0\\
t\rightarrow t_c
\end{matrix}}\left[ \dfrac{(t-t_c)}{\epsilon^{\frac{4}{7}} }\right], \;\;\; 
\]
where $u_c=u(x_c,t_c)$, exist and are bounded. In \cite{GK,GK08,DGK11,GK12} it 
was shown that it is indeed possible to study these scalings in $\epsilon$ 
numerically if the critical point $(x_{c},t_{c})$ is known. The above 
formulae make it clear that the error in the values $x_{c}$, $t_{c}$ 
must be smaller than the considered values of $\epsilon$ if the correct 
scalings are to be identified numerically. This implies that for a 
numerical study of critical scaling phenomena, a control of the error 
bounds for the critical point is crucial. 

It is mathematically and physically interesting to study 
similar problems in higher dimensions where much less about shock 
formation and dissipative or dispersive 
regularizations is known. A breakdown of regularity 
always indicates a limit of the applicability of the studied model. It is 
at such points that effects as dissipation and dispersion, which have 
been neglected in the simplified model, become important. The main 
technical problem in higher dimensions is that the dispersionless 
system even of completely integrable equations as 
the Kadomtsev-Petviashvili (KP) equation \cite{KP}, 
a $2+1$-dimensional variant of 
KdV, is not integrable in the usual sense. This means that there are 
in general
no standard solution generating techniques as hodograph methods or 
(linear) Riemann-Hilbert 
problems available in this context (but see \cite{MS08} for a 
nonlinear Riemann-Hilbert approach). The KP equations read
\begin{equation}\label{e1}
\partial_{x}\left(\partial_{t}u+6u\partial_{x}u+\epsilon^{2}\partial_{xxx}u\right)+\lambda\partial_{yy}u=0,\,\,\lambda=\pm1
\end{equation}
where $(x,y,t)\in\mathbb{R}_{x}\times\mathbb{R}_{y}\times\mathbb{R}_{t}$
and where $\epsilon\ll1$ is a small scaling parameter. 
The case $\lambda=-1$ corresponds to the KP I model with a 
\emph{focusing} effect, and the case $\lambda=1$
corresponds to the KP II model with a \emph{defocusing} effect.

The focus in this paper will be on the dispersionless variant of the 
KP equation, the dKP equation, also known as Khokhlov-Zabolotskaya 
equation \cite{ZK69}, which follows 
from (\ref{e1}) for $\epsilon=0$,
\begin{equation}\label{dKP}
\partial_{x}\left(\partial_{t}u+6u\partial_{x}u\right)+\lambda\partial_{yy}u=0,\,\,\lambda=\pm1 .
\end{equation}
It appears in many applications as a model for dissipation- and 
dispersionless, essentially one-dimensional waves with weak 
transverse effects,   for instance in nonlinear acoustic and in gaz 
dynamics, see \cite{ZK69,KP}. It also plays a role in the context of 
general relativity and differential geometry  \cite{AC,DMT01}.  
Alinhac and coworkers,  see for instance 
\cite{alinhac} and references therein, showed that it appears as a universal 
model in the geometric study of singularity formation in 
nonlinear wave equations. Existence of solutions up to break-up was 
proven, and the singularity was identified for a generic setting as a 
one-dimensional cusp. 

The dKP equation is not completely integrable in the sense that an
infinite number of conserved quantities exists (in fact it only has 
three). But it is integrable 
according to the definition of  \cite{FK04} that an infinite number of 
hydrodynamical reductions exists, see also \cite{Kod1, Kod2, Zak2}. 
As was shown in \cite{DMT01}, solutions can be constructed with 
Twistor methods in terms of Einstein-Weyl metrics. In \cite{MS06,MS07} 
solutions to the dKP equation were obtained via the solution of a 
nonlinear Riemann-Hilbert problem. This allowed the study of the 
long-time asymptotics and finite time break-up in  \cite{MS07, MS08}.
Another way to construct dKP solutions arises from the theory of 
Frobenius manifolds. In \cite{Rai12} the infinite dimensional 
Frobenius manifold corresponding to dKP was constructed. A common 
feature of the above solution generating techniques is that even the 
most explicit forms of the solutions require the solution of singular 
nonlinear integral equations and the inversion of implicit functions which makes it 
difficult to find explicit examples for the wave breaking. Since none of 
these approaches has been studied numerically so far, we solve here 
dKP directly up to the formation of a singularity.

As mentioned above, a 
precondition for the numerical study of scaling 
laws as for KdV is that the critical 
point and the critical solution can be determined 
with sufficient accuracy. It is the 
goal of this paper to provide the necessary tools for this in 
a rather general setting, and to apply them to the dKP equation for 
several classes of initial data.  
The applicability of this approach will be shown first for the 
$1+1$-dimensional case where exact solutions provide tests. Note that a dissipative regularization does 
not give precise information on the critical point as can be 
inferred from the example of the Burgers' equation with  small 
dissipation we will also consider. The same will be shown to be true for a dispersive 
regularization as provided by KdV. Thus it can be concluded that a 
direct study of the dispersion- and dissipationless equation with the 
techniques explored in this paper is more promising in the context 
of identifying singularity formation.

The basic idea of the used numerical approach
is to approximate the spatial dependence of the 
solutions by a discrete Fourier series. It is known that the 
asymptotic behavior (for large wave numbers $k$) of the Fourier 
coefficients is exponential for solutions analytic in the complex 
plane in the vicinity of the real axis. This dependence becomes
algebraic at the break-up 
of the solution. Thus the vanishing of the exponential decrease in 
the Fourier coefficients would indicate the appearance of the 
singularity. In practice there are several problems in this context: 
firstly it is difficult numerically to integrate the equation up to the 
critical point; moreover even if this can be done with sufficient 
accuracy, the Fourier coefficients might be polluted especially at 
the high wave numbers where the asymptotic behavior is to be read 
off; and last but not least, the asymptotic relations for the Fourier 
coefficients are established for a continuous Fourier transform, 
whereas numerically only a discrete Fourier transform is considered 
as an approximation. 

Asymptotic Fourier analysis was first applied numerically by Sulem, 
Sulem and Frisch \cite{SSF}, at the time of course with much lower 
resolution than is accessible today. Therefore the study in 
\cite{SSF} was mainly qualitative, but gave a proof of concept. 
In \cite{SSP}
they studied singular solutions to the two-dimensional cubic NLS equation. 
An application of the method to the 2-d Euler equations can be found 
in \cite{FMB, MBF}. It has also been applied to 
the study of complex singularities of the 3-d Euler equations
in \cite{CR}, in thin jets with surface tension \cite{PS98}, the 
complex Burgers' equation \cite{SCE96} and the Camassa-Holm equation \cite{RLSS}
to name just a few examples. The main goal in these 
works is to decide whether the studied solutions develop a 
singularity. However, the purpose of the present paper is to identify 
time and location of a singularity, which is known to appear in the 
solution, with prescribed accuracy with these methods.

The tracking of singularities can be of course 
also achieved with other techniques, for instance Pad\'e approximants 
in \cite{Wei}. The problem with the latter is, however, that this 
adds an additional dimension to the problem which is especially 
expensive if this technique is  to be applied in higher 
dimensions. Full singularity tracking for a function in 
$\mathbb{R}^{n}$ would imply the study of the function in 
$\mathbb{C}^{n}$, thus doubling the spatial dimensions. Therefore a 
quantitative numerical approach based on the asymptotic behavior of 
discrete Fourier series would be an economic way to study singularity 
formation in higher dimensions for functions analytic before a 
critical time $t_{c}$. It is one goal of this paper to present a 
careful investigation of the possibilities and the limitations of 
this approach, and to show that the method can provide information on 
critical points with the needed accuracy to study scaling laws. 
The techniques are developped for the example of the Hopf equation and 
various regularizations thereof, and will then be applied to the dKP 
equation for various initial data. The thus identified dKP solution 
at the critical point is used to study the scaling in the small 
dispersion parameter $\epsilon$ of corresponding KP solutions. It is 
found that the scaling of the difference between KP and dKP solutions is 
compatible with $\epsilon^{2/7}$. This suggests that the behaviour of 
the KP solutions in the singular direction
is similar to KdV solutions near the break-up of 
Hopf solutions. 

To establish asymptotic Fourier analysis as a quantitative tool to 
identify shock formation, we test the various numerical problems in the 
approach \cite{SSF} separately and 
show they can be controlled. Then we will establish that this method 
in fact can be used to determine the critical point and the critical 
solution to the Hopf equation. The paper is organized as 
follows: In Sec.~2 we collect some known facts about the asymptotic 
behavior of Fourier coefficients of a real function analytic in the 
vicinity of the real axis. We numerically integrate the 
Hopf equation up to the critical point for given initial data and 
compare it with the exact solution. For 
the latter we show  
how well the Fourier coefficients can be fitted to the expected 
asymptotic decrease. We perform this fitting during the 
numerical solution of the Hopf equation to determine the critical 
time. The same approach is applied to solutions to the 
Burgers's and KdV equation for the same initial data. In Sec.~3 we 
apply these techniques to the break-up in dKP solutions for two 
classes of initial data, localized in two dimensions or as the line 
solitons of KP, localized in one spatial direction and infinitely 
extended in the other. In both cases, 
the scaling of the difference between KP solutions in the small 
dispersion limit and dKP solutions at the 
critical point is studied. We add some 
concluding remarks in Sec.~4.

\section{1+1-dimensional case}
In this section we will address numerically $1+1$-dimensional 
examples from the family of Hopf equations. We first review some 
well known facts from asymptotic Fourier analysis. Then we will 
numerically integrate the Hopf equation for a concrete example 
up to shock formation. We discuss how the asymptotic behaviour of the 
Fourier coefficients can be used to identify the critical time and 
solution of the Hopf equation. A similar analysis is then performed 
for the Burgers' and the Korteweg-de Vries equation which can be seen 
as dissipative and respectively dispersive regularizations of the Hopf 
equation. 

\subsection{Asymptotic Fourier analysis}
In this subsection we review some well known facts on the asymptotic 
behavior of the Fourier transform of a real function which is analytic 
in a strip around the real axis. The resulting formulae will be then 
used to numerically track singularities of the solutions in the 
complex plane via their Fourier coefficients without having to deal 
with two real dimensions as in \cite{Wei}. This allows to 
determine when a singularity hits the real axis, i.e., when
the real solution becomes singular. 

Spectral methods are a powerful tool in the numerical solution of 
differential equations because of
their excellent approximation 
properties for smooth functions, see for instance \cite{Can}. It is well known that the error in 
the approximation of such a  function with a spectral approximant of 
order $N$, e.g.\ a polynomial of degree $N$, decreases faster than 
any power of $1/N$. This provides a huge advantage over methods with 
just an algebraic decrease of the error and is the basis of the 
efficiency of spectral methods, see also \cite{Fornb}. The most used such 
methods are the expansion of the studied function 
in terms of discrete Fourier series. For given analyticity 
properties of the studied function, precise mathematical statements on the 
asymptotic behavior  
of the Fourier coefficients exist. As we will detail below, these 
imply that the location of 
singularities in the complex plane can be obtained from a given Fourier 
series computed on the real axis. It should be possible in this way 
to numerically determine the type of the singularity, to trace it in 
the complex plane,  and to establish when the singularity hits the 
real axis as was done for the first time in \cite{SSF}.

The Fourier 
transform $\hat{u}$ of  $u(x)\in L^{2}(\mathbb{R})$ is defined as 
\begin{equation}
    \hat{u}(k) = \int_{\mathbb{R}}^{}u(x)e^{-ik x}dx
    \label{fourier}.
\end{equation}
A singularity in the complex plane of the form $u\sim 
(z-z_{j})^{\mu_{j}}$, $\mu_{j}\notin \mathbb{Z}$, 
with $z_{j}=\alpha_{j}-i\delta_{j}$ in the lower half 
plane ($\delta_{j}\geq 0$) implies with a steepest descent argument  for $k\to\infty$  
the following asymptotic behavior of the Fourier 
coefficients (for a detailed derivation see e.g.~\cite{asymbook}),
\begin{equation}
    \hat{u}\sim 
    \sqrt{2\pi}\mu_{j}^{\mu_{j}+\frac{1}{2}}e^{-\mu_{j}}\frac{(-i)^{\mu_{j}+1}}{k^{\mu_{j}+1}} e^{-ik\alpha_{j}-k\delta_{j}}
    \label{fourierasym}.
\end{equation}
Consequently for a single such singularity with positive $\delta_{j}$, the modulus of the Fourier 
coefficients decreases exponentially for large $k$. For 
$\delta_{j}=0$, i.e., a singularity on the real axis,  
the modulus of the Fourier coefficients has an algebraic dependence 
on $k$. 
If there are several singularities of this form at $z_{j}$, 
$j=1,\ldots,J$, there will be oscillations in the modulus of the 
Fourier coefficients for moderately large $k$.

To numerically compute a Fourier transform, it has to be approximated 
by a discrete Fourier series which can be done efficiently via a 
 fast Fourier transform (FFT), see e.g.~\cite{tref}. The discrete Fourier transform  of the 
 vector $\mathbf{u}$ with components $u_{j}=u(x_{j})$, where 
 $x_{j}=2\pi L j/N$, $j=1,\ldots,N$ (i.e., the Fourier transform on 
 the interval $[0,2\pi L]$ where $L$ is a positive real number) 
 will be always denoted by $v$ in 
 the following. There is no obvious analogue of 
 relation (\ref{fourierasym}) for a discrete Fourier series, but it 
 can be seen as an approximation of the latter, which is also the 
 basis of the numerical approach in the solution of the PDE. It is possible to 
 establish bounds for the discrete series, see for 
 instance \cite{arnold}. We will in this paper always fit the 
 discrete Fourier series for large $|k|$ to the 
 asymptotic formula (\ref{fourierasym}). 
 
 \begin{remark}\label{ms}
 Note that the minimal 
 distance in Fourier space is $m:=2\pi L/N$ with $N$ being the number of 
 Fourier modes and $2\pi L$ the length of the computational domain in 
 physical space. This defines the smallest distance which can be 
 resolved in Fourier space. All values of $\delta$ below this 
 threshold cannot be distinguished numerically from 0.  
 \end{remark}

\subsection{Numerical integration of the Hopf equation}\label{numhopf}
In this subsection we will study the integration of the Hopf equation 
for a concrete example up to the appearance of a gradient catastrophe. 
Since the method of characteristics gives an exact solution in this 
case in implicit form which can be obtained numerically in principle 
with machine precision, we can use this as a test case for the 
numerical integration of a nonlinear dispersionless equation up to the first 
critical point in the solution. We 
discuss how to obtain asymptotically reliable Fourier coefficients, and with which precision this can 
be achieved.

In the following we will consider the example 
$u_{0}(x)=\mbox{sech}^{2}x$ for which the gradient catastrophe 
appears at the critical time $t_{c}=\sqrt{3}/8\sim 0.2165$ and at the critical point $x_{c} = 
\sqrt{3}/2-\ln((\sqrt{3}-1)/\sqrt{2})\sim1.5245$. These initial data 
are motivated by the well known KdV soliton. They have the advantage 
that they belong to the class of rapidly decreasing 
functions which, for sufficiently large $L$, can be analytically 
continued as periodic functions within numerical precision (the 
quantity $L$ and thus the length of the interval $[-\pi L,\pi L]$ is chosen large 
enough that $u_{0}(\pm \pi L)$ is smaller than machine precision 
($10^{-16}$ in our case). This ensures that there is no Gibbs phenomenon and that the 
analytic behavior of the Fourier coefficients is not affected by 
effects at the boundary of the computational domain. 

\begin{figure}[htb!]
\begin{center}
      \includegraphics[width=0.6\textwidth]{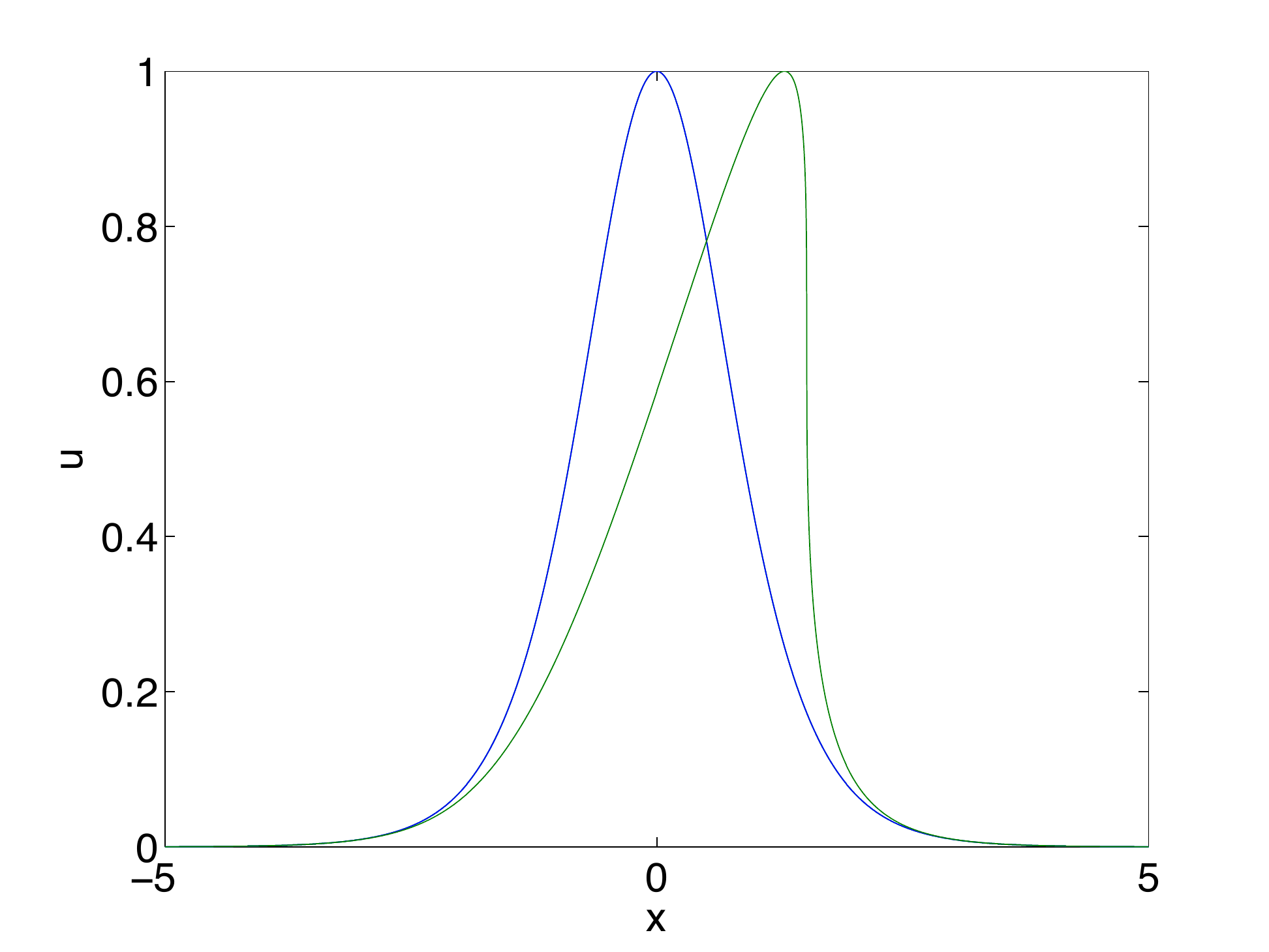}
     \caption{Solution to the Hopf equation (\ref{Hopf}) propagating to 
     the right for 
     $u_{0}(x)=\mbox{sech}^{2}x$ at $t=0$ in blue and at the critical 
     time $t_{c}=\sqrt{3}/8\sim 0.2165$ in green.}
       \label{hopfexact}
\end{center}
\end{figure}

Numerically we 
solve equation (\ref{Hopfsol}) iteratively,
$$\xi_{n+1}=x-6tu_{0}(\xi_{n})$$
with the initial guess $\xi_{0}=x$. If some relaxation is used, the 
iteration converges and is stopped once the residual of 
$6tu_0(\xi)+\xi-x$ is smaller than $10^{-13}$ in the $L_{\infty}$ norm.  
The iteration is  done in a vectorized and thus  efficient way, i.e., for all 
points at the same time. For the studied example the solution to the 
Hopf equation can be seen in Fig.~\ref{hopfexact} at the initial and 
at the critical time. 

To directly integrate the Hopf equation numerically, we use a 
standard Fourier spectral method in $x$: for $x\in [-\pi,\pi]L$, 
we take 
$L=5$, and the modulus of the Fourier coefficients computed via a discrete Fourier 
transformation decreases to machine precision for more than $N=2^{8}$ 
Fourier modes.  Treating the spatial 
dependence in (\ref{Hopf}) with a discrete Fourier transform, we 
obtain for the Hopf equation a system of ordinary differential 
equations (ODEs) in $t$ for the Fourier coefficients. Since the 
latter system is not stiff in contrast to the KdV equation treated in 
a later subsection, we can use standard time integrators here. We choose 
the well-known explicit fourth order Runge-Kutta scheme for 
convenience.

The accuracy of the numerical solution is tested via comparison with 
the exact solution. In addition we trace a conserved quantity of the 
Hopf equation, the energy,
$$E[u]=\int_{\mathbb{T}}^{}u^{3}dx,$$
which, when computed for the numerical solution, will depend on time due to 
unavoidable numerical 
errors. As shown for instance in \cite{ckkdvnls,KR}, the conservation of 
energy can be used as an indicator of numerical accuracy. We test 
this for the Hopf equation below by tracing the quantity 
$\Delta E = |1- 
\frac{E[u](t)}{E[u](0)}|$, where $E[u](t)$ is the numerically computed energy. 

\begin{figure}[htb!]
\begin{center}
      \includegraphics[width=0.6\textwidth]{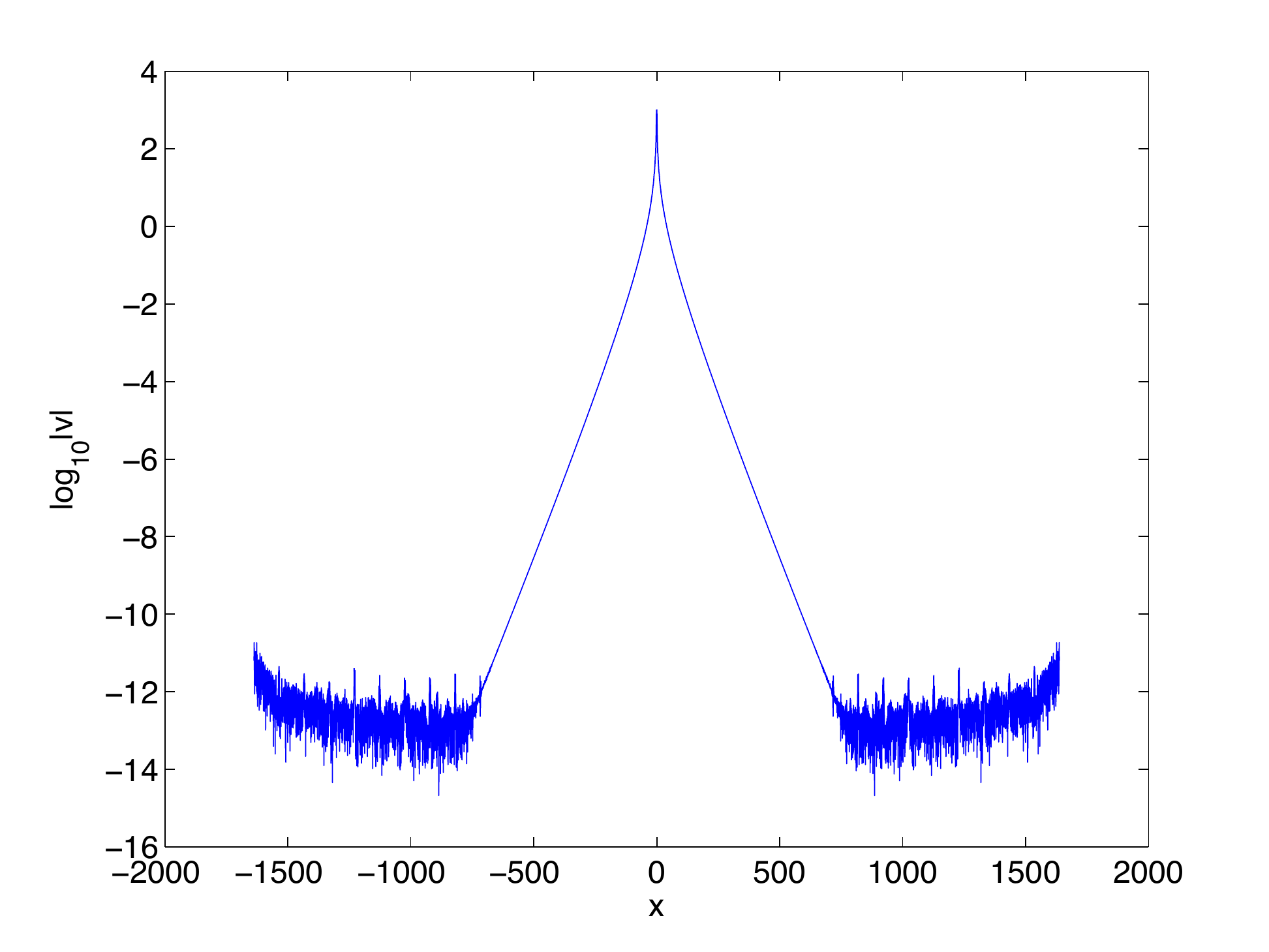}
     \caption{Modulus of the Fourier coefficients for the solution 
     to the Hopf equation (\ref{Hopf}) for 
     $u_{0}(x)=\mbox{sech}^{2}x$ at $t=0.18$.}
       \label{hopffourier_18}
\end{center}
\end{figure}

Since the solution of the Hopf equation is very different for times 
$t\ll t_{c}$ and for $t\sim t_{c}$, we will first numerically solve the 
Hopf equation up to $t=0.18 \ll t_{c}$. For $N=2^{14}$ Fourier modes, 
the solution at time $t=0.18$ is fully resolved as can be seen from 
Fig.~\ref{hopffourier_18}, where the modulus of the Fourier 
coefficients is shown. The modulus goes down to machine precision. 
For a time step $\Delta_{t}=7.2*10^{-5}$ the $L_{\infty}$ norm of the 
difference between numerical solution and exact solution is  
$6.4*10^{-11}$, the quantity $\Delta E\sim 10^{-14}$. 
For 
$\Delta_{t}=3.6*10^{-5}$, the respective numbers are 
$||u_{num}-u_{exact}||_{\infty}=3.7*10^{-12}$ and $\Delta E\sim 
10^{-15}$, for $\Delta_{t}=1.8*10^{-5}$ one finds 
$||u_{num}-u_{exact}||_{\infty}=8.4*10^{-13}$ and $\Delta E\sim 
10^{-14}$. Thus $\Delta E$ 
overestimates the numerical precision by roughly two to three orders of 
magnitude (once it is of the order of machine precision, it is 
obviously dependent on rounding errors and can even increase with 
smaller resolution as above). It
can be in fact used as an indicator of numerical 
accuracy in cases for which no exact solution is known, provided that 
there is sufficient spatial resolution. To sum up, the solution can be obtained up to machine precision for 
times $t\ll t_{c}$. 

For larger $t$, time steps of different size might be convenient. 
Thus we run the code with a certain time step up to 
time $t=0.18$ and use the result as initial data for a second time 
evolution up to the critical time $t_{c}$ (since we have an exact 
solution for the studied example, it could be also used as exact initial 
data at that time; but 
since the goal is to develop and test methods for a general case for which no 
exact solution is known, this is not done here).  The important point 
is that the solution can be 
obtained with machine precision at any time $t\ll t_{c}$ 
numerically. 

\begin{figure}[htb!]
  \includegraphics[width=0.5\textwidth]{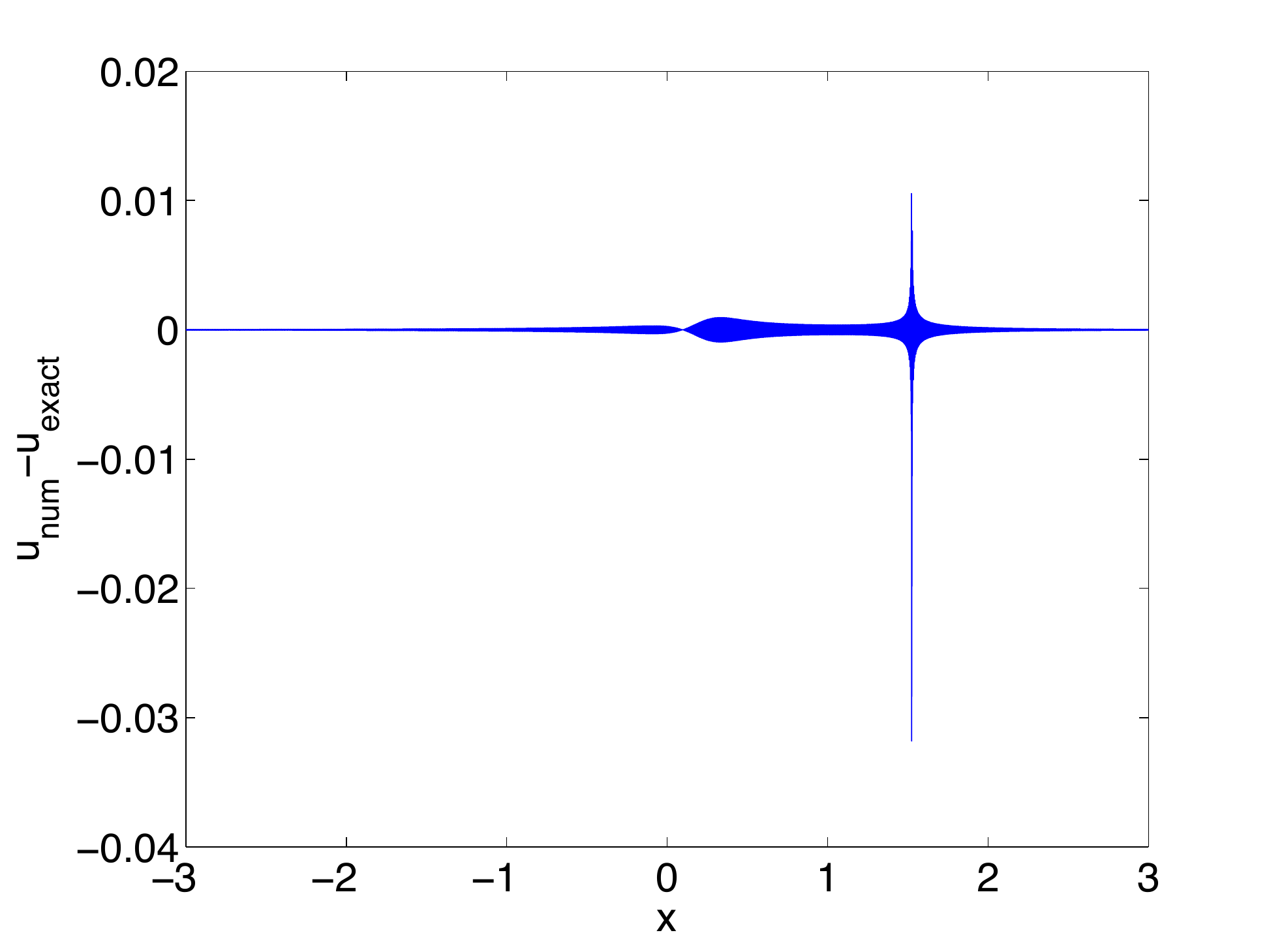}
  \includegraphics[width=0.5\textwidth]{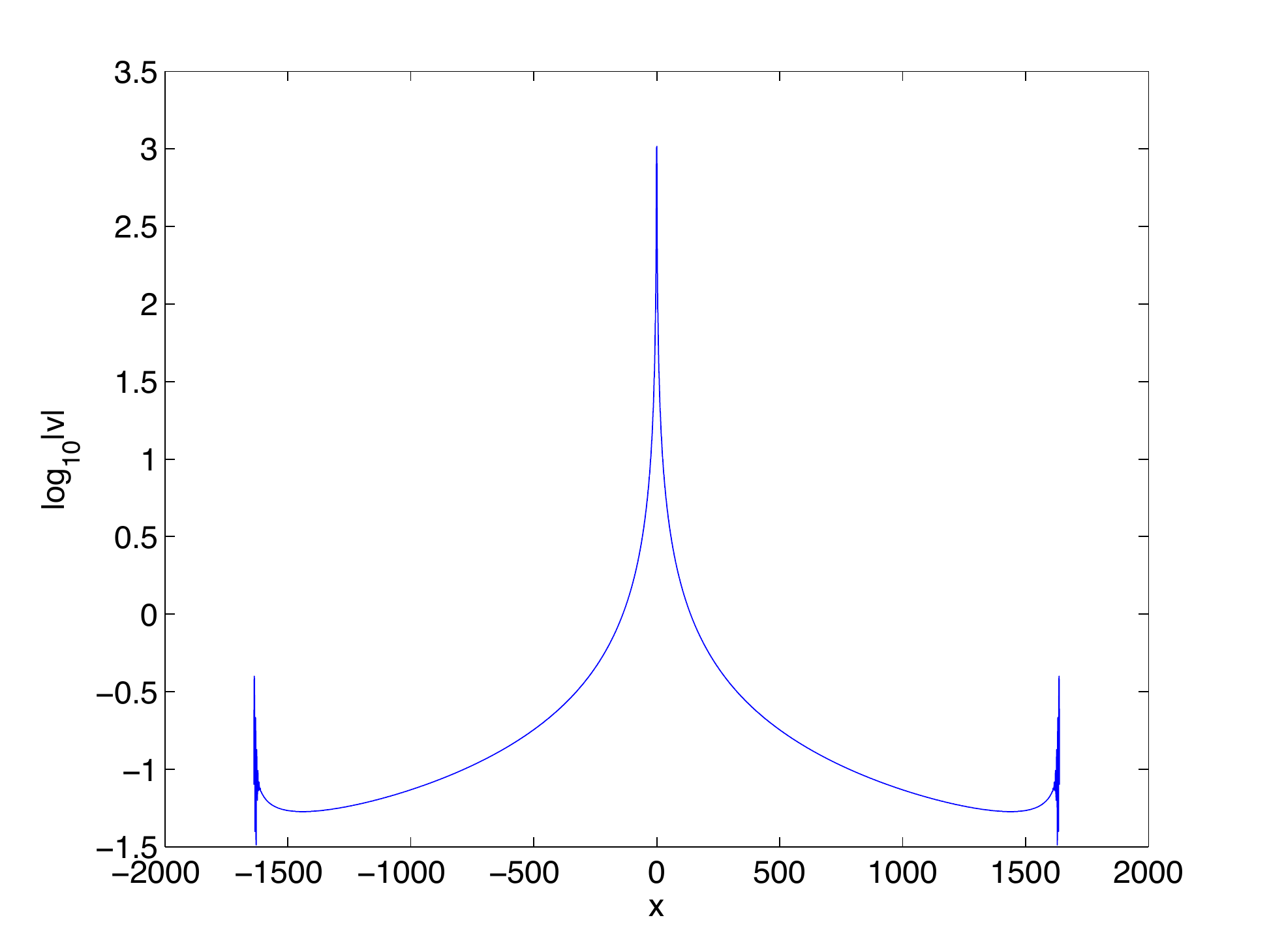}
 \caption{Difference of the numerical and exact solution at time 
 $t_{c}$ for $N=2^{14}$ Fourier modes on the left, and the modulus of the 
 Fourier coefficients for the numerical solution on the right.}
   \label{hopfc14}
\end{figure}

For the studied example, one gets for $N=2^{14}$ with 
$\Delta_{t}=7.3*10^{-5}$ an 
$L_{\infty}$ norm of the difference between numerical and exact 
solution of the order of 0.032. This error (as well as the $L_{2}$ 
norm of the difference between numerical and exact solution) does not 
decrease for smaller time steps. This shows that the error, which is 
always biggest at the critical point as can be seen in 
Fig.~\ref{hopfc14}, is due to a lack of resolution in $x$ and not in 
time. This lack of resolution is clear from the Fourier coefficients 
which can be also seen in Fig.~\ref{hopfc14}. Visibly the Fourier 
coefficients for high wave numbers do 
not show the expected decrease 
with $k$ due to numerical errors at the critical time. This will 
be important for the asymptotic analysis in the following section.

Note that the quantity $\Delta E$ for the computed energy decreases 
to machine precision if the time step decreases to 
$\Delta_{t}=9.13*10^{-6}$. This shows that this quantity is only a useful indicator of the 
numerical accuracy for sufficient spatial solution, i.e., if the 
Fourier coefficients decrease to the wanted accuracy. Lower values of 
$\Delta E$ than the modulus of the Fourier coefficient of the 
highest wave numvers are 
meaningless.
To obtain higher accuracies, it appears necessary to allow for higher 
spatial resolution. For $N=2^{15}$ Fourier modes and 
$\Delta_{t}=3.7*10^{-5}$, we 
get 
$ \| u_{num} - u_{exact} \|_{\infty} = $ 
$0.025$, 
for $N=2^{16}$ and $\Delta_{t}=1.8*10^{-5}$ 
$ \| u_{num} - u_{exact} \|_{\infty} = $ we find
$0.016$. The time steps have been chosen in a way 
that the error does not decrease further if $\Delta_{t}$ is decreased. 

\subsection{Asymptotic fitting of the Fourier coefficients}
In this subsection we will study the asymptotic behavior of the Fourier 
coefficients for the solution to the Hopf equation discussed in the 
previous subsection. We show how to apply the asymptotic formula 
(\ref{fourierasym}) to the exact solution before
and at the critical time. This will be done here first for 
the exact solution to avoid the problems with errors in the 
numerically determined Fourier coefficients at the critical time for 
high wave numbers. 

\begin{figure}[htb!]
  \includegraphics[width=0.5\textwidth]{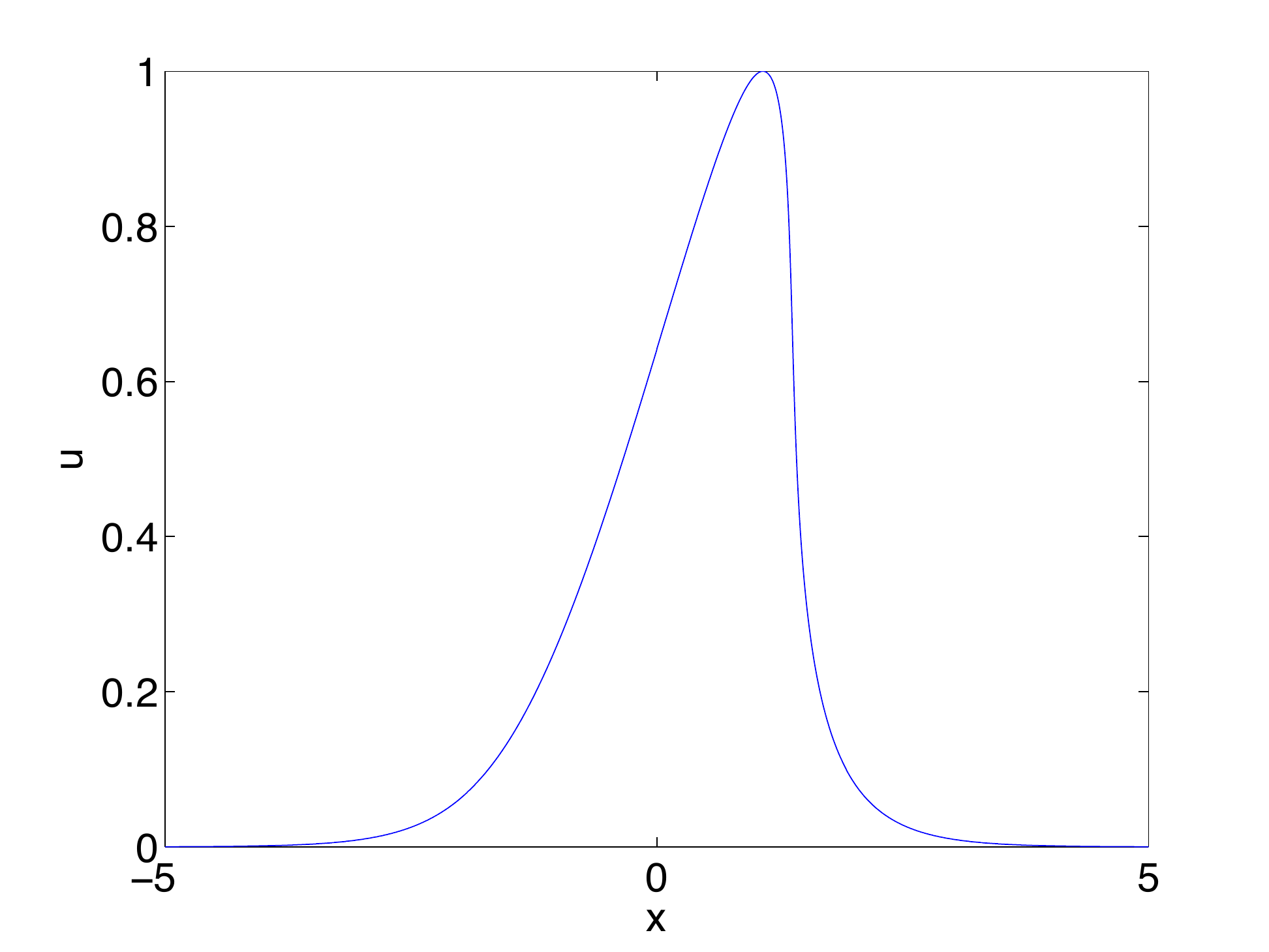}
  \includegraphics[width=0.5\textwidth]{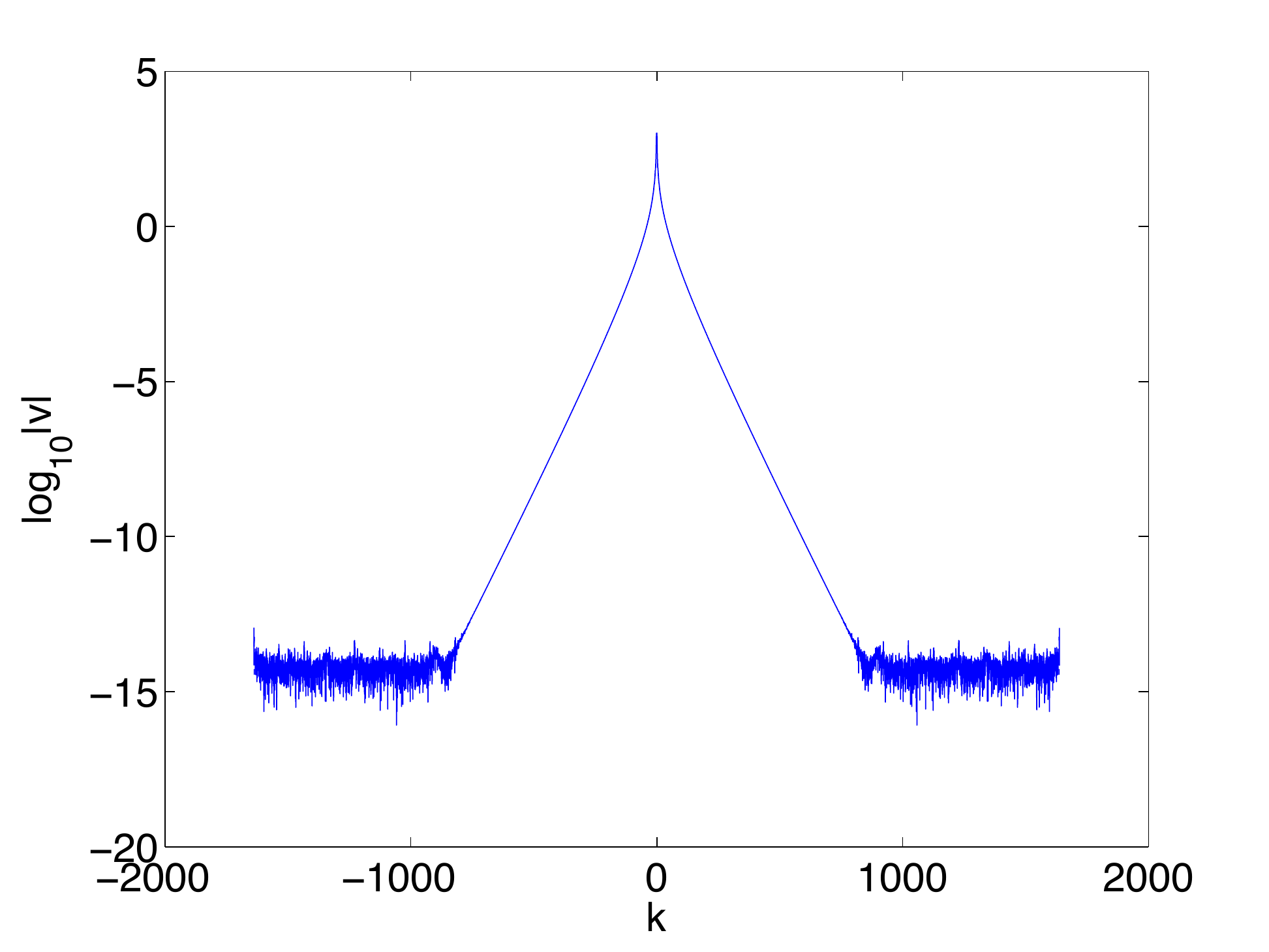}
 \caption{Exact solution of the Hopf equation for the initial data  
 $u_{0}(x)=\mbox{sech}^{2}x$ for $t=0.18\ll t_{c}$ on the left, and 
 the modulus of the 
 Fourier coefficients on the right.}
   \label{hopfex18}
\end{figure}

As in the previous subsection, we discuss the Hopf equation (\ref{Hopf}) 
for the initial data $u_{0}(x)=\mbox{sech}^{2}x$. For a time $t\ll 
t_{c}$, the solution and the Fourier coefficients computed via 
FFT with $N=2^{14}$ Fourier modes on the interval $x\in[-\pi,\pi]L$ 
with $L=5$ can be seen in 
Fig.~\ref{hopfex18}. The latter show the expected exponential 
decrease up to a level of roughly $10^{-13}$ where the error 
saturates because of the finite numerical precision.

\begin{figure}[htb!]
  \includegraphics[width=0.5\textwidth]{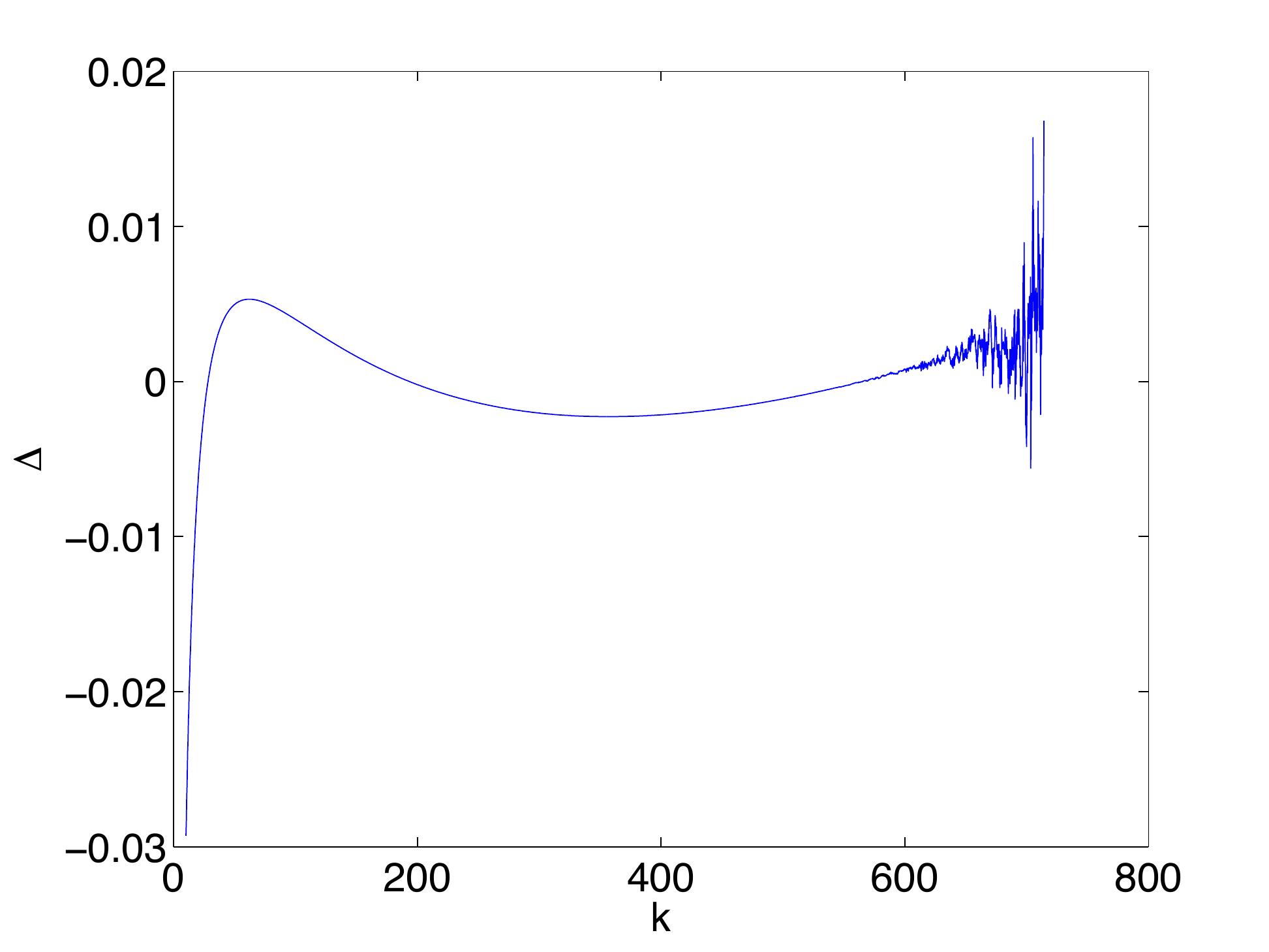}
  \includegraphics[width=0.5\textwidth]{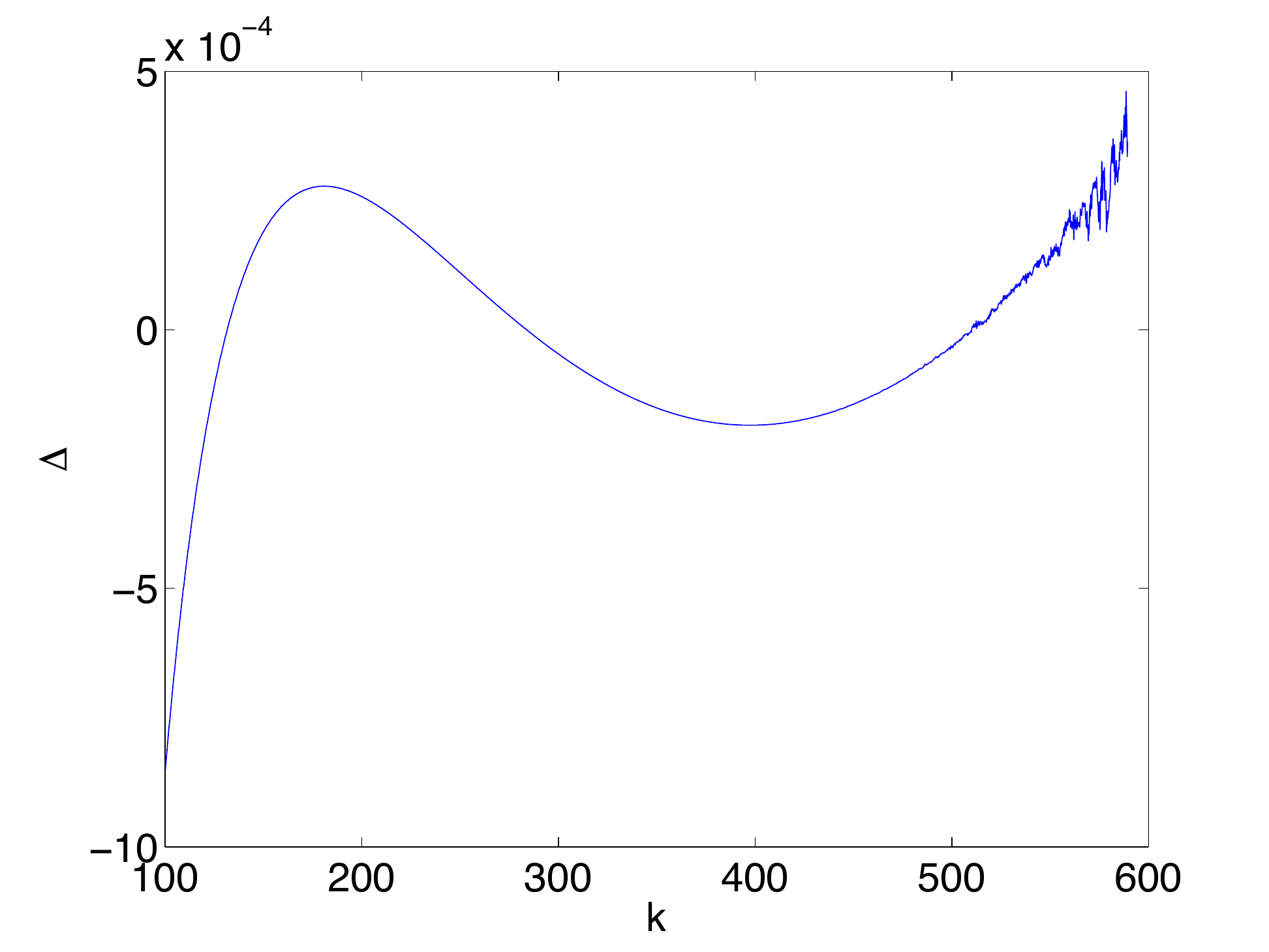}
 \caption{Difference between the modulus of the 
 Fourier coefficients for the situation of Fig.~\ref{hopfex18}
 and the fitted asymptotic formula 
 (\ref{fourierasym}) for $k>10$ and $|v|>10^{-12}$ on the left, and 
 for $k>100$ and $|v|>10^{-10}$ 
 on the right.}
   \label{hopfex18fit}
\end{figure}

In this case, only one singularity is expected in the lower half 
plane. To test the asymptotic formula (\ref{fourierasym}) in practice, we first do a least 
square fitting for $\ln |v|$, which should be according to 
(\ref{fourierasym}) of the form
\begin{equation}
    \ln |v|\sim A- B\ln k-k\delta.
    \label{abd}
\end{equation}
The fitting is done for a given range of wave numbers $k$ (we only 
consider positive $k$). For obvious reasons, $k$ has to be limited to values 
for which $|v|$ is larger than 
the rounding errors. We choose this threshold to be $|v|>10^{-12}$. In 
principle formula (\ref{fourierasym}) holds only for $k\gg 1$. However 
if we do the fitting for all positive $k$, we get $A = 6.6874$, 
$B=1.4309$ and $\delta=0.0349$ as well as $\Delta=2.05$, where  
we define the quantity $\Delta$ 
 as the 
 $L_{\infty}$ norm of the difference between the solution and the 
the asymptotic formula (\ref{abd}),
\begin{equation}
    \Delta:= \|   \ln |v| - (A - B \ln k - k \delta)   \|_{\infty}.
    \label{Delta}
\end{equation}
As expected the difference is 
maximal near $k=0$. For $k_{min}=10$ we get $A = 6.68693$, 
$B=1.4731$ and $\delta=0.0348$ and $0.0293$ for the quantity $\Delta$.
The difference between $\ln |v|$ and 
the fitted formula is shown in Fig.~\ref{hopfex18fit}. The main 
difference appears for small $k$. From the figure it can be concluded 
that a minimal value of $k_{min}=100$ might be optimal. On the other 
hand it can be also seen that rounding errors lead to some fuzzy 
behavior of the Fourier coefficients at high wave numbers. This 
suggests to consider only values $|v|>10^{-10}$. With these choices, 
we obtain $A = 6.9311$, 
$B=1.4863$ and $\delta=0.0347$. The difference between $\ln |v|$ and 
the fitted formula can be seen in Fig.~\ref{hopfex18fit}, it is of 
the order $10^{-4}$. Note that the value of $\delta$ is not much 
affected by the different choices of these thresholds. This just 
reflects the fact that the fitting to an exponential is much more stable than 
the fitting to an algebraic function.

The $L_{\infty}$ norm $\Delta$ (\ref{Delta}) of the difference between $\ln |v|$ and 
the fitted formula (\ref{fourierasym}) can be used to indicate the 
quality of the fitting. An additional test is provided from the 
relation between $A$ and $B$.  Since we compute the Fourier 
coefficients via a discrete Fourier transform whereas formula (\ref{fourierasym}) holds for 
a standard continuous Fourier transform, an additional term $\ln 
\frac{2\pi L}{N}$ appears in the formula, $A=\frac{1}{2}\ln 
(2\pi)+(B-\frac{1}{2})\ln(B-1)-(B-1)+\ln \frac{2\pi L}{N}$. We get for 
the right hand side of this relation $6.5352$, a value close to the 
computed $A = 6.9311$. However it is clear that the quantity $\Delta$ 
is the more convenient indicator of the quality of the fitting. Thus 
it will be used to this end in the following.

\begin{figure}[htb!]
\begin{center}
      \includegraphics[width=0.6\textwidth]{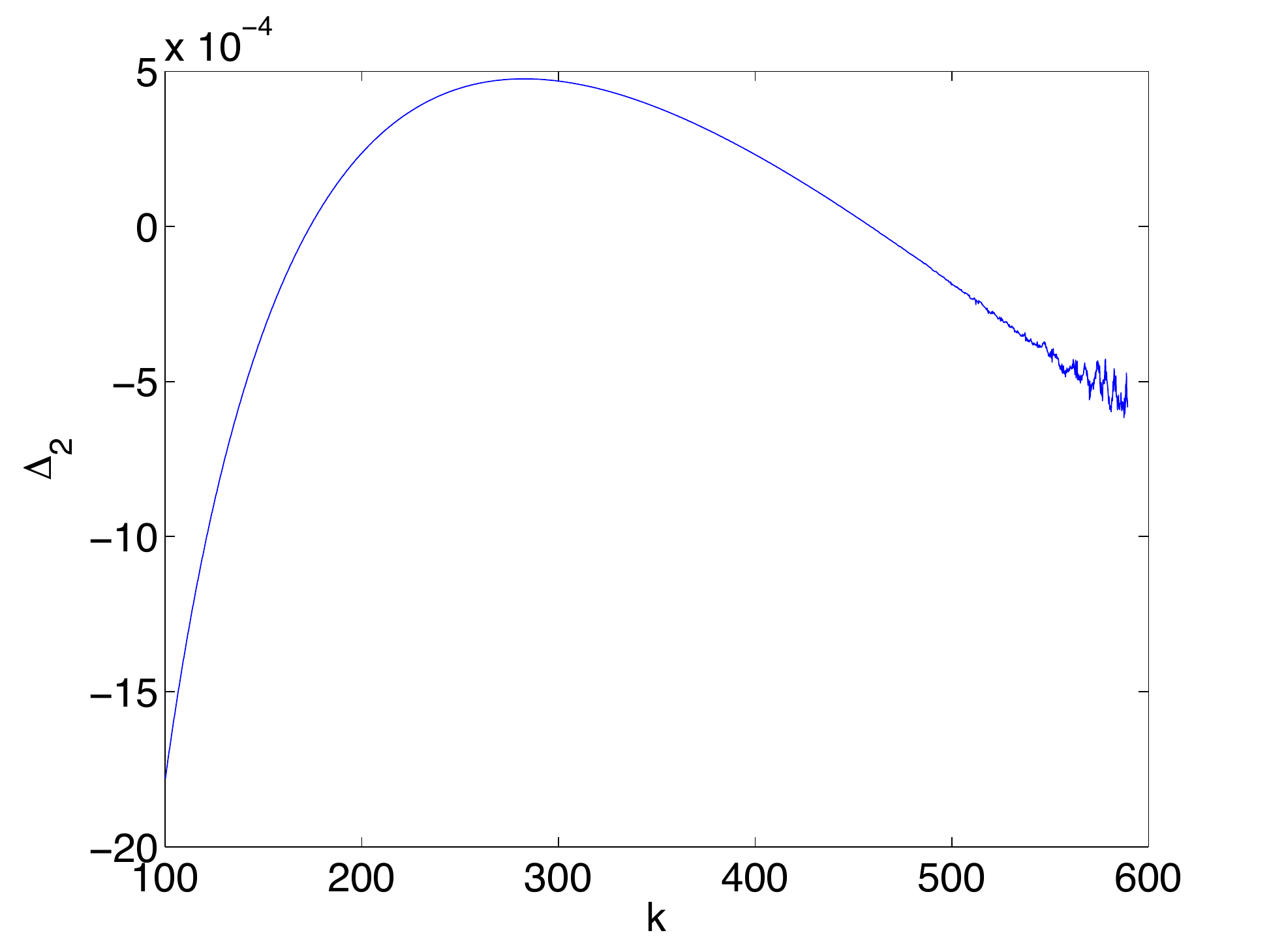}
     \caption{Difference $\Delta_{2}$ between the imaginary part  of the 
     logarithm of the
     Fourier coefficients for the situation of Fig.~\ref{hopfex18}
     and the fitted asymptotic formula 
     (\ref{fourierasym})
     for $k>100$ and $|v|>10^{-10}$.}
       \label{hopfex18fita}
\end{center}
\end{figure}

The above fitting has only been carried out for the absolute value of 
the Fourier coefficients. To determine $\alpha$, the real part of the 
location of the singularity, we fit the imaginary part of the 
logarithm of $v$, 
\begin{equation}
    \phi:=\Im \ln v\sim C-\alpha k
    \label{phi}.
\end{equation}
Since the logarithm is 
branched in Matlab at the negative real axis with jumps of $2\pi$, the computed $\phi$ 
will in general have many jumps. Thus one has first to construct a 
continuous function from the computed $\phi$. The analytic 
continuation will be done in 
the following way: starting from the first value (largest $k$), we 
check for all other values of $\phi(k_{j})$ whether 
$|\phi(k_{j+1})-\phi(k_{j})|>|\phi(k_{j+1})-\phi(k_{j})\pm\pi|$. If 
this is the case, we put $\phi(k_{j+1})\to \phi(k_{j+1})\pm \pi$.  
The result of this procedure will be a continuous function which will 
be fitted with a least square approach to a linear function. For the 
above example, we get $\alpha=1.3754$ and $C=136.5425$. The 
$L_{\infty}$ norm of the
difference of the function $\phi$ and the fitted 
straight line,    
\begin{equation*}
    \Delta_{2}=||\Im \ln v- C+\alpha k||_{\infty},
\end{equation*}
can be seen in Fig.~\ref{hopfex18fita}. It is of the 
same order of magnitude as the quantity $\Delta$ which shows the self 
consistency of the approach. 

\begin{figure}[htb!]
\begin{center}
      \includegraphics[width=0.6\textwidth]{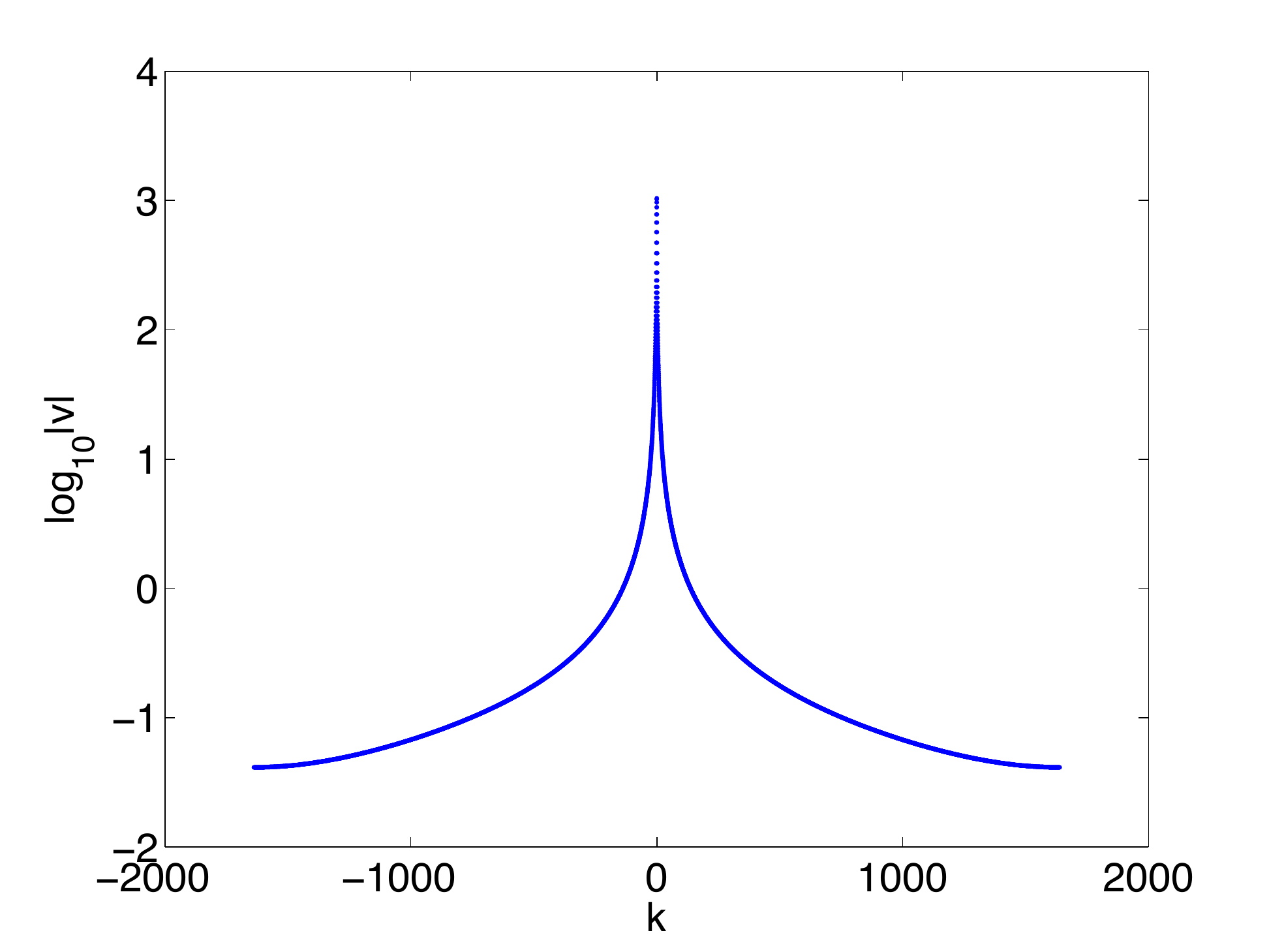}
     \caption{Modulus of the 
     Fourier coefficients of the exact solution of the Hopf equation for the initial data  
     $u_{0}(x)=\mbox{sech}^{2}x$ at the critical time $t_{c}$.}
       \label{hopfexcfourier}
\end{center}
\end{figure}

The above situation is typical for a single singularity in the 
complex plane away from the real axis. This picture changes if 
$\delta=0$, i.e., when the singularity hits the real axis which 
happens for the considered example at the critical time $t_{c}$. In 
this case the Fourier coefficients as expected no longer show an 
exponential decrease as can be seen in Fig.~\ref{hopfexcfourier}. 
Note the difference to the numerically computed Fourier coefficients 
at high wave numbers in Fig.~\ref{hopfc14}.

\begin{figure}[htb!]
  \includegraphics[width=0.5\textwidth]{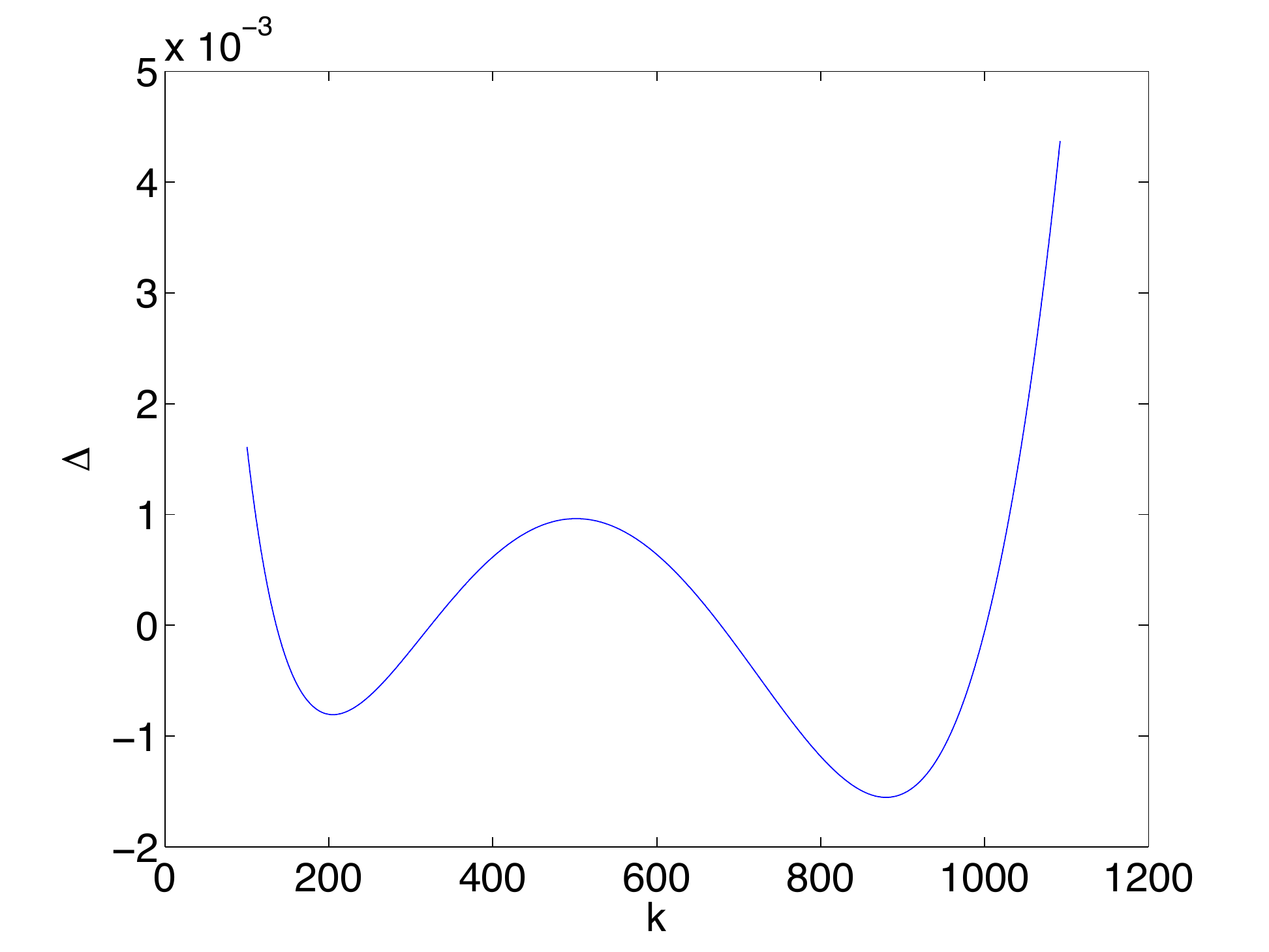}
  \includegraphics[width=0.5\textwidth]{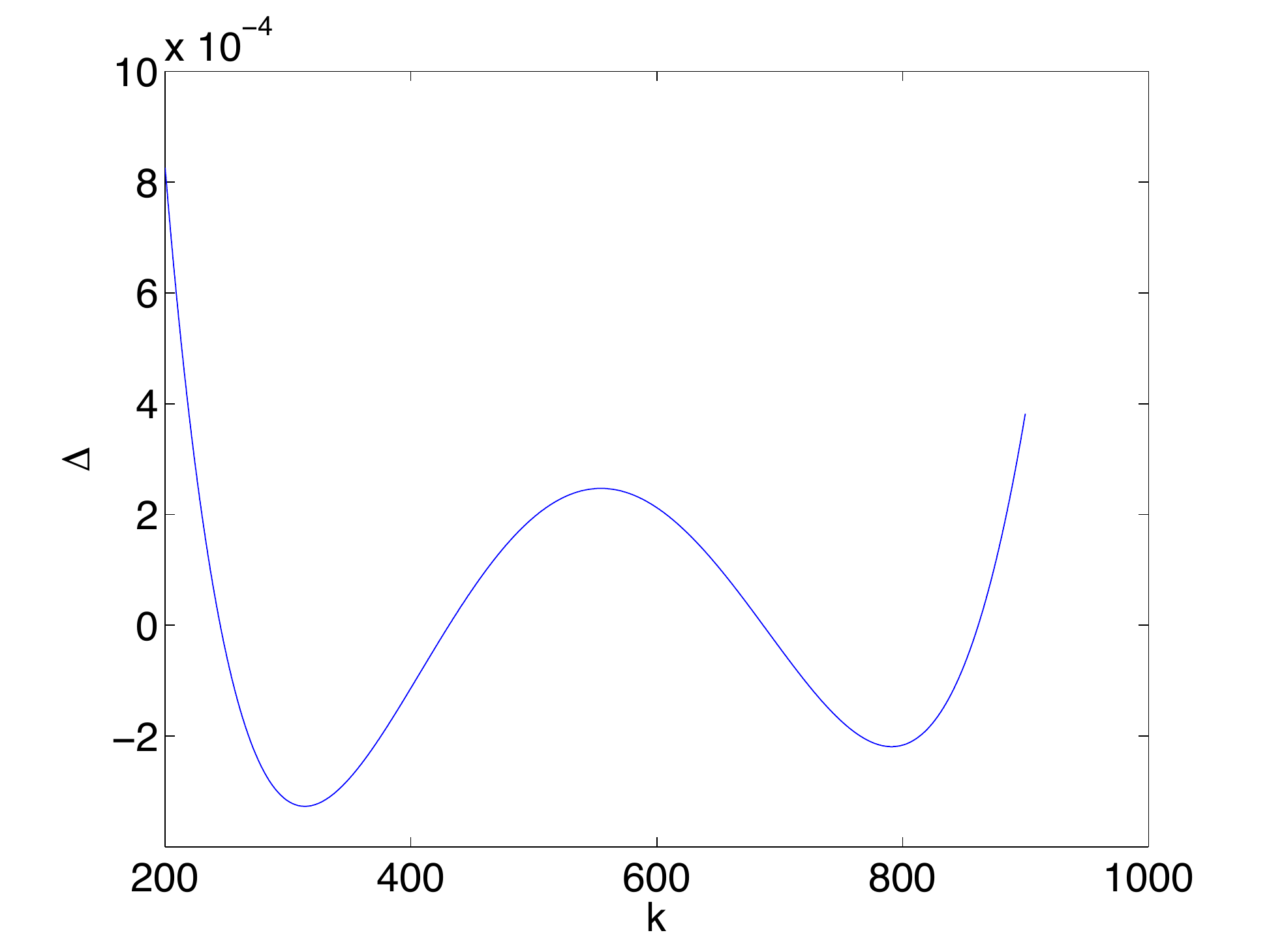}
 \caption{Difference between the modulus of the 
 Fourier coefficients for the situation of Fig.~\ref{hopfexcfourier}
 and the fitted asymptotic formula 
 (\ref{fourierasym}) for $k>100$ on the left, and 
 for $200<k<900$  
 on the right.}
   \label{hopfexcfit}
\end{figure}

The asymptotic fitting of $\ln |v|$ as in Fig.~\ref{hopfex18fit} for 
$k>100$ yields $A = 7.0283$, 
$B=1.4249$ and $\delta=-0.0001$. The difference between $\ln |v|$ and 
the fitting curve can be seen in Fig.~\ref{hopfexcfit}. Since there 
is no exponential decay in this case, the fitting is much more 
sensitive to the considered cutoffs, for low $k$ due to the 
condition $k\gg1$ for the asymptotic formula to be valid, for large $k$ due to rounding 
errors. This is also obvious from Fig.~\ref{hopfexcfit} where the 
errors are maximal at both ends of the shown spectrum.  For 
$200<k<900$ we get $A = 6.4573$, 
$B=1.3084$ and $\delta=0.0001$ with a fitting error as shown in 
Fig.~\ref{hopfexcfit}. In both cases $\delta$ is of the order of 
$10^{-4}$. The quantity $B$ should be equal to $4/3$ and is as before 
much more sensitive to a choice of the limits $k_{min}$ and $k_{max}$ 
for the fitting. 

By fitting the quantity $\phi$ for $200<k<900$, we get $\alpha=1.5251$ and 
$C=306.3523$ with a fitting error of the same order. 
The former is very close to the exact value $x_{c} = 
\sqrt{3}/2 - \ln((\sqrt{3}-1)/\sqrt{2})\sim 1.5245$.
If we repeat the above fitting for $N=2^{15}$ Fourier modes for 
$200<k<N/3/L$ ($L=5$ as before), we obtain $A = 7.4012$, 
$B=1.3597$ and $\delta=-3.9054*10^{-5}$ with a $\Delta=0.012$ and 
$\alpha=1.5248$ and $C=306.3584$. Thus a higher number of Fourier 
modes makes the fitting more stable and reliable.

To sum up, we have shown in this subsection that the quality of the 
fitting depends on the choice of the minimal and maximal values for 
the wave numbers. Generally a choice of $k_{min}=100$ or higher 
appears to be recommended. However this has the disadvantage that 
many Fourier modes are excluded, especially at times $t\ll t_{c}$ 
where only few modes contribute. Rounding errors enforce a choice of 
$k_{max}$ of the order of $N/3/L$. By studying the difference between 
the considered quantity and the fitting curve, these values can be 
optimized as outlined in the next subsection. A proper choice of these thresholds is more important at 
the critical time when the singularity hits the real axis since there 
is no exponential decay in this case.  

\subsection{Asymptotic behaviour of the Fourier coefficients for a numerical solution to the Hopf equation}

In the previous subsection we have studied the
 asymptotic fitting of the Fourier coefficients on a finite interval of the wave numbers and 
 have fixed the limits of this interval essentially by hand. Here we will 
 do this in a more systematic way by choosing a suitable lower 
 threshold (which depends on the studied problem), and by varying the 
 upper limit to reach a prescribed difference between the Fourier 
 coefficients and the fitted curve with a maximal number of coefficients. 
 Since this difference would 
 obviously be minimal if upper and lower limit coincide, we impose the 
 additional condition that at least half of the Fourier coefficients 
 with $|v|>10^{-10}$ will be included in the fitting to include most of the 
available information on the solution. With this 
 approach and the preliminary studies of the previous sections, we 
 are then able to perform the asymptotic fitting during the actual 
 numerical computation of the Hopf solution. We show that the vanishing 
 of the computed $\delta$ 
can be in fact used to determine the critical time $t_{c}$ and the 
type of the singularity via the exponent $\mu$.  The goal is to obtain an 
error in the fitting of the same order as the numerical error.

The fitting method used in this paper is obviously not the 
only possible one. For instance in \cite{CR} and related publications 
the \emph{sliding fit} approach is applied  where only few values of $k$ 
are taken into account in the procedure (typically between $4$ and 
$100$  for a resolution of  1024 modes, see for example \cite{RLSS}). 
The quality of the fitting in this approach is roughly determined by 
the independence of the computed fitting parameters on the minimal 
value $k_{min}$ of the $k$ used in the fitting, and on the 
discretization size $2\pi L/N$. Note that the sliding procedure 
confirms the results found here, for $k_{min} \leq 100$  and up to 
$600$ points considered. Nevertheless, this approach is 
computationally expensive which implies that the time of the 
singularity formation has to be previously roughly determined in 
another way to be able to perform the study close to it. Since it is 
the aim of this paper to identify the break-up numerically via the asymptotic behavior of the Fourier coefficients, even in cases where the critical time is not known, we do not use the sliding fit here.

As in the previous section, the asymptotic fitting of $\ln(|v|)$ is done on the interval $k_{min}<k<k_{max}$, where we put
\begin{eqnarray*}
k_{min}=10 \,\,\mbox{and}\,\,\,
k_{max}=\min(k_1, k_2), \,\,\,\mbox{with}\,\, \\
k_1=\underset{|v|>1e-10}{\max} (k), \,\,\,
k_2=\underset{ \Delta < p}{\max} (k)
\end{eqnarray*}
where $p$ denotes a prescribed value for $\Delta$ in (\ref{Delta})
on the fit interval.
The choice $k_{min}=10$ allows the fitting even for short 
times where there are only few Fourier coefficients with 
$|v|>10^{-10}$. 
It satisfies the condition $k_{min}\gg 1$ whilst maximizing the number 
of the numerically computed Fourier coefficients included in the 
fitting procedure. 
Indeed, for given $p$, larger values of $k_{min}$ would lead to a 
too small number of coefficients used for the fitting.

\begin{figure}[htb!]
  \includegraphics[width=0.45\textwidth]{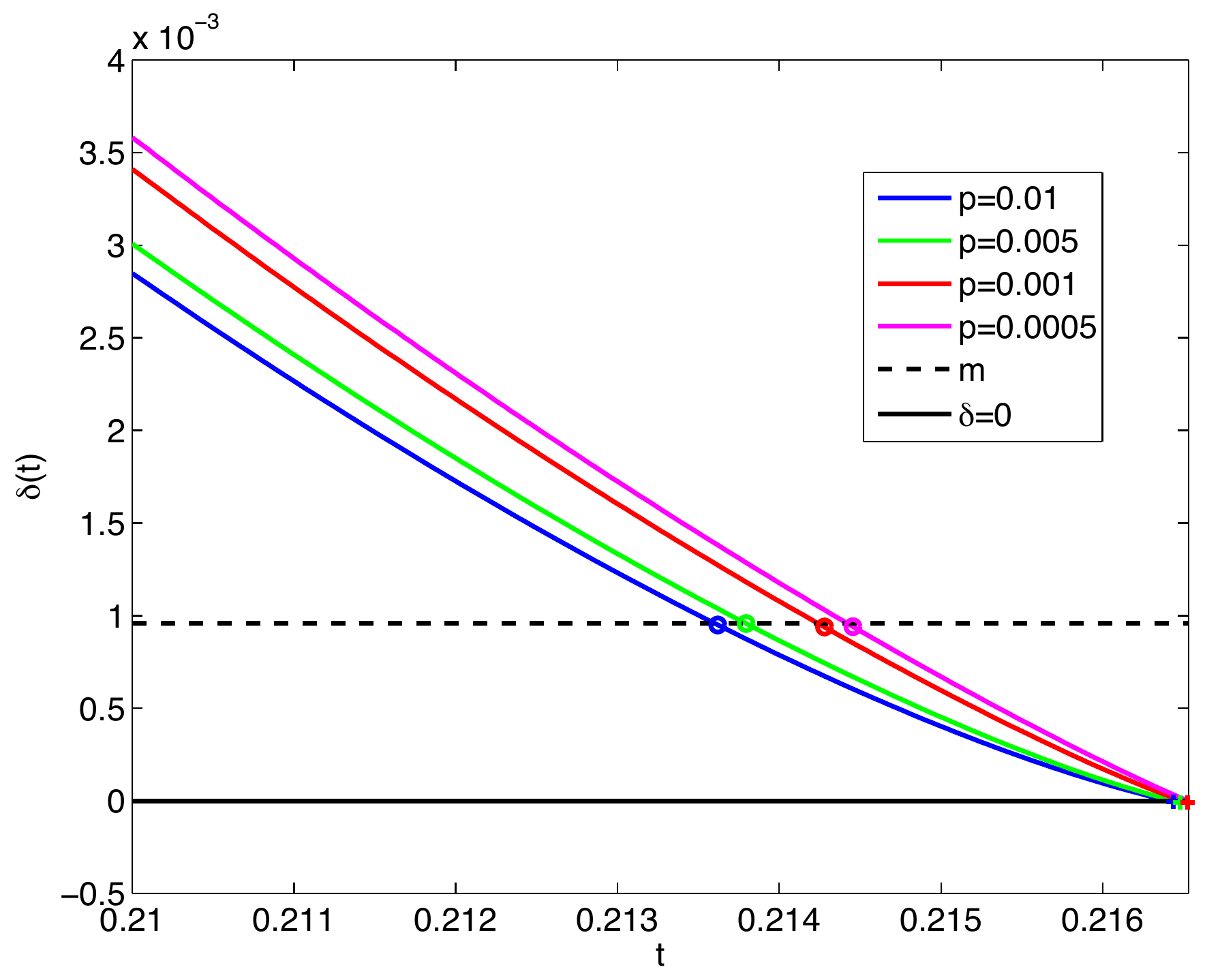}   
  \includegraphics[width=0.45\textwidth]{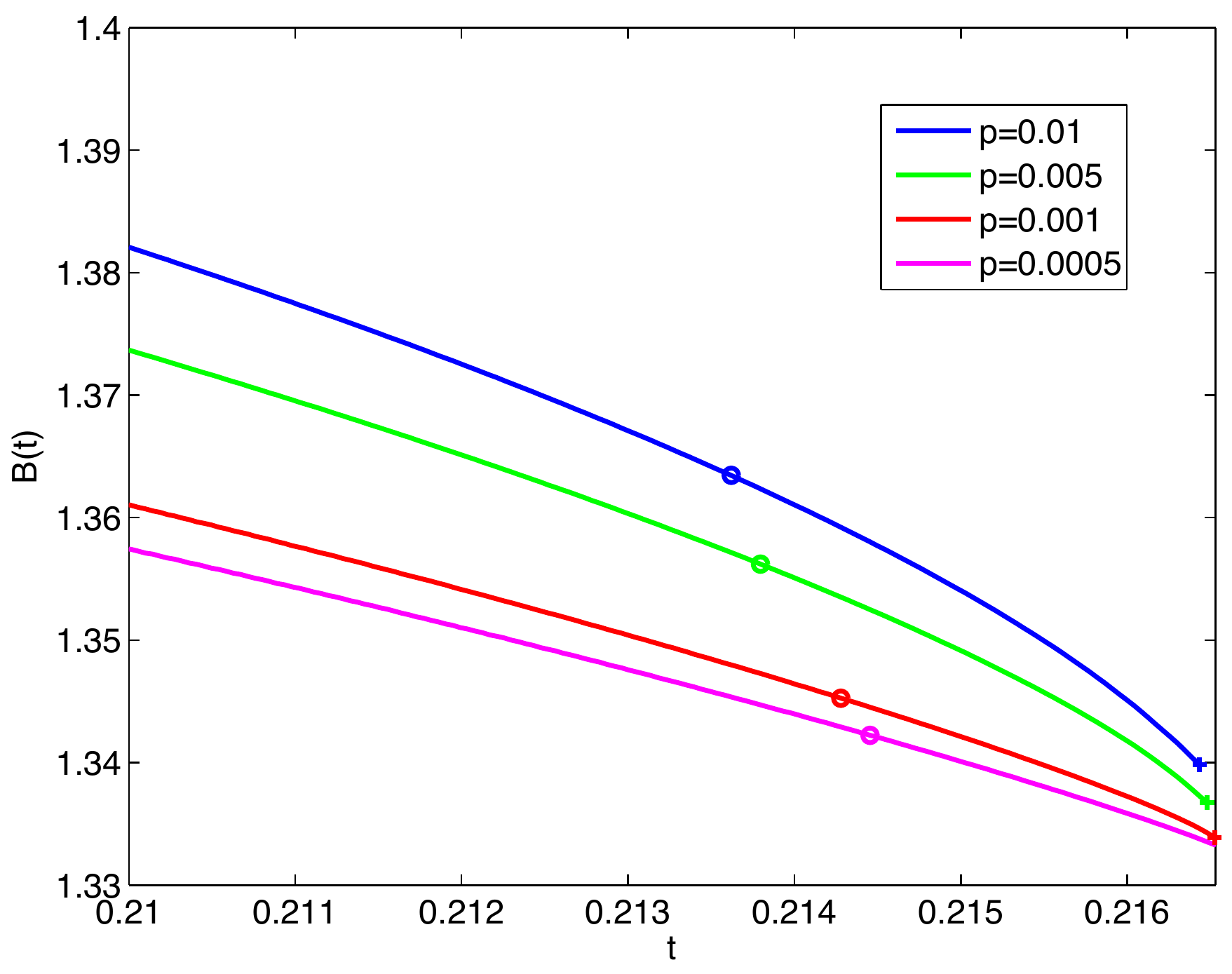}
 \caption{Time evolution of the fitting parameters $\delta$ and $B$ for the numerical solution to the Hopf equation 
 with $u(x,0)=sech^2(x)$, with $N=2^{15}$ and different values of $p$ 
 between $t=0.21$ and $t_c=\sqrt{3}/8$. The dashed line in the left 
 figure corresponds 
 to the minimal resolution in Fourier space discussed in 
 remark~\ref{ms}.}
   \label{HRKdap}
\end{figure}

For $N=2^{15}$, we observe the influence of $p$ on the time of vanishing of $\delta$ and the corresponding value of $B$ in Fig. \ref{HRKdap} and in Table \ref{tab1}. 
Due to the large computational domain considered ($x \in [-\pi, 
\pi]L$, $L=5$), a sufficient resolution has to be used to get a 
sufficiently small value of $m$ ($m:=2\pi L/N$, $m\sim 
10^{-3}$ for $N=2^{15}$ and $L=5$). 
\begin{table}[htb!]
\centering
\begin{tabular}{|c|c|c|c|c|}
\hline
$\,\,\,\,\,\,p=10^{-2}\,\,\,$ & $t$ & $\delta$ & $B$ & $\,\,\,\alpha\,\,\,\,\,\,$ \\
\hline 
$ t_{\delta<m}$ &  $ 0.2136 $ & $ 9* 10^{-4}  $ & $ 1.363  $ & $ 1.5128 $  \\
$ t_{\delta<0}$ &  $ 0.2164 $ & $ -5* 10^{-6}  $ & $ 1.339  $ & $ 1.5242 $    \\
\hline
\end{tabular}
\\
\vspace{0.3cm}
\begin{tabular}{|c|c|c|c|c|}
\hline
$p=5* 10^{-3}$ & $t$ & $\delta$ & $B$ & $\alpha$\\
 \hline
$ t_{\delta<m}$ &  $ 0.2138 $ & $ 9* 10^{-4}  $ & $ 1.356  $ & $ 1.5134$   \\
$ t_{\delta<0}$ & $ 0.2165 $ & $ -7* 10^{-6}  $ & $ 1.337  $ & $ 1.5244 $   \\
\hline
\end{tabular}
\\
\vspace{0.3cm}
\begin{tabular}{|c|c|c|c|c|}
\hline
$\,\,\,\,\,\,p= 10^{-3}\,\,\,$ & $t$ & $\delta$ & $B$ & $\,\,\,\alpha\,\,\,\,\,\,$  \\
\hline
$ t_{\delta<m}$ &  $ 0.2143 $ & $ 9* 10^{-4}  $ & $ 1.345  $ & $ 1.5149 $    \\
$ t_{\delta<0}$ & $ 0.2165 $ & $ -7* 10^{-6}  $ & $ 1.334  $ & $ 1.5245 $  \\
\hline
\end{tabular}
\\
\vspace{0.3cm}
\begin{tabular}{|c|c|c|c|c|}
\hline
$p=5* 10^{-4}$ & $t$ & $\delta$ & $B$ & $\alpha$\\
\hline
$ t_{\delta<m}$ &  $ 0.2145 $ & $ 9* 10^{-4}  $ & $ 1.342  $ & $ 1.5153 $   \\
$ t_{\delta<0}$ & $ 0.2166 $ & $ -1* 10^{-5}  $ & $ 1.333  $ & $ 1.5243 $   \\
\hline
\end{tabular}
\caption{Values of the fitting parameters for the numerical solution to the Hopf equation 
 with $u(x,0)=sech(x)^2$, with $N=2^{15}$ and different values of 
 $p=0.01, 0.005, 0.001, 0.0005$. The exact values are 
 $t_{c}\sim0.2165$, $B=4/3$ and $x_{c}\sim1.5245$.}
\label{tab1}
\end{table}

If $p$ decreases to $10^{-3}$, one finds that $\left(t_c, B(t_c)\right) \to 
\left(\frac{\sqrt{3}}{8},\, \frac{4}{3} \right)$. But for smaller values of $p$, for example for
$p=5*10^{-4}$, the values for both $t_c$ and  $\alpha(t_c)$ become 
worse approximations of the exact ones.
In fact, to get such small fitting errors, we have to consider too few 
Fourier coefficients which explains the lack of precision for these. Thus we 
add the requirement that at least half of the total number of the Fourier coefficients available have to be used for the fitting. 
For $N=2^{15}$, this implies that the `minimal' achievable fitting error whilst using enough Fourier coefficients is $p=0.01$, which corresponds roughly to 
the precision with which we compute the numerical solution, see 
Sec.~\ref{numhopf}.

\begin{figure}[htb!]
  \includegraphics[width=0.45\textwidth]{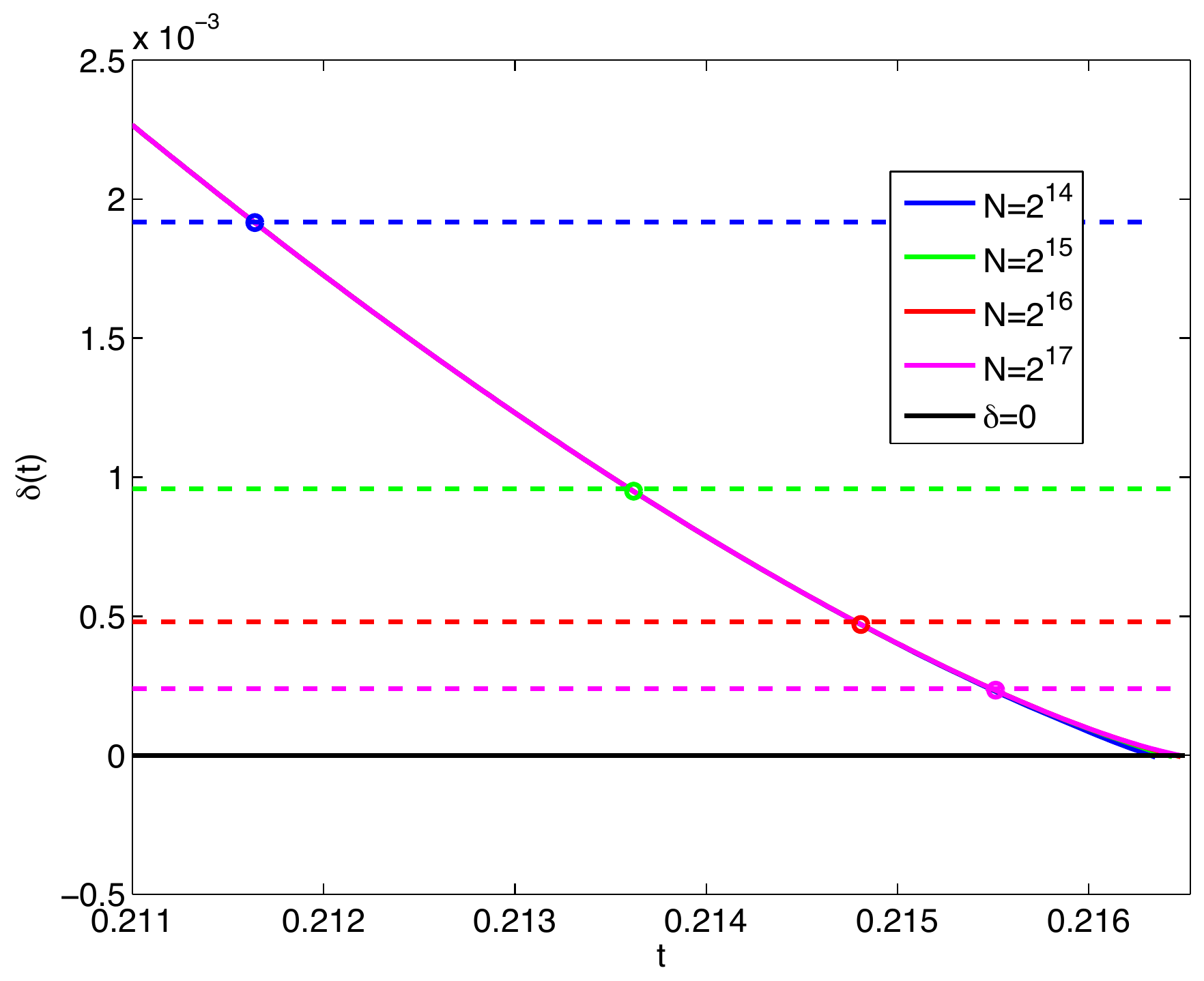}\includegraphics[width=0.45\textwidth]{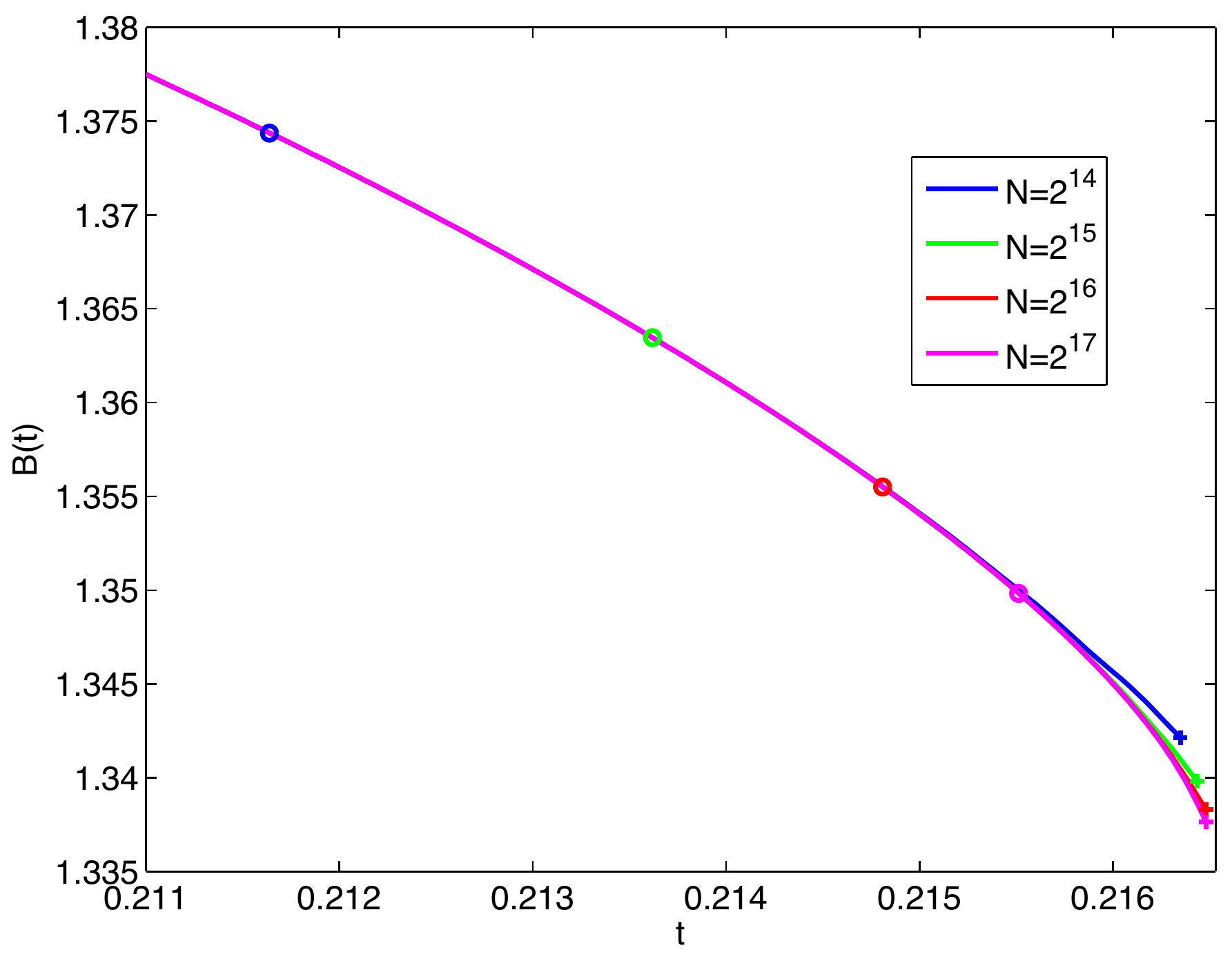} 
 \caption{Time evolution of the fitting parameters $\delta$ and $B$ for the numerical solution to the Hopf equation 
 with $u(x,0)=sech(x)^2$, for different values of $N$ between $t=0.21$ and $t_c=\sqrt{3}/8$.}
   \label{HRKp01}
\end{figure}

Thus for different values of $N$, and the related precision with which we compute the numerical solution, 
we observe in Fig. \ref{HRKp01} the time evolution of $\delta$ and $B$, and give the values of 
$\delta, B$ and $\alpha$ at two different times in Table \ref{tab2}.
\begin{table}[htb!]
\centering
\begin{tabular}{|c|c|c|c|c|}
\hline
$N=2^{14}$ & $t$ & $\delta$ & $B$ & $\alpha$ \\
 \hline
$ t_{\delta<m}$ & $ 0.2116 $ & $ 0.0019  $ & $ 1.37  $ & $ 1.5048 $    \\
$ t_{\delta<0}$  & $ 0.2164 $ & $ -6* 10^{-6}  $ & $ 1.345  $ & $ 1.5238 $  \\
\hline
\end{tabular}
\\
\vspace{0.4cm}
\begin{tabular}{|c|c|c|c|c|}
\hline
$N=2^{15}$ & $t$ & $\delta$ & $B$ & $\alpha$\\
 \hline
$ t_{\delta<m}$ &  $ 0.2136 $ & $ 9* 10^{-4}  $ & $ 1.363  $ & $ 1.5128 $   \\
$ t_{\delta<0}$ & $ 0.2164 $ & $ -5* 10^{-6}  $ & $ 1.339  $ & $ 1.5242 $   \\
\hline
\end{tabular}
\\
\vspace{0.4cm}
\begin{tabular}{|c|c|c|c|c|}
\hline
$N=2^{16}$ & $t$ & $\delta$ & $B$ & $\alpha$\\
 \hline
$ t_{\delta<m}$ &  $ 0.2148$ & $ 4.7* 10^{-4}  $ & $ 1.355  $ & $ 1.5176 $   \\
$ t_{\delta<0}$ & $ 0.2165 $ & $ -5* 10^{-6}  $ & $ 1.338  $ & $ 1.5244 $   \\
\hline
\end{tabular}
\\
\vspace{0.4cm}
\begin{tabular}{|c|c|c|c|c|}
\hline
$N=2^{17}$ & $t$ & $\delta$ & $B$ & $\alpha$\\
 \hline
$ t_{\delta<m}$ &  $ 0.2155$ & $ 2.3* 10^{-4}  $ & $ 1.349  $ & $ 1.5204 $   \\
$ t_{\delta<0}$ & $ 0.2165 $ & $ -7* 10^{-7}  $ & $ 1.337  $ & $ 1.5244 $   \\
\hline
\end{tabular}
\caption{Values of the fitting parameters for the numerical solution to the Hopf equation 
 with $u(x,0)=sech(x)^2$, for different values of $N=2^{14}, 2^{15}, 2^{16}$ and $2^{17}$.}
\label{tab2}
\end{table}
As $N$ increases, $\left(t_c, B(t_c)\right) \to 
\left(\frac{\sqrt{3}}{8},\, \frac{4}{3} \right)$ and $\alpha(t_c) \to 
1.5245$ without further restrictions on $k$. The time dependence of 
the fitting parameters is very similar for different $N$, which indicates that the fitting is reliable.\\
We conclude from this example that the  method can 
be  used in practice to determine the critical time $t_{c}$ and the 
type of the singularity via the value of $B(t_{c})$ if sufficient 
spatial resolution is used. The choice of the minimal wave number $k_{min}$ for 
the fitting, and the achievable minimal fitting error depend on the 
studied problem.

\subsection{Burgers' equation}
The Burgers' equation (\ref{burgers}) with small dissipation $0<\epsilon\ll 1$,
can be viewed as a purely dissipative regularization of the Hopf 
equation. An asymptotic description at the critical point was given 
in \cite{DE13}. Since the system of ODEs resulting from the approximation 
of the solution via a discrete Fourier series in the spatial 
coordinate is now stiff, we use a stiff integrator, a fourth order 
\emph{exponential time differencing} scheme by Cox and Matthews 
\cite{CM} as implemented in \cite{Schme,ckkdvnls}. This method is very efficient for both Burgers' and KdV 
equation, see e.g.~\cite{ckkdvnls}. 
In both cases we consider as before the initial 
data $u_{0}(x)=\mbox{sech}^{2}x$ for $x\in[-5\pi,5\pi]$.

\begin{figure}[htb!]
  \includegraphics[width=0.5\textwidth]{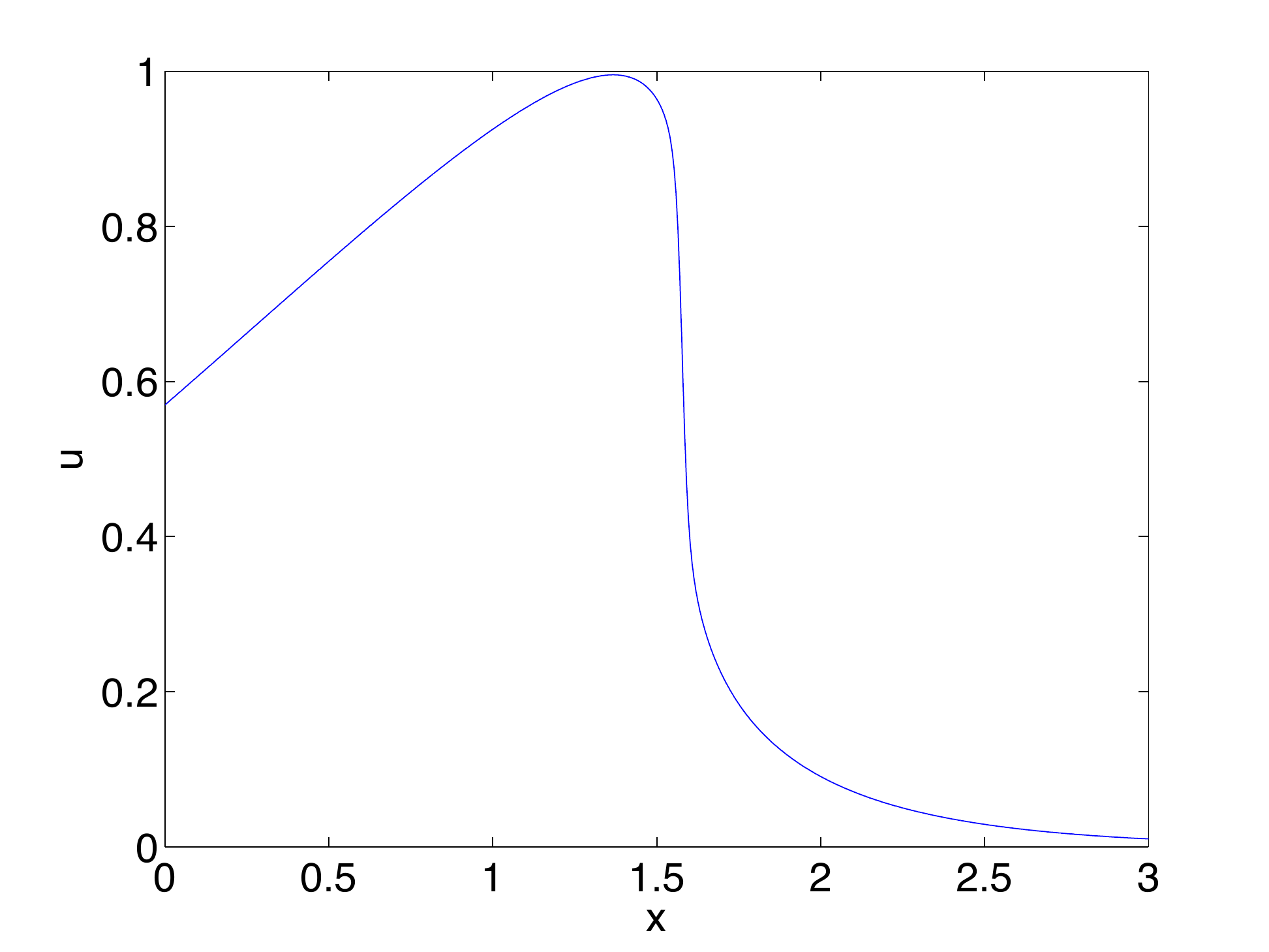}
  \includegraphics[width=0.5\textwidth]{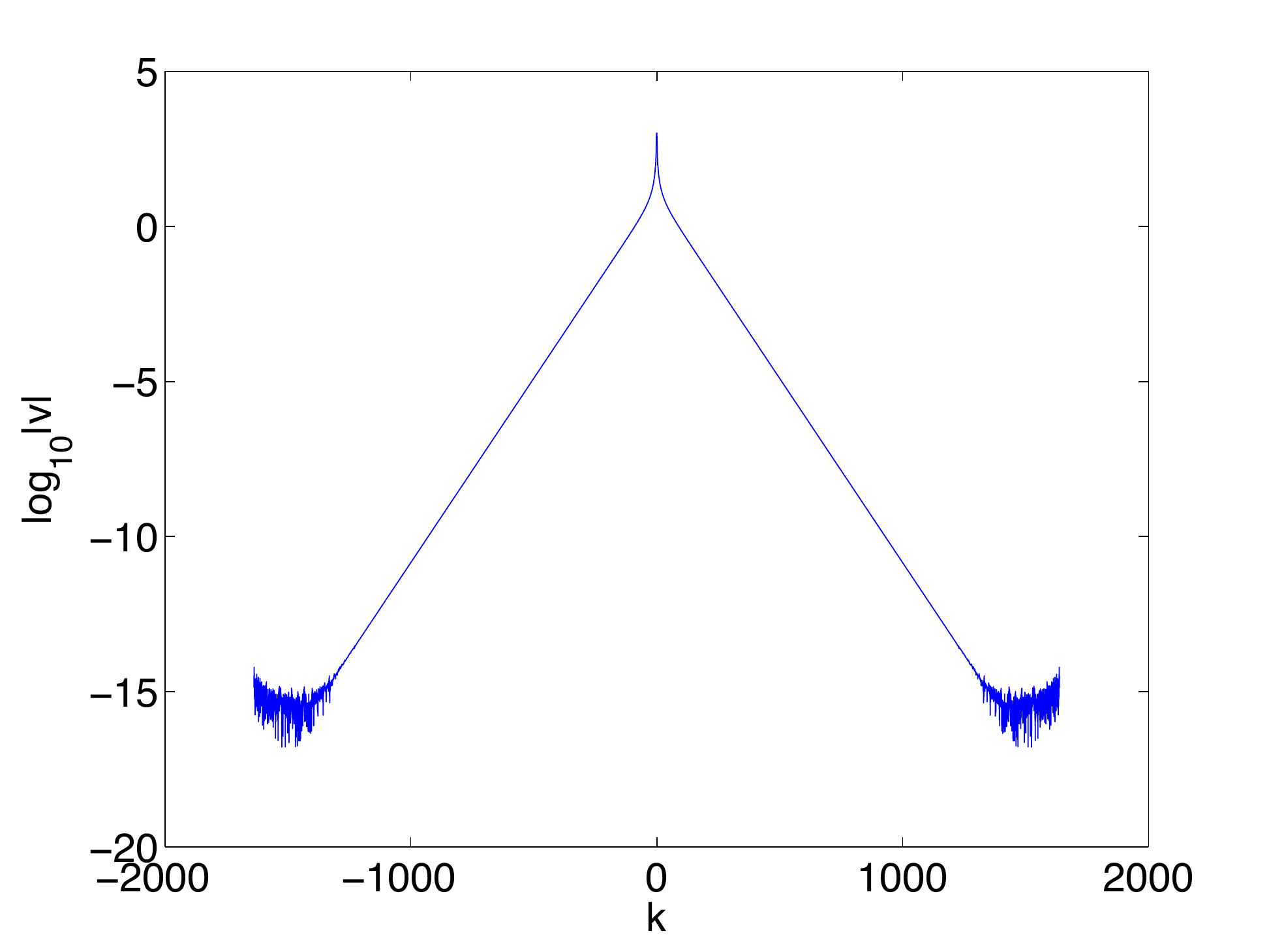}
 \caption{Solution to the Burgers' equation for the initial data 
 $u_{0}(x)=\mbox{sech}^{2}x$ for $t=0.23>t_{c}$ and $\epsilon=0.01$ 
 on the left, and 
 the modulus of the corresponding Fourier  coefficients 
 on the right.}
   \label{burgers01}
\end{figure}

For the Burgers' equation the shock of the Hopf 
solution for $t>t_{c}$ will be regularized  as can be seen for $\epsilon=0.01$ in 
Fig.~\ref{burgers01} at time 
$t=0.23>t_{c}$. There is a zone with a strong, but finite gradient of 
order $1/\epsilon$.

The Fourier coefficients for this solution are also shown in 
Fig.~\ref{burgers01}. It can be seen that they decrease (essentially 
exponentially) to machine precision. Doing a least square fitting 
for the modulus of the Fourier coefficients with 
$|v|>10^{-10}$ for $k>100$ (this is the result of the procedure 
outlined in the previous section with a fitting error of the order 
$10^{-2}$), we find  $A = 2.4858$, 
$B=0.0273$ and $\delta=0.0183$. This shows that the parameter 
$\delta$ indicating the distance between the nearest singularity to 
the real axis stays finite even for $t>t_{c}$ as expected. The 
difference between the modulus of the Fourier coefficients and the 
fitted asymptotic formula is of the order of $10^{-3}$.
Fitting the quantity $\phi$ in (\ref{phi}), we get $\alpha=1.5778$ 
and $C=155.5075$. This implies that the real part of the location of 
the singularity is shifted towards larger values than for the critical 
value $x_{c}$ of the Hopf solution which is not surprising since the 
Hopf shock also propagates with $t$. The difference between $\phi$ and the 
fitted asymptotic formula is of the order of $10^{-3}$, too. 

\begin{figure}[htb!]
  \includegraphics[width=0.5\textwidth]{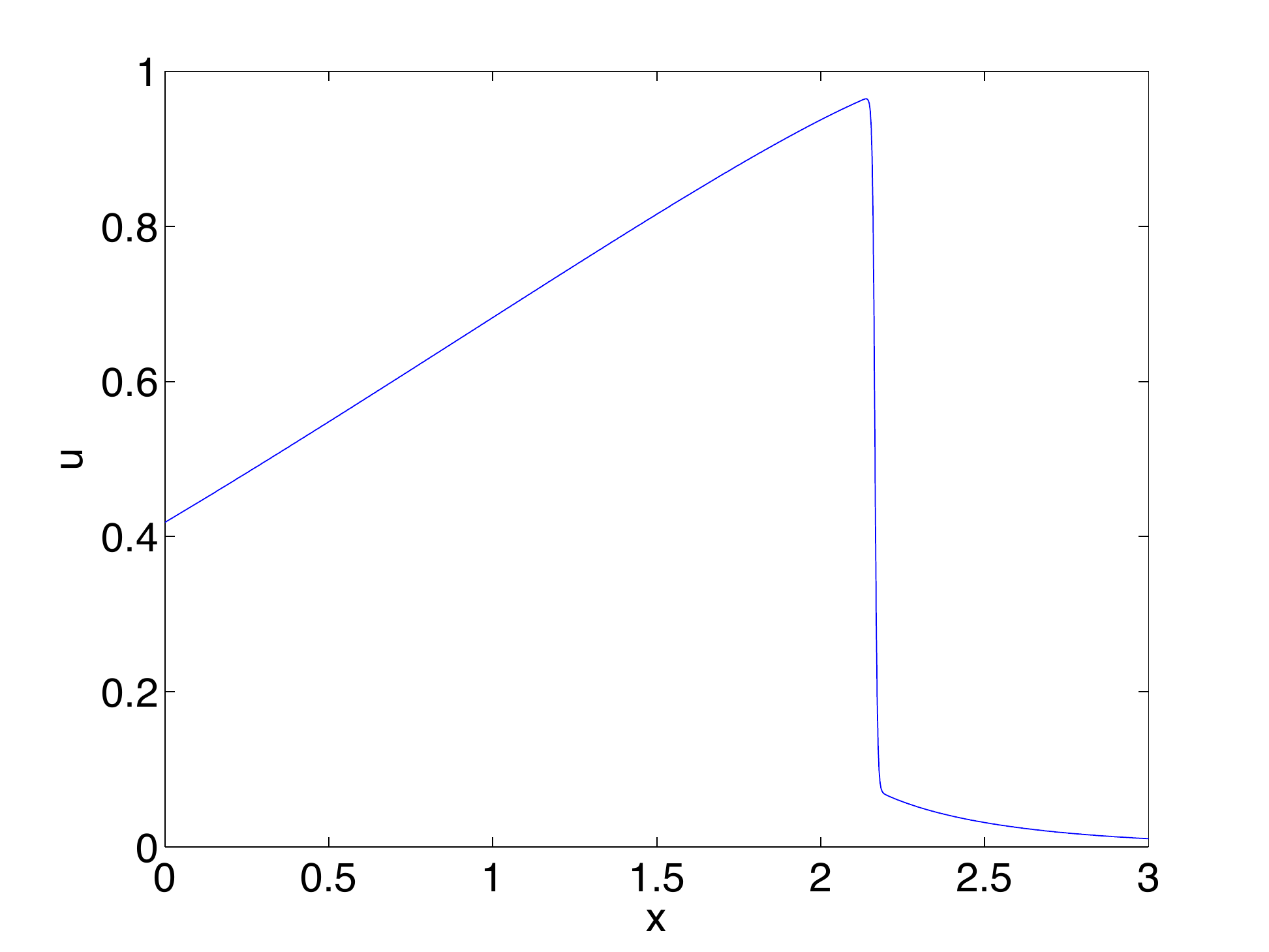}
  \includegraphics[width=0.5\textwidth]{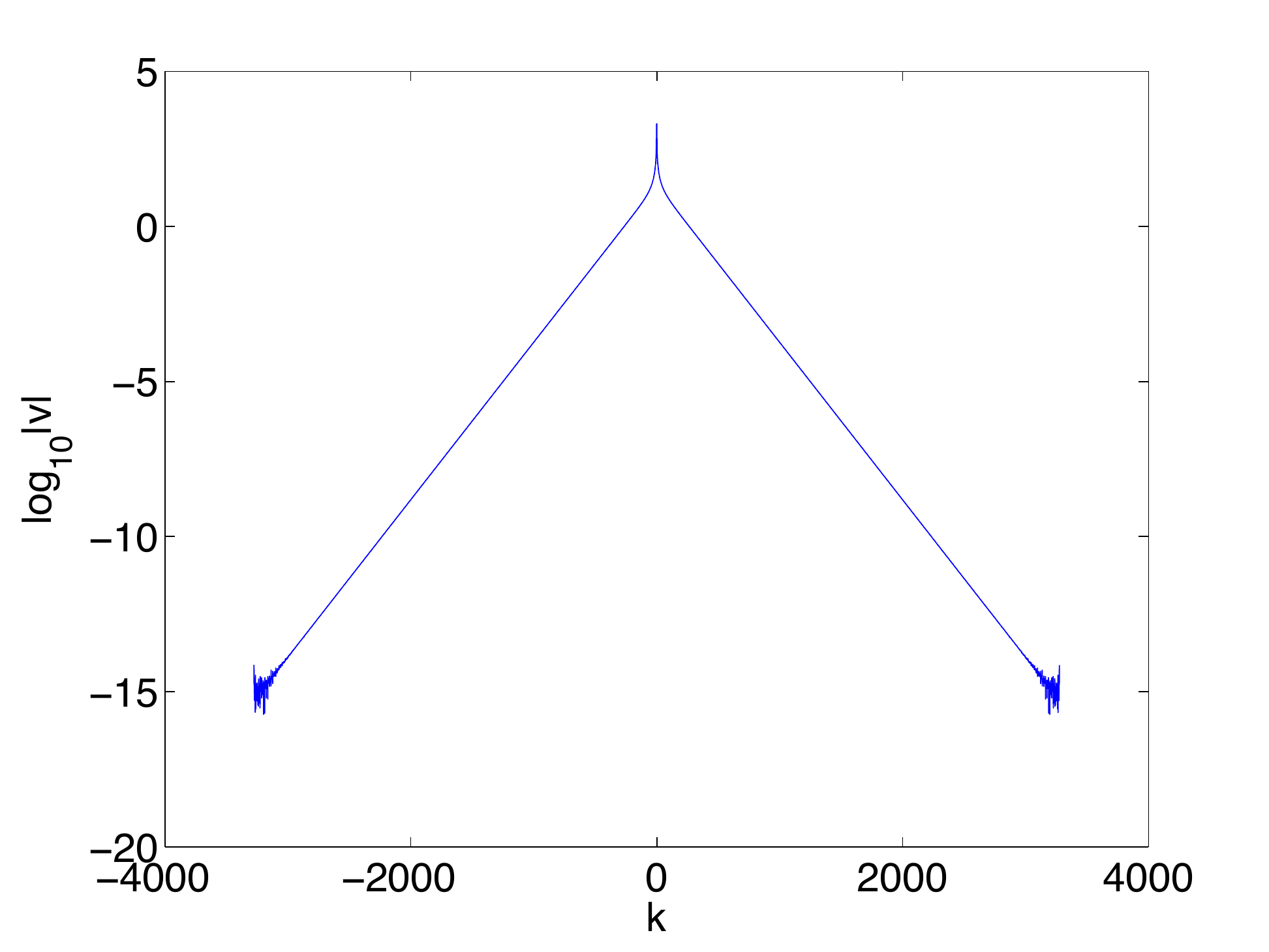}
 \caption{Solution to the Burgers' equation for the initial data 
 $u_{0}(x)=\mbox{sech}^{2}x$ for $t=0.4\gg t_{c}$ and $\epsilon=0.01$ 
 on the left, and 
 the modulus of the corresponding Fourier  coefficients 
 on the right.}
   \label{burgers04}
\end{figure}

For times $t\gg t_{c}$ the solution in the vicinity of 
the Hopf shock becomes steeper as shown in Fig.~\ref{burgers04}, which 
means that more Fourier modes have to be used in order that the 
modulus of the coefficients decreases to machine precision ($N=2^{15}$ instead of 
$N=2^{14}$ before).
A least square fitting as 
before for the modulus of the Fourier coefficients for $k>100$ 
gives  $A = 3.2397$, 
$B=0.0265$ and $\delta=0.0117$. The 
difference between the modulus of the Fourier coefficients and the 
fitted asymptotic formula is of the order of $10^{-2}$.
Fitting the quantity 
$\phi$ in (\ref{phi}), we get $\alpha=2.1657$ 
and $C=218.3404$ with a similar quality of the fitting. 

\begin{figure}[htb!]
  \includegraphics[width=0.5\textwidth]{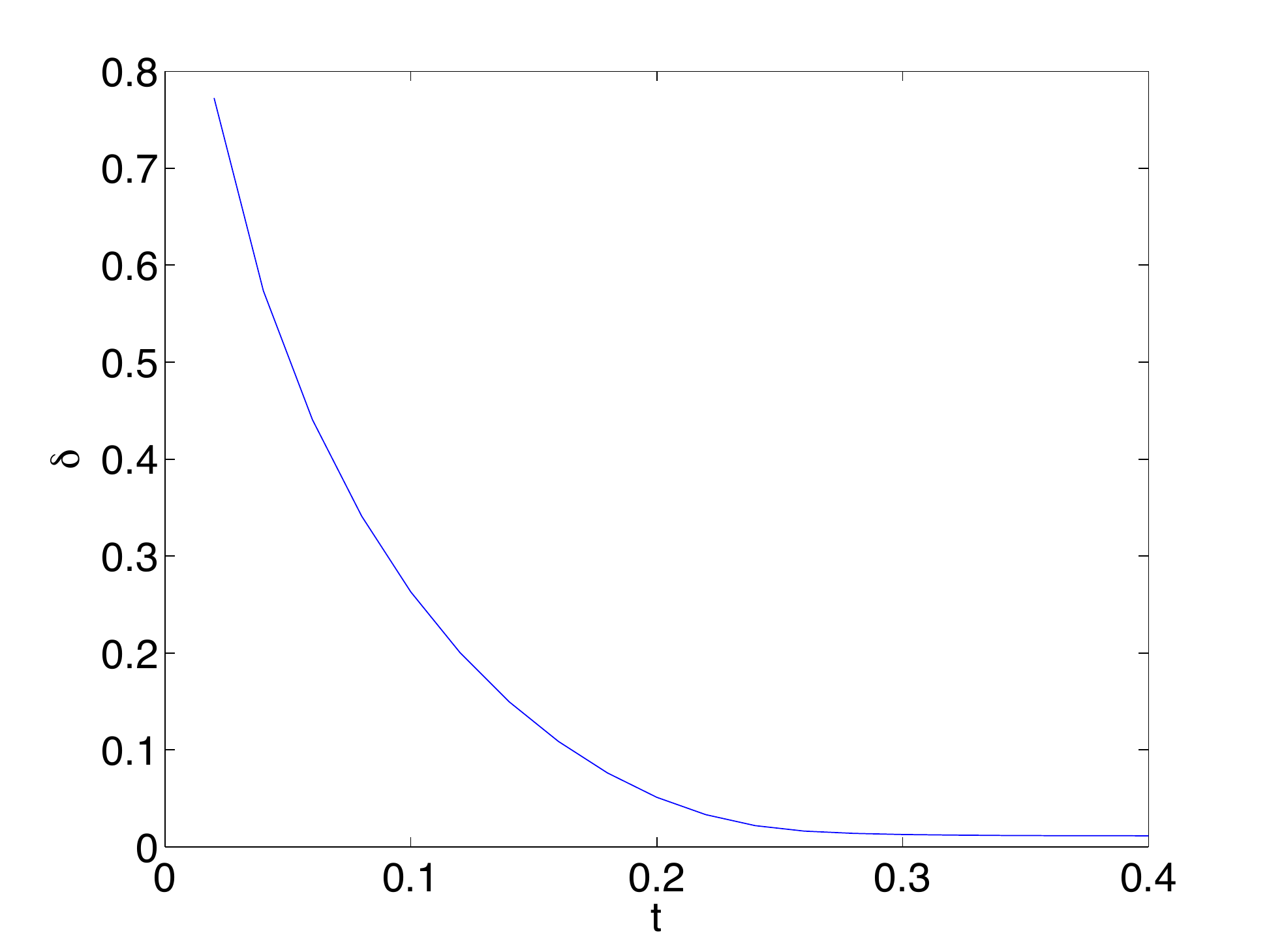}
  \includegraphics[width=0.5\textwidth]{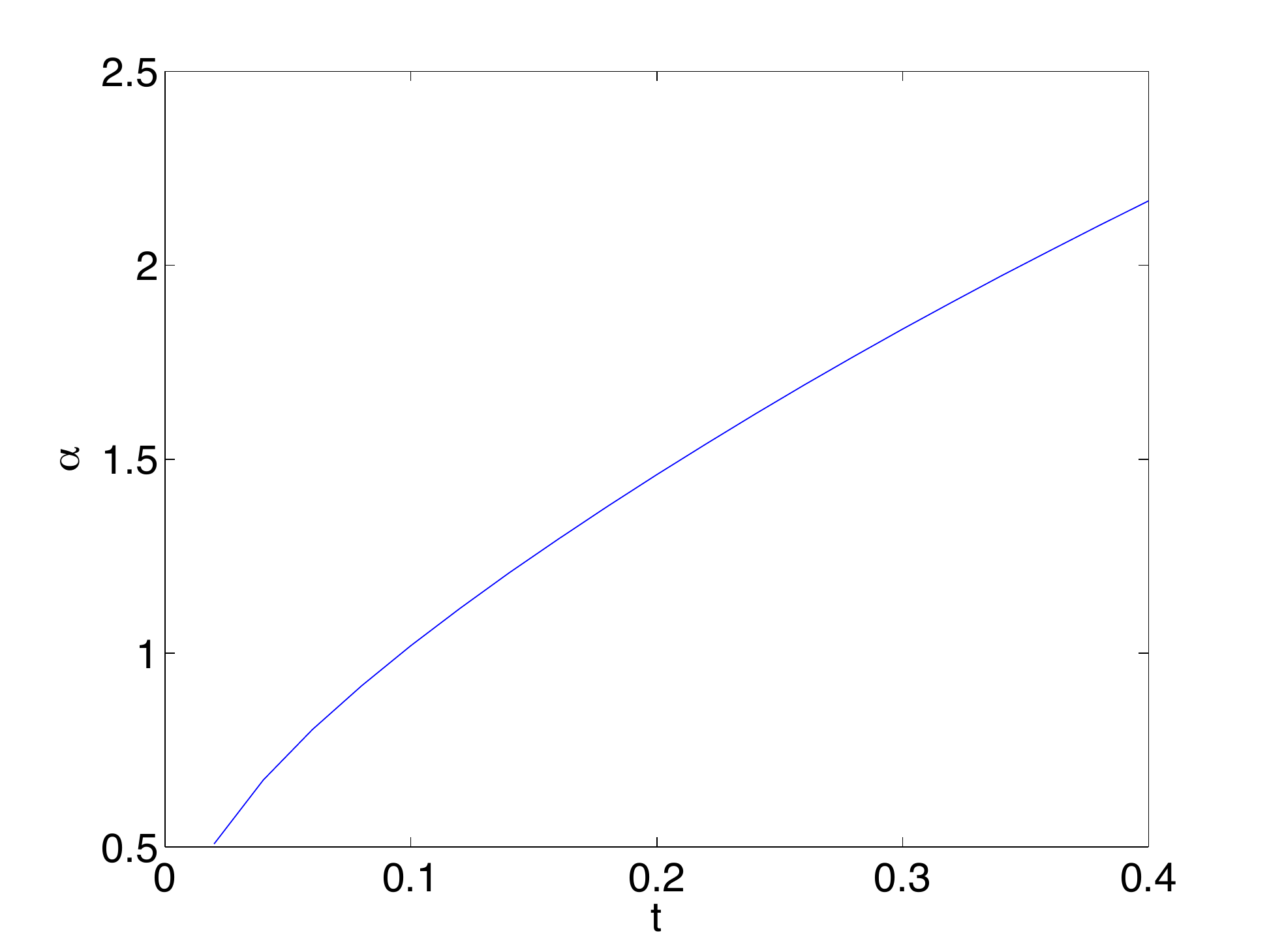}
 \caption{Fitting parameters for the solution to the Burgers' equation for the initial data 
 $u_{0}(x)=\mbox{sech}^{2}x$ for  $\epsilon=0.01$ in dependence on 
 the time for $\delta$
 on the left, and $\alpha$
 on the right.}
   \label{burgersalpha}
\end{figure}

If we trace the quantities $\alpha$ and $\delta$ during the 
computation as a function of time for $k>10$ (for $k>100$ there are 
not enough Fourier coefficients above the level of rounding errors
for small times), we get the expected behavior as 
can be seen in Fig.~\ref{burgersalpha}.  
The imaginary part $\delta$ of the singularity decreases strongly 
with $t$ before the critical time of the Hopf solution and 
comes close to the axis there. After $t_{c}$ the decrease becomes 
very slow, and the singularity stays at a finite distance from the 
real axis. This shows that the time $t_{c}$ cannot be really inferred 
from the Burgers' equation with small dissipation in a reliable way. The transition in all characteristic quantities of the 
solution is not distinguished at this time, and even for $t>t_{c}$ the 
quantity $\delta$ continues to decrease. Therefore to identify $t_{c}$ 
and $x_{c}$, it appears better 
to study the Fourier coefficients of the Hopf solution without 
dissipation to determine the critical point. 
The real part of the singularity $\alpha$ can also be seen 
in Fig.~\ref{burgersalpha}. It describes in some sense the 
motion of the dissipative shock. 

\subsection{Korteweg-de Vries equation}
\begin{figure}[htb!]
  \includegraphics[width=0.5\textwidth]{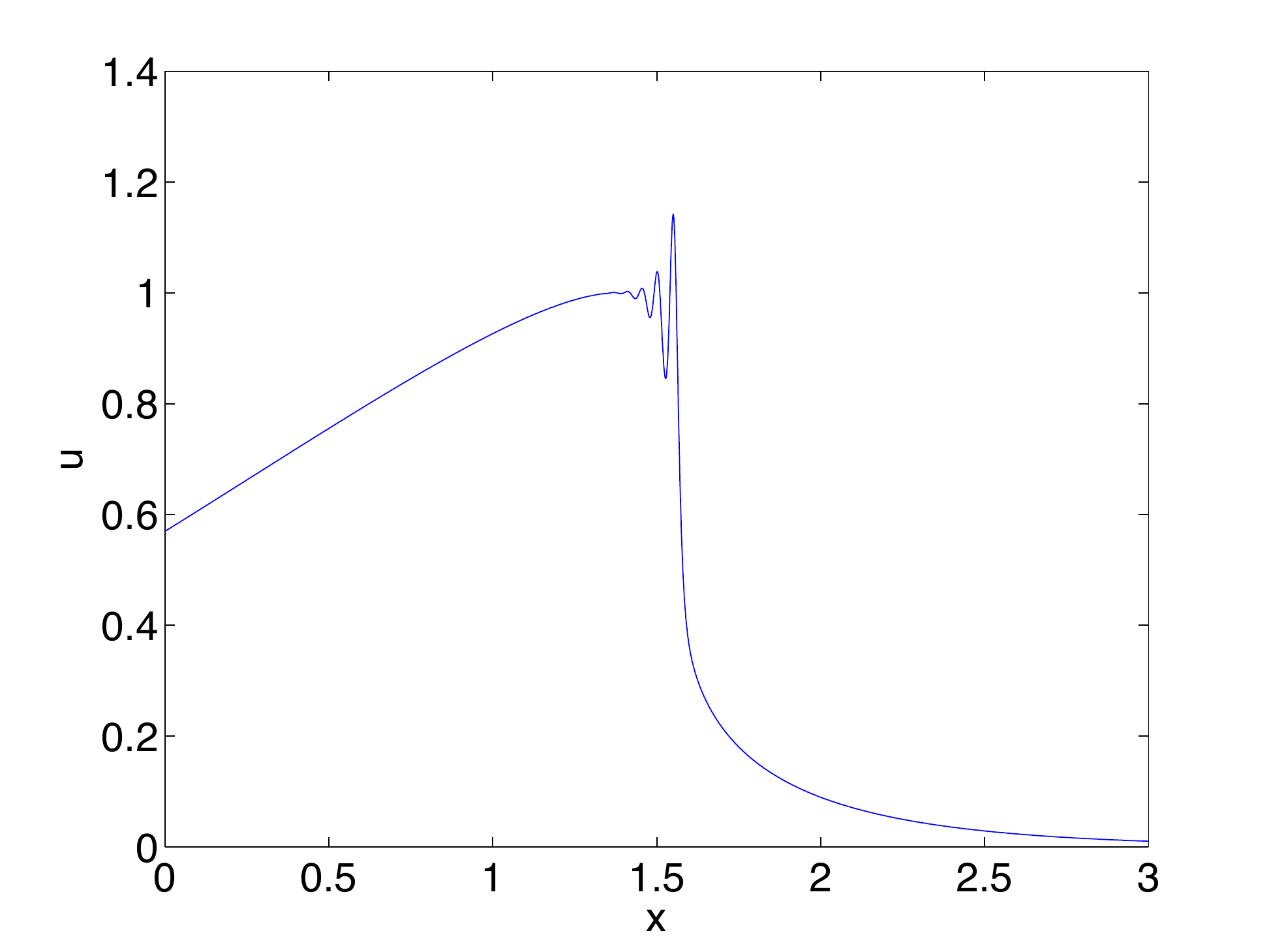}
  \includegraphics[width=0.5\textwidth]{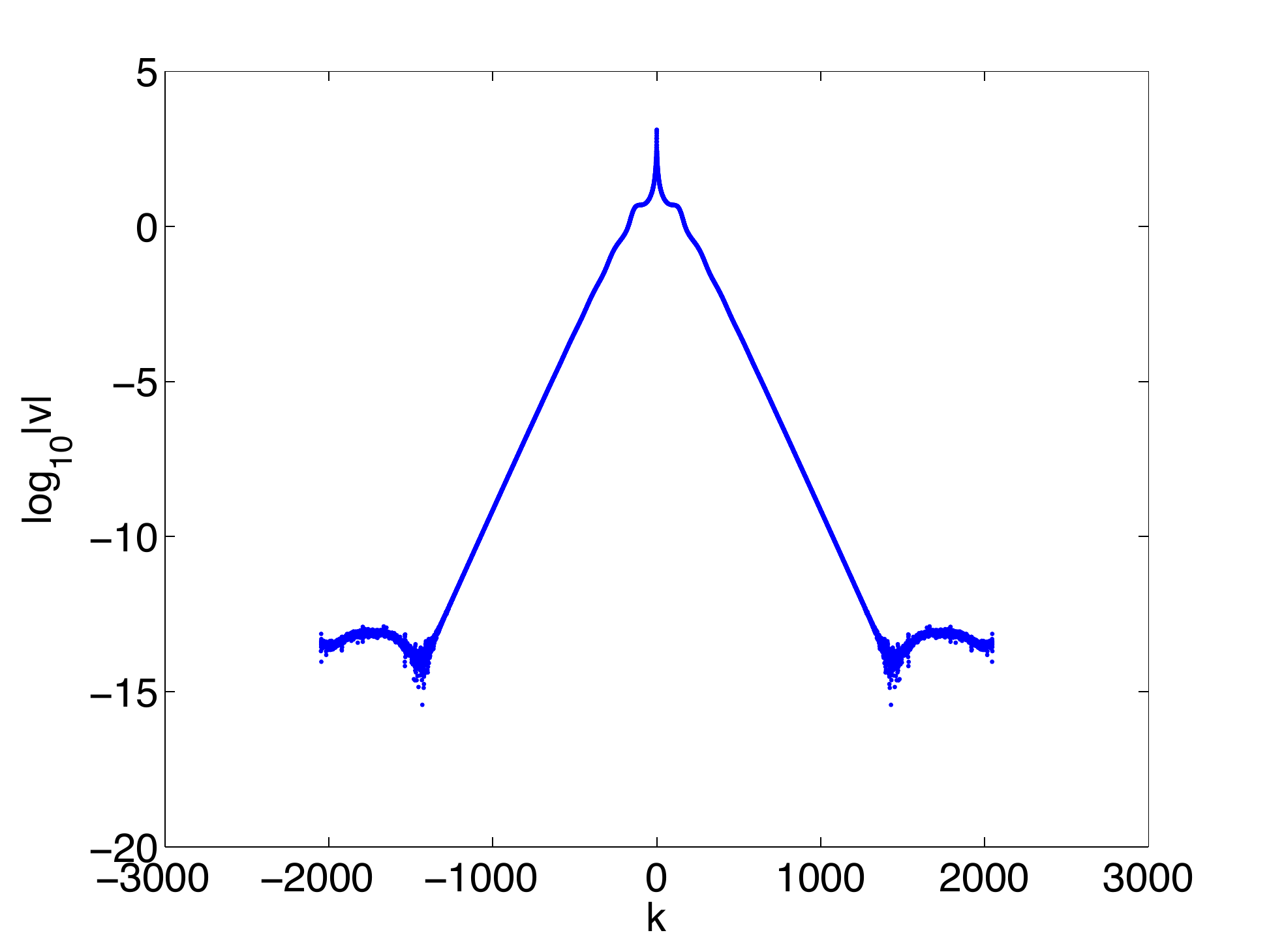}
 \caption{Solution to the KdV equation for the initial data 
 $u_{0}(x)=\mbox{sech}^{2}x$ for $t=0.23>t_{c}$ and $\epsilon=0.01$ 
 on the left, and 
 the modulus of the corresponding Fourier  coefficients 
 on the right.}
   \label{kdv1e4}
\end{figure}

The KdV equation (\ref{KdV}) with small dispersion $\epsilon\ll 1$
can be seen as a purely dispersive regularization of the Hopf 
equation. For the initial data studied there, there will be a zone of 
rapid modulated oscillations forming near the shock of the 
corresponding Hopf solution for $t$ greater than the critical time 
$t_{c}$ as can be seen for $\epsilon=0.01$ in Fig.~\ref{kdv1e4} at time 
$t=0.23>t_{c}$.

\begin{figure}[htb!]
  \includegraphics[width=0.5\textwidth]{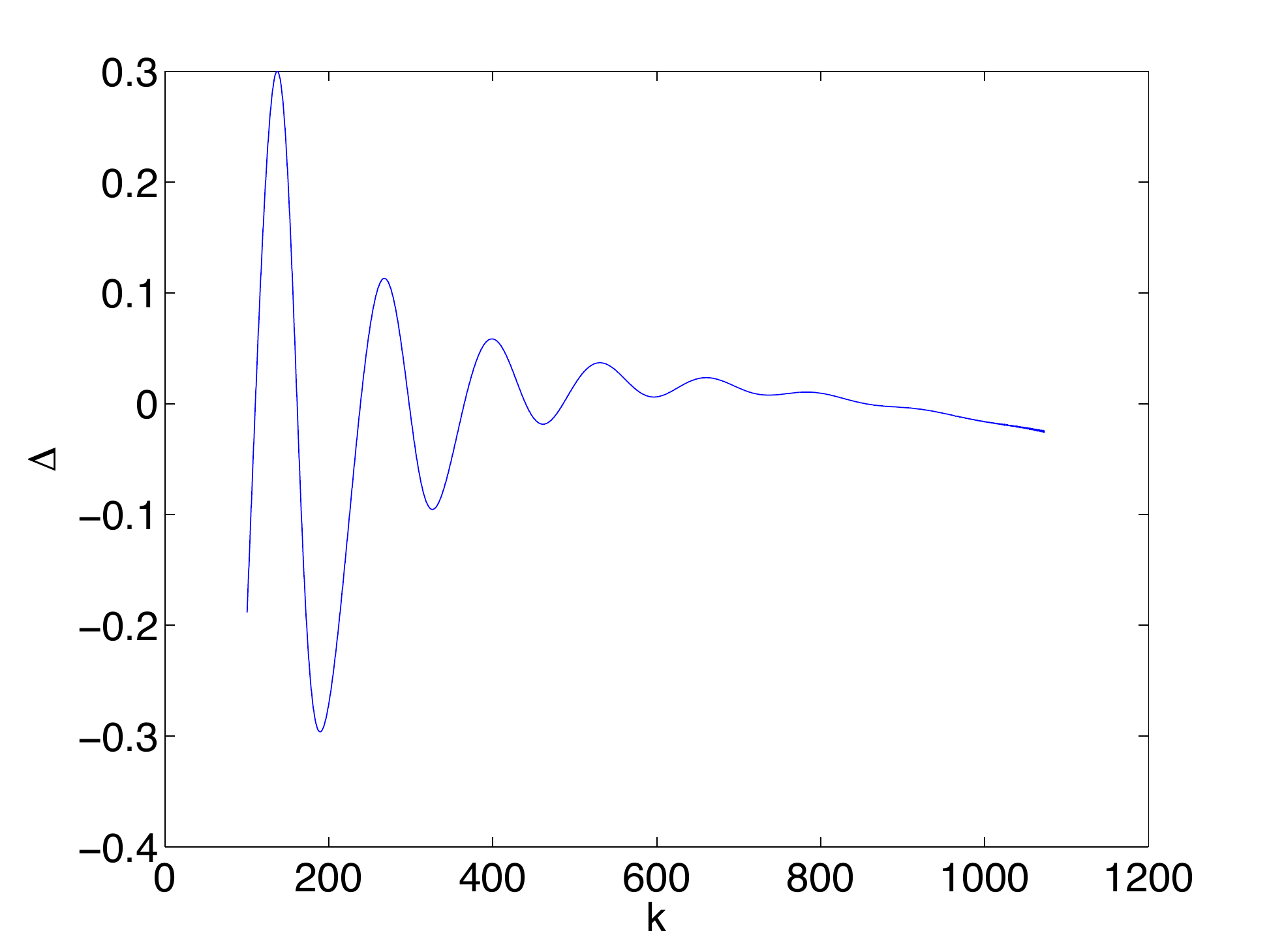}
  \includegraphics[width=0.5\textwidth]{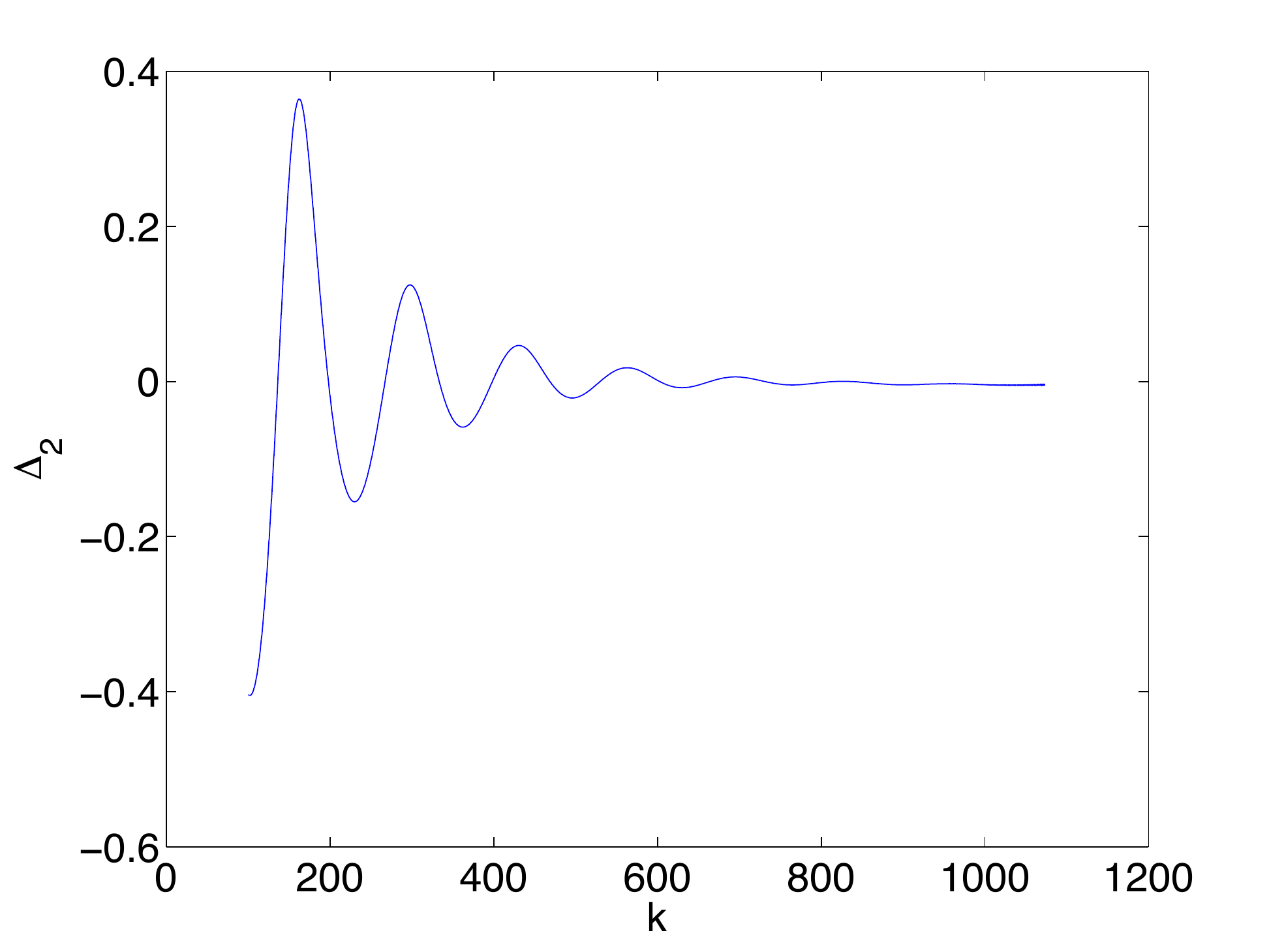}
 \caption{Difference between the modulus of the 
 Fourier coefficients for the situation of Fig.~\ref{kdv1e4}
 and the fitted asymptotic formula 
 (\ref{fourierasym}) for $k>10$ on the left, and 
 for the quantity $\phi$  (\ref{phi})
 on the right.}
   \label{kdv1e4fit}
\end{figure}

The Fourier coefficients for this solution are also shown in 
Fig.~\ref{kdv1e4}. It can be seen that they decrease (essentially 
exponentially) to machine precision. Doing a least square fitting as 
before for the modulus of the Fourier coefficients with 
$|v|>10^{-10}$ for $k>100$ (no 
further restrictions), we find  $A = 0.7638$, 
$B=0.8165$ and $\delta=0.0274$. This shows that the parameter 
$\delta$ indicating the distance between the nearest singularity to 
the real axis stays finite even for $t>t_{c}$ as expected. The 
difference between the modulus of the Fourier coefficients and the 
fitted asymptotic formula is shown in Fig.~\ref{kdv1e4fit}. It can be 
seen that there are damped oscillations for small $k$ which die out 
for large $|k|$ and which indicate the presence of several singularities in 
the complex plane in the vicinity of the real axis. A fitting to 
more than one singularity is problematic given the limited number of 
Fourier coefficients. Our study is here of a more qualitative 
nature. 

Fitting the quantity $\phi$ in (\ref{phi}), we get $\alpha=1.5517$ 
and $C=153.9262$. This implies that the real part of the location of 
the singularity is shifted towards larger values than for the critical 
value $x_{c}$ of the Hopf solution which is again not surprising since the 
Hopf shock propagates. The difference between $\phi$ and the 
fitted asymptotic formula is also shown in Fig.~\ref{kdv1e4fit}. 
There are damped oscillations in this difference for small $k$, too. 

\begin{figure}[htb!]
  \includegraphics[width=0.5\textwidth]{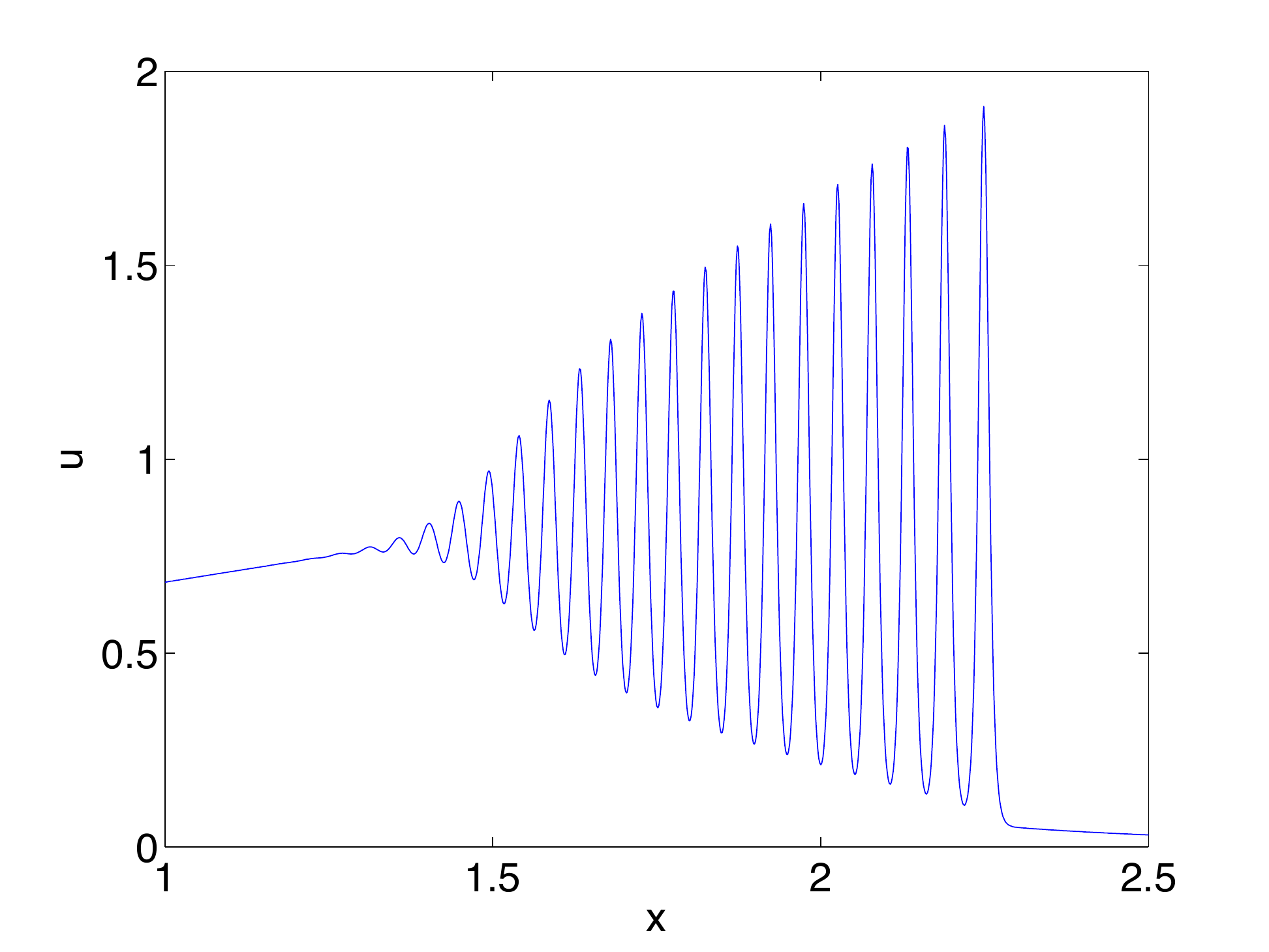}
  \includegraphics[width=0.5\textwidth]{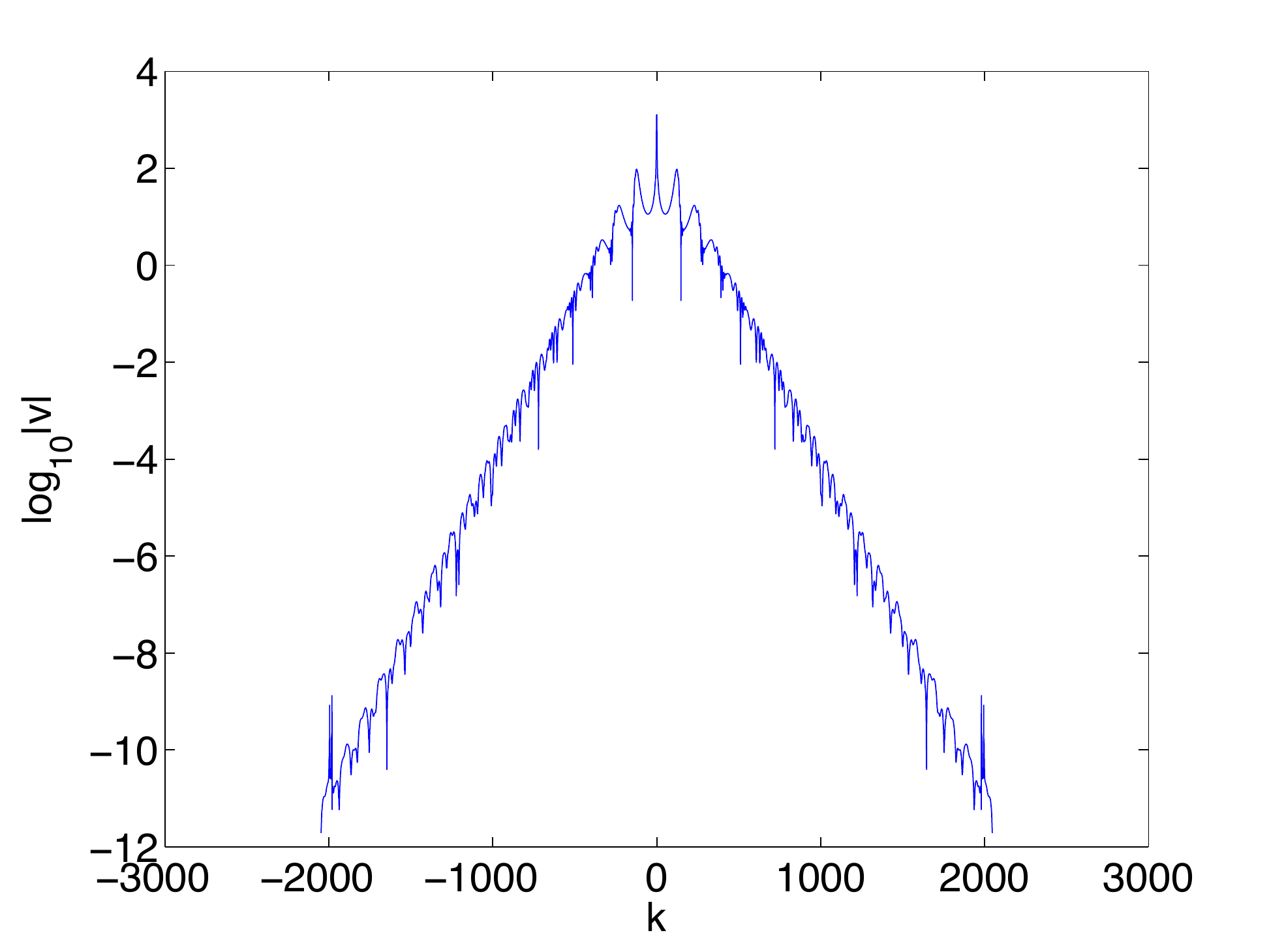}
 \caption{Solution to the KdV equation for the initial data 
 $u_{0}(x)=\mbox{sech}^{2}x$ for $t=0.4\gg t_{c}$ and $\epsilon=0.01$ 
 on the left, and 
 the modulus of the corresponding Fourier  coefficients 
 on the right.}
   \label{kdv1e404}
\end{figure}

For times $t\gg t_{c}$ there are many oscillations in the oscillatory 
zone as shown in Fig.~\ref{kdv1e404}. The higher number of 
oscillations also implies more oscillations in the Fourier 
coefficients as can be seen in Fig.~\ref{kdv1e404}.

\begin{figure}[htb!]
  \includegraphics[width=0.5\textwidth]{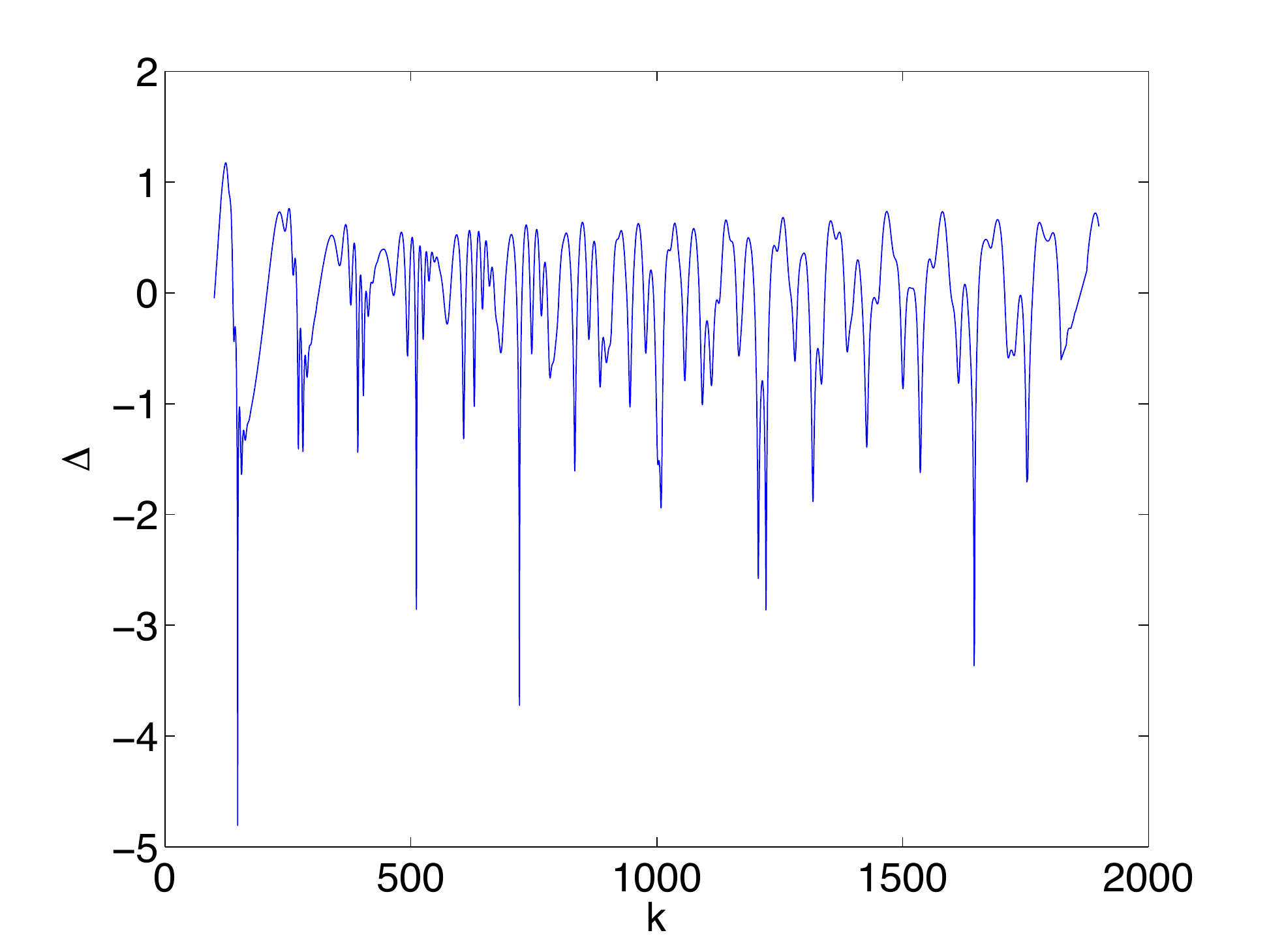}
  \includegraphics[width=0.5\textwidth]{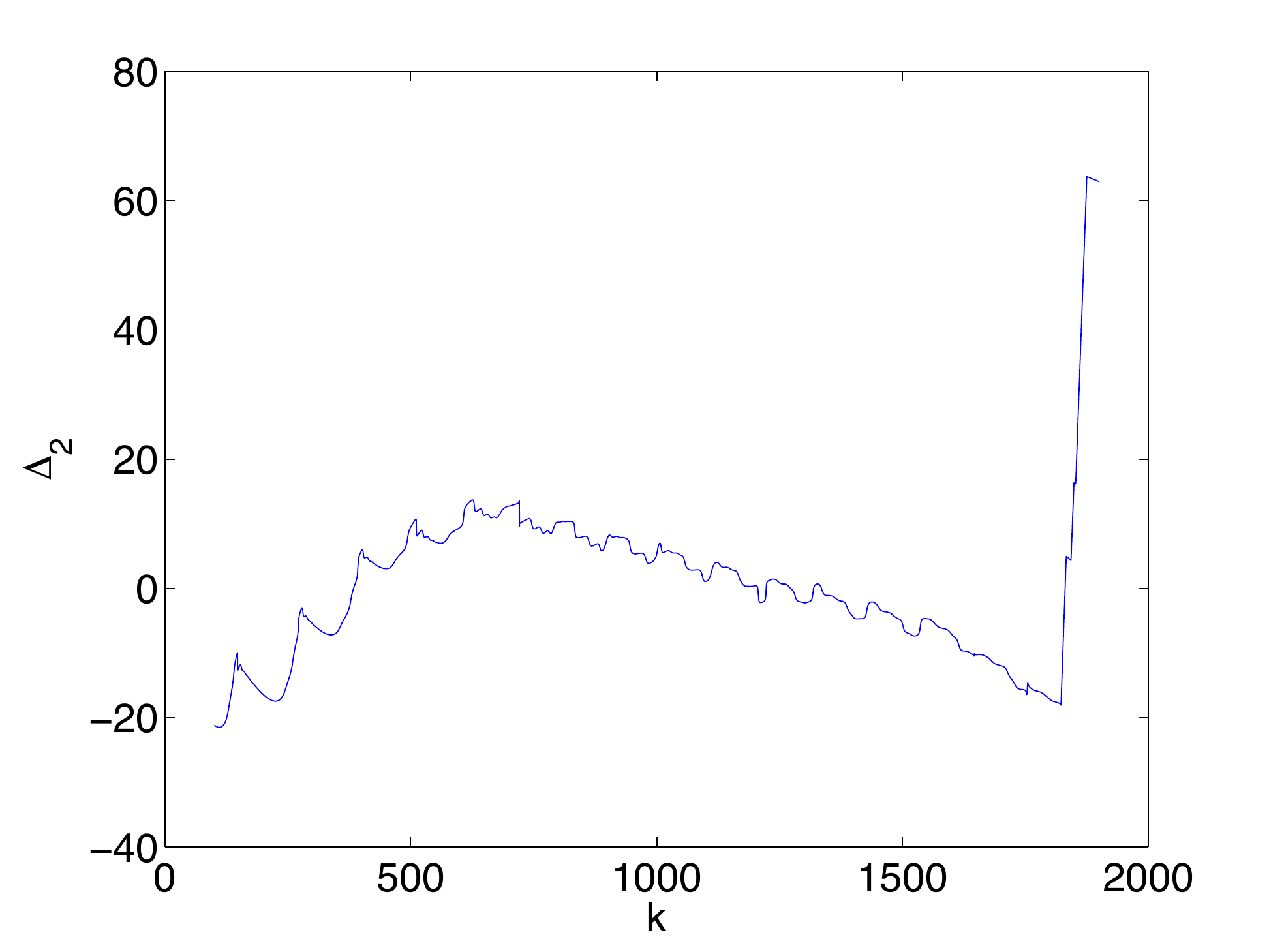}
 \caption{Difference between the modulus of the 
 Fourier coefficients for the situation of Fig.~\ref{kdv1e404}
 and the fitted asymptotic formula 
 (\ref{fourierasym}) for $k>10$ on the left, and 
 for the quantity $\phi$  (\ref{phi})
 on the right.}
   \label{kdv1e404fit}
\end{figure}

Consequently a least square fitting as 
before for the modulus of the Fourier coefficients for $100<k<1930$ 
(to exclude the last big oscillation in the Fourier coefficients 
which is presumably a numerical artifact), we find  $A = 0.7981$, 
$B=1.5684$ and $\delta=0.0164$. The 
difference between the modulus of the Fourier coefficients and the 
fitted asymptotic formula is shown in Fig.~\ref{kdv1e404fit}. In 
contrast to the case $t\sim t_{c}$ shown in Fig.~\ref{kdv1e4fit} where 
there are essentially just two big oscillations, there are many more 
here. For the modulus of the Fourier coefficients this implies with 
the asymptotic formula (\ref{fourierasym}) that one gets in the first 
case just standard damped oscillations of harmonic type.  In the 
latter case, several singularities in the complex plane approach the 
real axis which leads to more complicated oscillations (several 
phases) as can be seen in Fig.~\ref{kdv1e404fit}. Fitting the quantity 
$\phi$ in (\ref{phi}), we get $\alpha=2.1692$ 
and $C=241.7659$. The additional structure in the oscillations of the 
Fourier coefficients implies that a linear fitting of the phase 
$\phi$ is not efficient as can be seen in Fig.~\ref{kdv1e404fit}.

\begin{figure}[htb!]
  \includegraphics[width=0.5\textwidth]{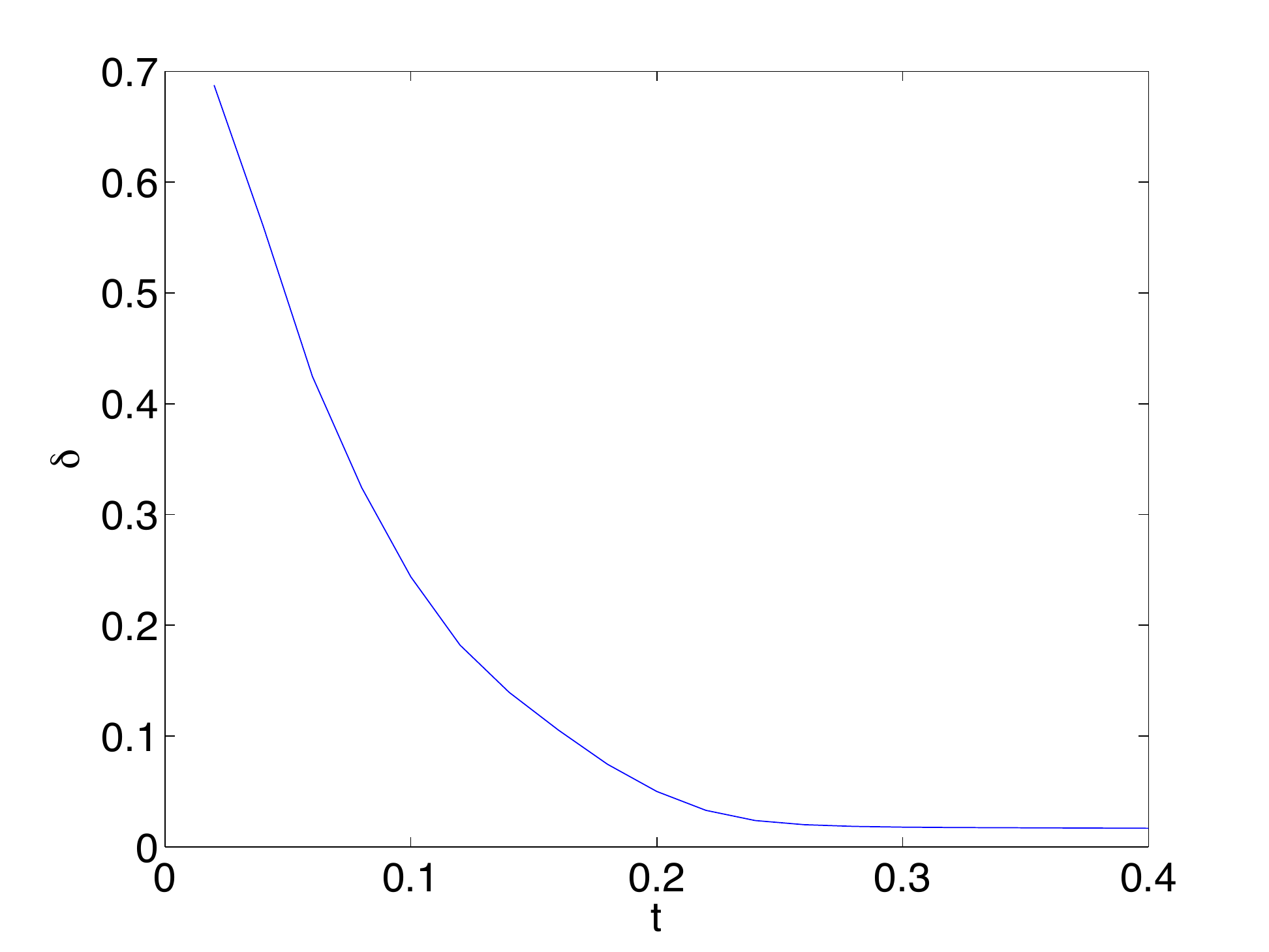}
  \includegraphics[width=0.5\textwidth]{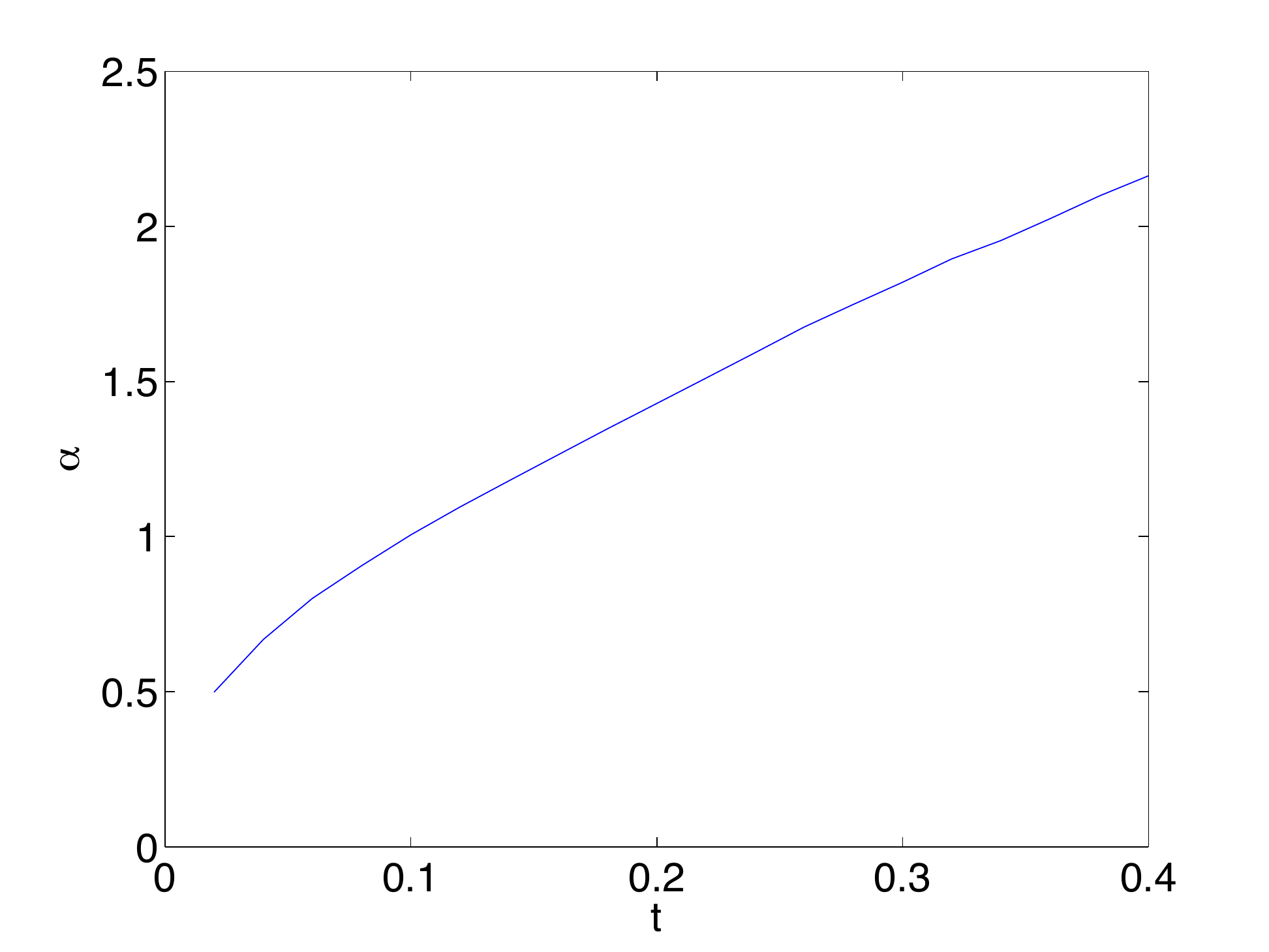}
 \caption{Fitting parameters for the solution to the KdV equation for the initial data 
 $u_{0}(x)=\mbox{sech}^{2}x$ for  $\epsilon=0.01$ for $\delta$
 on the left, and $\alpha$
 on the right.}
   \label{kdv1e4alpha}
\end{figure}

If we trace the quantities $\alpha$ and $\delta$ during the 
computation as a function of time for $k>10$ (again there are 
not enough non-trivial Fourier coefficients for $k>100$  at small times), we get the expected behavior as 
can be seen in Fig.~\ref{kdv1e4alpha}.  The used fitting just takes 
care of the singularity in the complex plane closest to the real axis. 
The imaginary part $\delta$ of the singularity decreases strongly 
with $t$ before the critical time of the Hopf solution and 
comes close to the axis there. 
After $t_{c}$ the decrease becomes 
very slow, and the singularity stays at a finite distance from the 
real axis. But again there is no clear indication
of the precise values of the critical quantities. The real part of the singularity $\alpha$ can also be seen 
in Fig.~\ref{kdv1e4alpha}. It gives in some sense the almost constant 
speed of the dispersive shock. 

\section{Kadomtsev-Petviashvili equations}
In this section we will study shock formation in the dKP equation 
with the methods discussed for the example of the Hopf equation in 
the previous section. This will be done  for both the dKP I and dKP II 
equations for two classes of initial data. The first class is 
motivated by the line solitons of KP, i.e., solutions exponentially 
localized in one spatial dimension and infinitely extended in the 
other. The second class of initial data are localized in both spatial 
directions. Then we will study for several values of the small 
dispersion parameter $\epsilon$ how the difference between the dKP 
solution and the corresponding KP solutions scales at the critical 
point and for $t\ll t_{c}$ in dependence of $\epsilon$.

\subsection{Theoretical preliminaries}
We fist collect some known analytic facts on KP and dKP equations 
 we will need in the following. 
The KP equation (\ref{e1}) is not in the standard form 
for a Cauchy problem since $t$ is not a timelike coordinate if it is 
discussed as a  standard second order PDE, see e.g.~\cite{KSM}. 
To be able to treat initial value problems in $t$, equation (\ref{e1}) can be 
rewritten in \emph{evolutionary 
form}, 
\begin{equation}
   \partial_{t}u+6u\partial_{x}u+\epsilon^{2}\partial_{xxx}u=
   -\lambda\partial_{x}^{-1}
   \partial_{yy}u,\,\,\lambda=\pm1
    \label{kpevol}.
\end{equation}
The antiderivative $\partial_{x}^{-1}$ in (\ref{kpevol}) is to be understood as 
the Fourier multiplier with the singular symbol $-i/k_{x}$. In the 
numerical computation this multiplier is regularized as described in \cite{KR}.  
Equations (\ref{kpevol}) and (\ref{e1}) are equivalent for certain 
classes of boundary conditions as periodic or rapidly decreasing at 
infinity, i.e., the ones studied in the present paper. The 
evolutionary form of dKP follows from (\ref{kpevol}) for $\epsilon=0$.

The divergence structure of the KP equations (\ref{e1}) has the well 
known consequence that 
\begin{equation}
    \int_{\mathbb{T}}^{}\partial_{yy}u(x,y,t)dx=0,\quad \forall t>0
    \label{const},
\end{equation}
even if this constraint is not verified at the initial time. As shown in \cite{FS,MST} 
the solution to a 
Cauchy problem not satisfying the constraint will not be smooth in 
time for $t=0$. This leads to numerical problems which could have a 
negative impact on the determination of break-up singularities in dKP 
equations, see \cite{KSM}.  Thus we consider only initial data subject to the 
constraint (\ref{const}). The first example with period $2\pi L_{y}$ 
in $y$ is 
\begin{equation}
u(x, y, 0) = \exp \left(- (x-\cos(y/L_y))^2 \right),
\label{uini2}
\end{equation}
whereas the second is just the derivative of a rapidly decreasing (in 
both spatial dimensions)
function as studied in \cite{KSM,KR},
\begin{equation}
u(x, y, 0) = - \partial_x \, \mbox{sech}^2 (R), \,\, R=\sqrt{x^2+y^2}. 
\label{uini1}
\end{equation}
It is known (see e.g.~\cite{KSM} for references and examples) 
that solutions to KP equations for initial data in the 
Schwarzian space of rapidly decreasing functions generically do not 
stay in this space, but develop tails with an algebraic fall off to 
infinity. These tails will lead to a Gibbs phenomenon in a periodic 
setting. To reduce the latter, we will choose a larger computational 
domain than would be necessary if the solutions remained in the 
Schwarzian space as for instance in the case of KdV. We always give 
the Fourier coefficients to show that these algebraic tails do not 
affect the results in our examples. 

The dKP equation has only three conserved quantities and thus 
does not belong to the family of completely 
integrable equations having an infinite number of conserved 
quantities. Therefore standard solution generating techniques as dressing 
transformations \cite{NMPZ} and linear Riemann-Hilbert problems 
cannot be applied in this context. In \cite{DMT01} it was shown that 
solutions to the dKP equation can be found in terms of Einstein-Weyl 
geometries which can be constructed with Twistor methods, see for 
instance  \cite{Pen76,WRS77}. However this approach is rather implicit if 
one wants to construct solutions for a given Cauchy problem and has 
not yet been used to this end. In \cite{FK04} it was shown that the 
dKP equation allows for an infinite number of hydrodynamic reductions 
and is thus integrable in this sense. Alinhac and coworkers 
\cite{alinhac} found that the 
dKP equation appears as a universal model in the geometric analysis 
of blow-up in nonlinear wave equations. They showed that finite time wave 
breaking is generic for solutions to this equation, proved existence of 
regular solutions up to the break-up time $t_{c}$ and classified the 
singularity. Localized initial data become singular 
in a point $(t_{c},x_{c},y_{c})$ and  the type of 
singularity is a one-dimensional cusp as for the Hopf equation, the 
second spatial dimension remains regular. However 
these methods do not provide formulae to determine  
$(t_{c},x_{c},y_{c})$ for given initial data.

The most explicit results for break-up in dKP solutions so far have 
been given by Manakov and Santini \cite{MS08}. They showed in 
\cite{MS06} that the dKP equation can be solved in implicit form  
\begin{equation}
    u = F (x - 2ut, y, t),
    \label{uF}
\end{equation}
similar to the solution of the Hopf equation via the method of 
characteristics. But $F$ is here an integral over one component of the solution of a 
nonlinear vector Riemann-Hilbert problem, see \cite{MS08} for 
details. The latter is equivalent to a systems of two coupled 
nonlinear, singular integral equations. If the latter system of 
integral equations is solved for a given normalization condition, the 
dKP solution follows from the implicit relation (\ref{uF}). To get a 
numerical solution to dKP in this way is not straight forward and has 
not been achieved for general initial data so far. But the form of 
the solution allows the study of singularity formation which has been 
done in \cite{MS08,MS12}. It is shown that there will be generically 
a gradient catastrophe at $(t_{c},x_{c},y_{c})$ 
in one spatial direction, whereas the solution 
remains smooth in the second spatial direction. Formulae are given 
for the solution for small $|t-t_{c}|$, $|x-x_{c}|$ and $|y-y_{c}|$, 
and the type of singularity is identified again as cubic in the singular 
direction.  
In particular, 
they show that at the breaking time $t_c$, if the initial condition 
$u_0(x,y)$ for the dispersionless KP equation is even in $y$, then 
the first breaking point is on the $x$-axis, i.e., $y_c=0$, and that 
the $y$-derivative of $u$ is zero at $(x_c,0)$. Our numerical 
results  for the case of localized initial data in Sect. 3.3. are in 
accordance with this prediction. But it is not clear how to obtain in 
the formalism \cite{MS08}  the values for the critical point $(t_{c},x_{c},y_{c})$ for given 
initial data. Therefore we solve the dKP equation in this paper 
directly.

The above mentioned theoretical results in \cite{alinhac,MS08} and references therein 
indicate that break-up in dKP solutions is one-dimensional. 
We choose 
initial data symmetric with respect to the $y$-coordinate to ensure a 
blow-up of the $x$-gradient only. However this is without loss of 
generality since the same approach as applied below can be used for 
general initial data. In this case one first has to determine 
numerically the direction of the strongest gradient which will 
indicate where break-up occurs. Since the rest is as for the symmetric 
initial data, we concentrate on this conceptionally simpler case. 

\subsection{Numerical approach}
We use here again a Fourier spectral method for the spatial coordinates 
and as for the Burgers' and KdV equations in the previous section 
a fourth order exponential time differencing scheme
 derived by Cox and Matthews in \cite{CM} for the time integration. 
We consider periodic (up to numerical precision) 
solutions in $x$ and $y$, i.e. solutions on $\mathbb{T}^2 \times \mathbb{R}$.  
The computations are carried out with $N_x \times N_y$ points for $(x, y) \in [-L_x\pi, L_x\pi] \times  [-L_y\pi, L_y\pi] $.

As discussed in the previous section, for a reliable identification of 
the shock formation via the asymptotic behavior of the Fourier 
coefficients, sufficient spatial resolution is crucial. In practice 
we needed $2^{14}$ to $2^{16}$ modes for the Hopf equation in $1+1$ 
dimensions.
To allow such high resolution simulations, we had to 
parellelize the codes also for the asymptotic analysis of the Fourier 
coefficients.
A prerequisite for parallel numerical algorithms is that sufficient 
independent computations can be identified for each processor,
that require only small amounts of data to be communicated between 
independent computations.
To this end,
we perform a data decomposition, which makes it possible to 
do basic operations on each object in the data domain (vector, matrix...) 
 to be executed safely in parallel by the available processors.

Our domain decomposition is implemented by developing a
code describing the local computations and local data structures for a single process. 
Global arrays are divided in the following way:
denoting by 
$x_n = 2 \pi n L_x/N_x,\,\, y_m = 2 \pi m L_y/N_y$, $n=-N_x/2,...,N_x/2, \,\, m=-N_y/2,...,N_y/2,$
the respective discretizations of $x$ and $y$ in the corresponding 
computational domain, 
$u$ is then represented by an $N_x \times N_y$ matrix. The latter is 
distributed among processors such that 
each processor $P_i, i=1...n_p$, ($n_p$ denoting the number of processors we can access) will receive $N_x \times \frac{N_y}{n_p}$ elements of 
$u$ corresponding to the elements 
\begin{equation}
u\left(1:N_x, (i-1).\frac{N_y}{n_p}+1 : i.\frac{N_y}{n_p} \right) 
\end{equation}
in the global array, and then 
each parallel task works on a portion of the data.

While processors execute
an operation, they may need values from other processors. The above domain decomposition 
has been chosen such that the distribution of operations is balanced 
and that the communication
is minimized. 
The access to remote elements has been implemented
via explicit communications, using  sub-routines of the MPI library.
The computation of the 2d-FFT is performed by 
using a transposition approach which takes advantage of the existing 
sequential FFT routines (more precisely, we use serial FFTW routines 
\cite{FJ}). The 
asymptotic fitting of the Fourier coefficients, here in one spatial 
direction, requires in addition 2 local communications, a core being reserved to do the required computations. More precisely, the 
processor which `takes care' of $v(k_x, 0, t)$ sends these data to 
the reserved processor,  which  receives the data (thus 2 local 
communications per time iteration, a SEND, and a RECEIVE), and then 
it alone performs the least square fit required, which takes less 
time than the computation of the solution for one iteration, 
performed by the other processors.  
The main 
technical problem in the use of  ETD schemes is the efficient and accurate 
numerical evaluation of the functions 
$$    \phi_{i}(z) = 
    \frac{1}{(i-1)!}\int_{0}^{1}e^{(1-\tau)z}\tau^{i-1}d\tau, \quad 
    i=1,2,3,4,$$
i.e., functions of the form $(e^{z}-1)/z$ and higher order 
generalizations thereof, as explained in \cite{KassT} and in 
\cite{Schme}. We use the approach given in \cite{Schme} in our 
implementation. 
The computation of these functions takes 
only negligible time for the $2+1$-dimensional equations studied here, 
especially since it has to be done only once 
during the time evolution. In addition each processor only computes 
the portion of these functions that the further local computations require, 
according to the domain decomposition specified above. 

To control the numerical accuracy we trace as for the Hopf equation 
in the previous section and as in \cite{KR} the 
conservation of the $L_{2}$ norm of $u$, the \emph{mass} of the 
solution. It is exactly conserved for the KP equations, but due to 
unavoidable numerical errors the computed mass $m(t)$ will numerically depend on time. 
Since mass conservation is not implemented in the code, it provides a 
valid check of the numerical accuracy for sufficient resolution in 
Fourier space.  It was shown in \cite{KR} that the quantity
\begin{equation}
    \Delta_{E}=\frac{m(t)-m(0)}{m(0)}
    \label{deltaE}
\end{equation}
typically overestimates the $L_{\infty}$ norm of the difference 
between numerical and exact solution by two orders of magnitude. 

\subsection{Shock formation in dKP solutions for 
initial data localized in one spatial dimension}
We will first study numerically the appearence of break-up singularities in 
solutions to the dKP  equations for the initial data  (\ref{uini2})
with the methods explained above. Whereas initial 
data localized in both spatial dimensions will have a gradient 
catastrophe in just one point, initial data localized in one 
dimension and infinitely extended in the other appear to blow up on a 
curve. 

\begin{figure}[htb!]
\centering
\includegraphics[width=\textwidth]{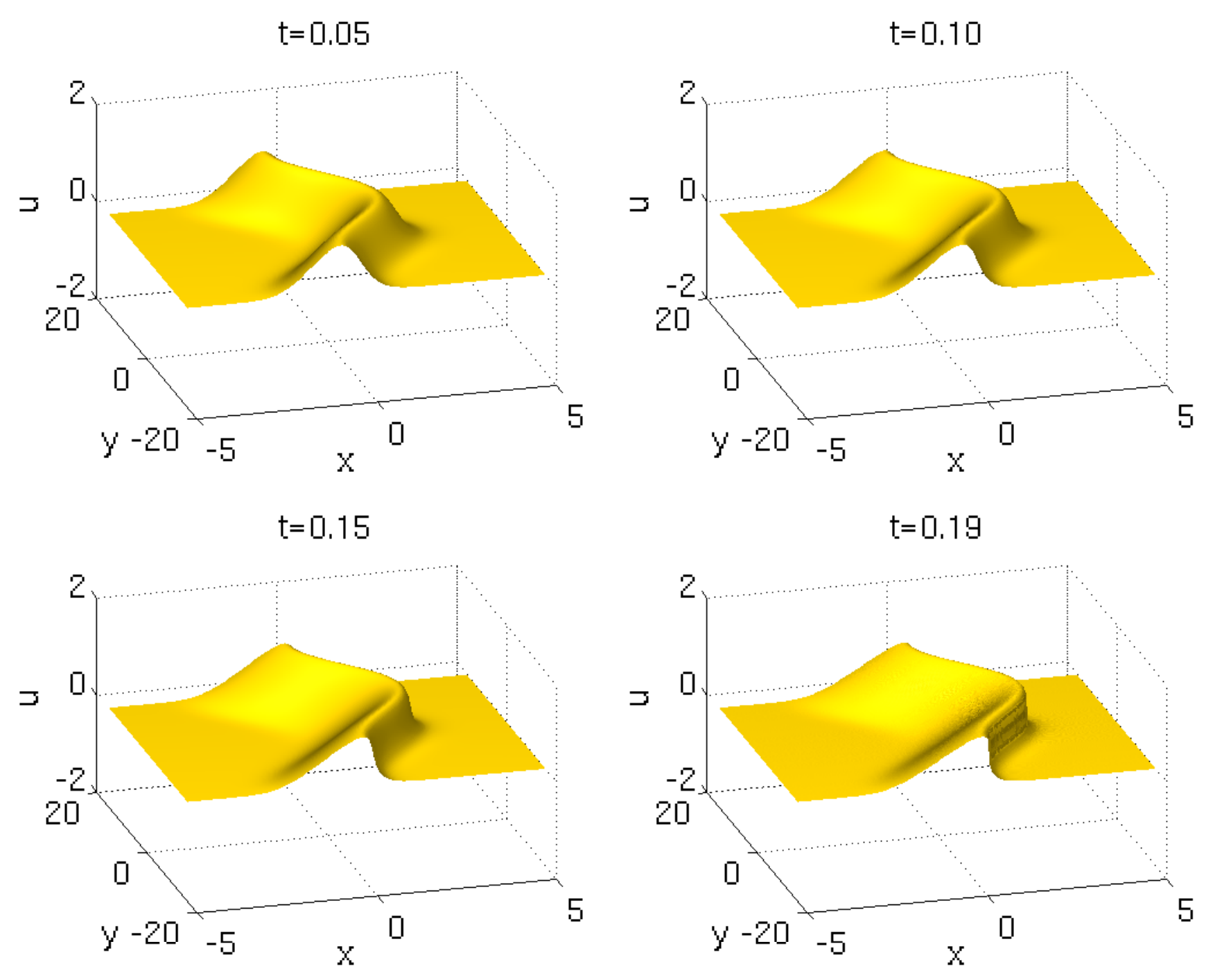} 
\caption{Solution to the dKP I equation for the initial data (\ref{uini2}), for several values of t}
\label{u1uts2d}
\end{figure}

To determine a solution for the initial data (\ref{uini2}), we use $N_x=N_y=2^{14}, L_x=L_y=5$. We observe that the 
solution develops a shock in the $x$-direction at time $t \sim 
0.19$, as can be seen in Fig. \ref{u1uts2d}.

\begin{figure}[htb!]
\centering
  \includegraphics[width=0.45\textwidth]{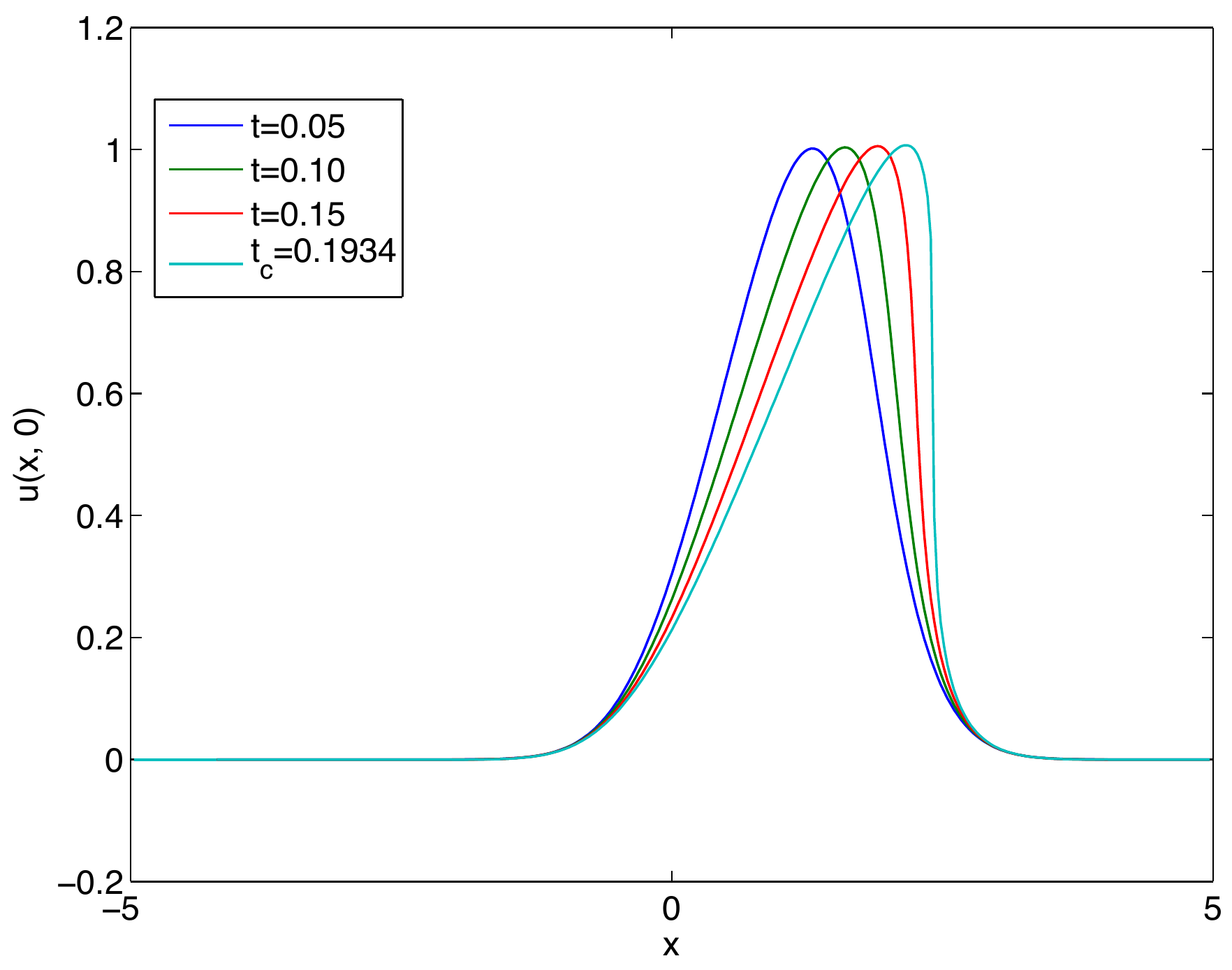} 
  \includegraphics[width=0.45\textwidth]{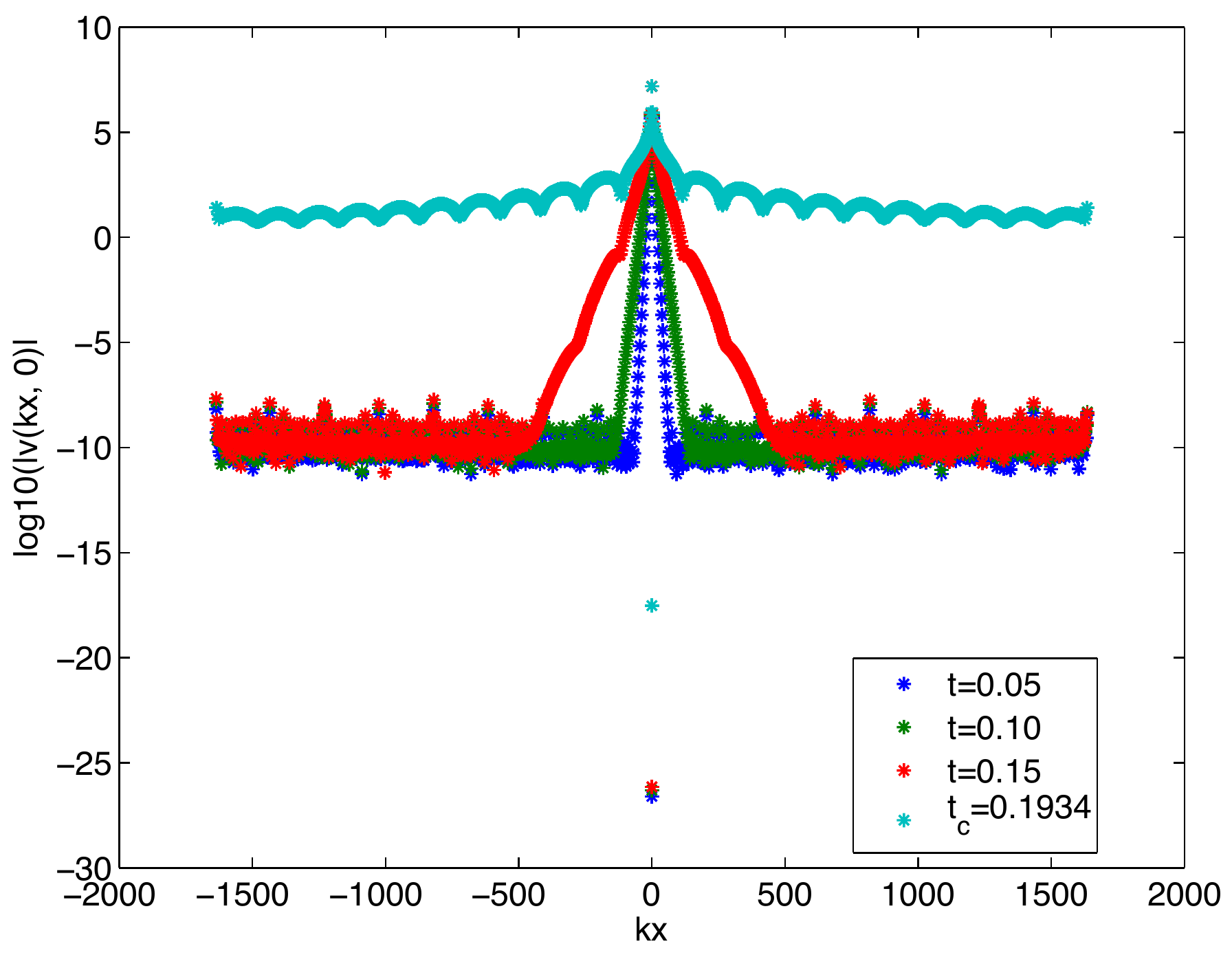}
 \caption{Solution to the dKP I equation for the initial data 
 (\ref{uini2}), plotted on the $x$-axis for several values of t, and the corresponding 
 Fourier coefficients, plotted on the $k_x$-axis.}
 \label{u1uts}
\end{figure}

\begin{figure}[htb!]
\centering
  \includegraphics[width=0.45\textwidth]{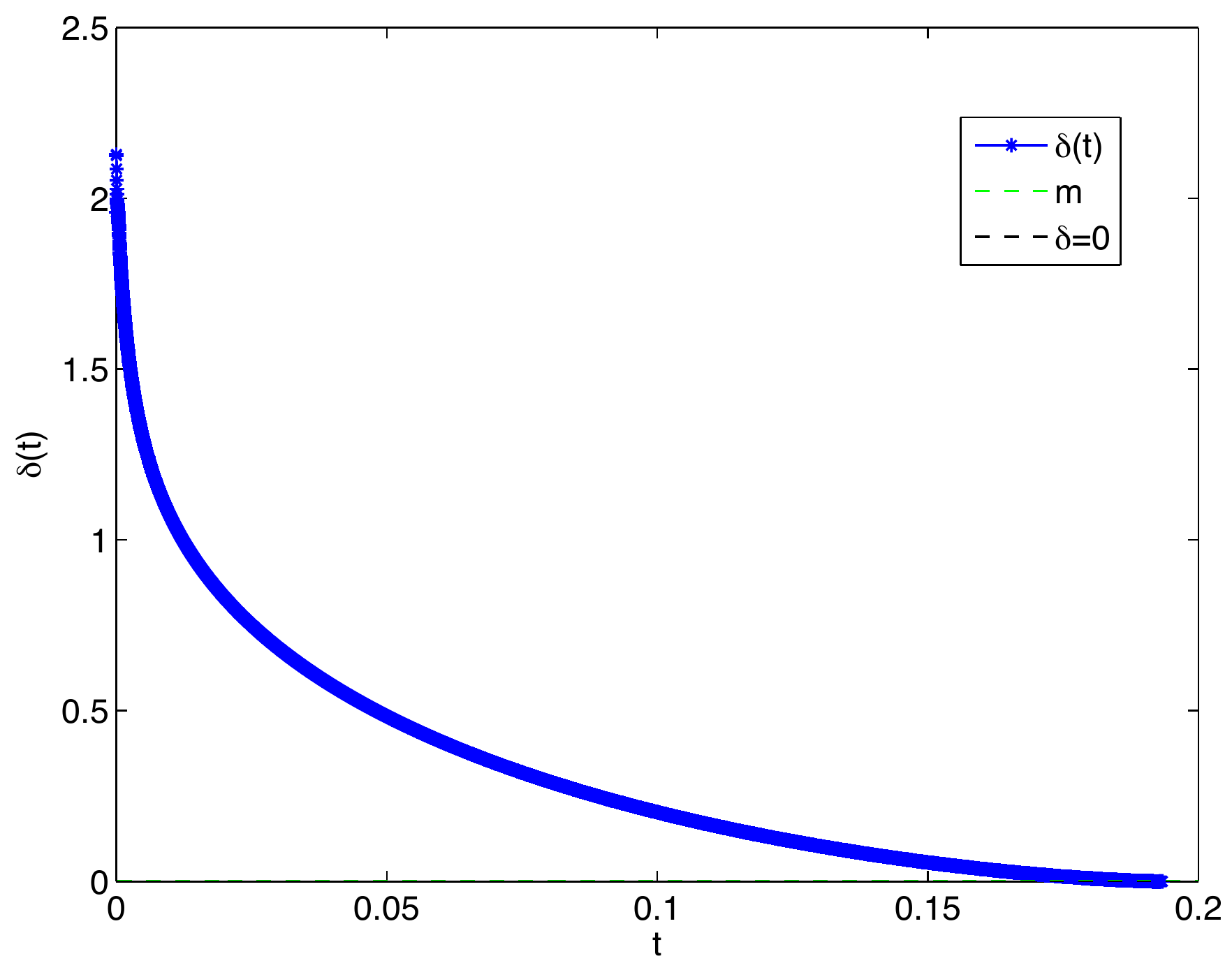}
  \includegraphics[width=0.45\textwidth]{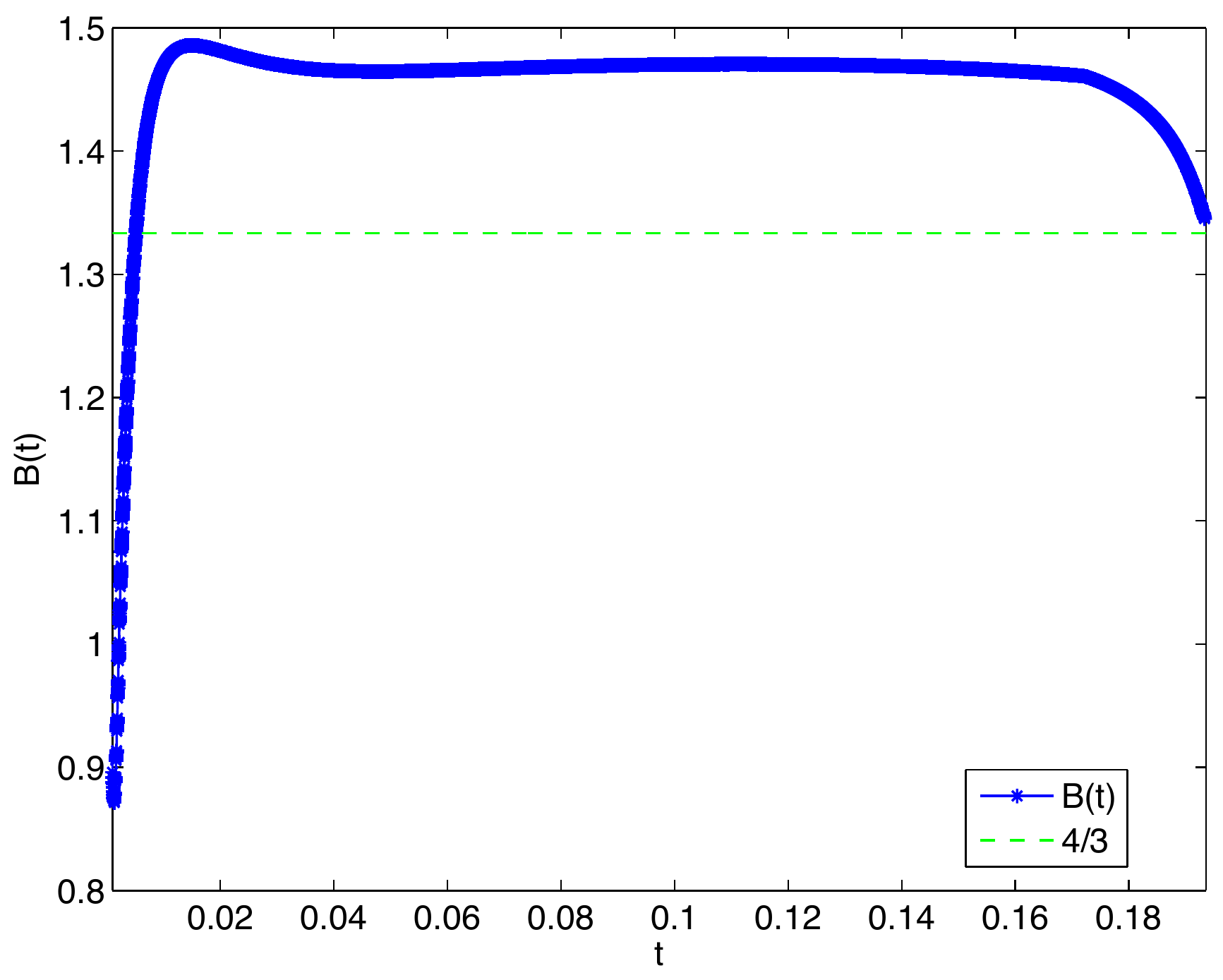}
 \caption{Time dependence of the fitting parameters $\delta$ (left), 
 $B$ (right) for the numerical solution to the dKP I 
 equation for initial data (\ref{uini2}).
The fitting  is done for  $5<k_x<\max(k_x)/2$. }
 \label{u1f1}
\end{figure}

\begin{figure}[htb!]
\centering
  \includegraphics[width=0.45\textwidth]{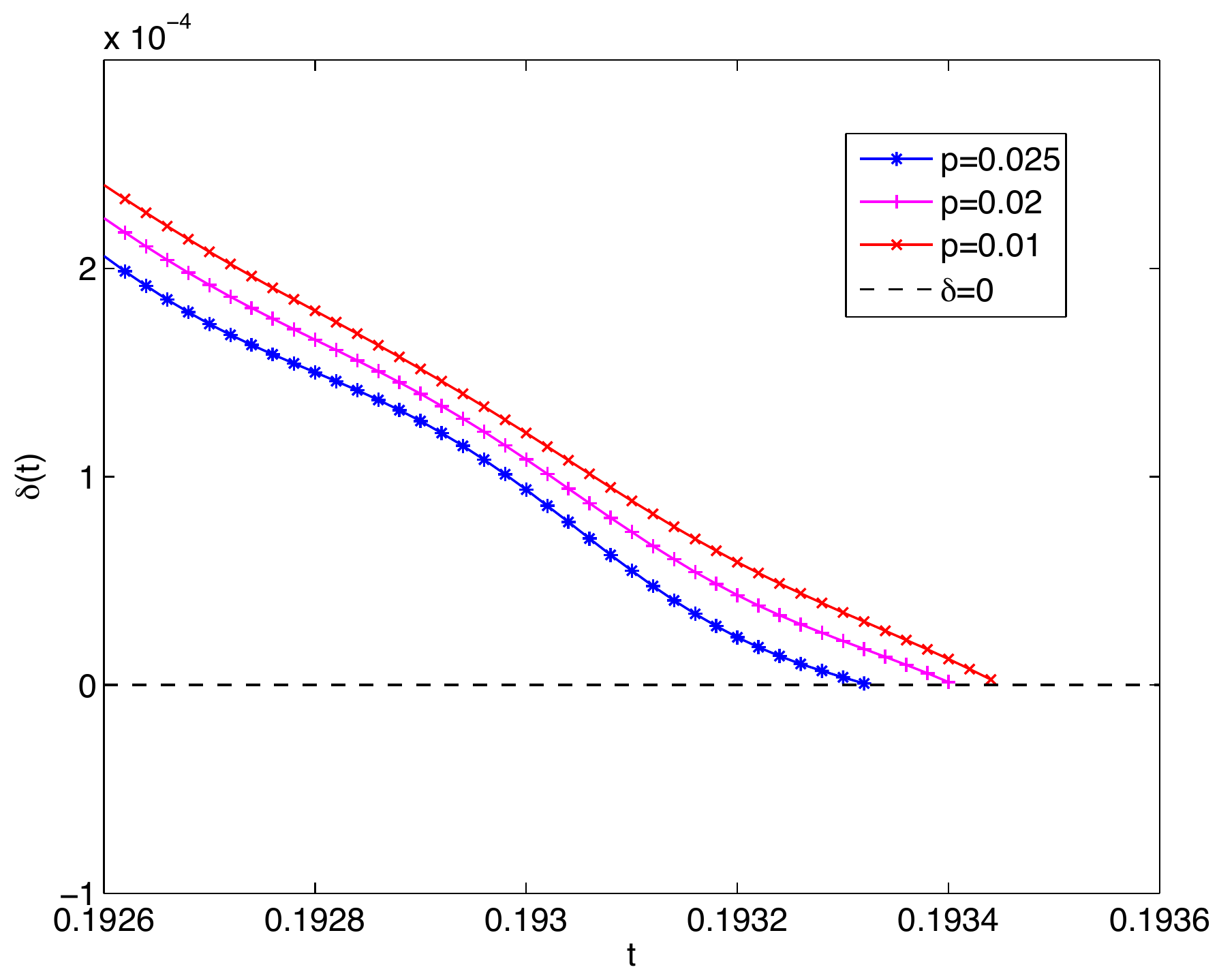}\includegraphics[width=0.45\textwidth]{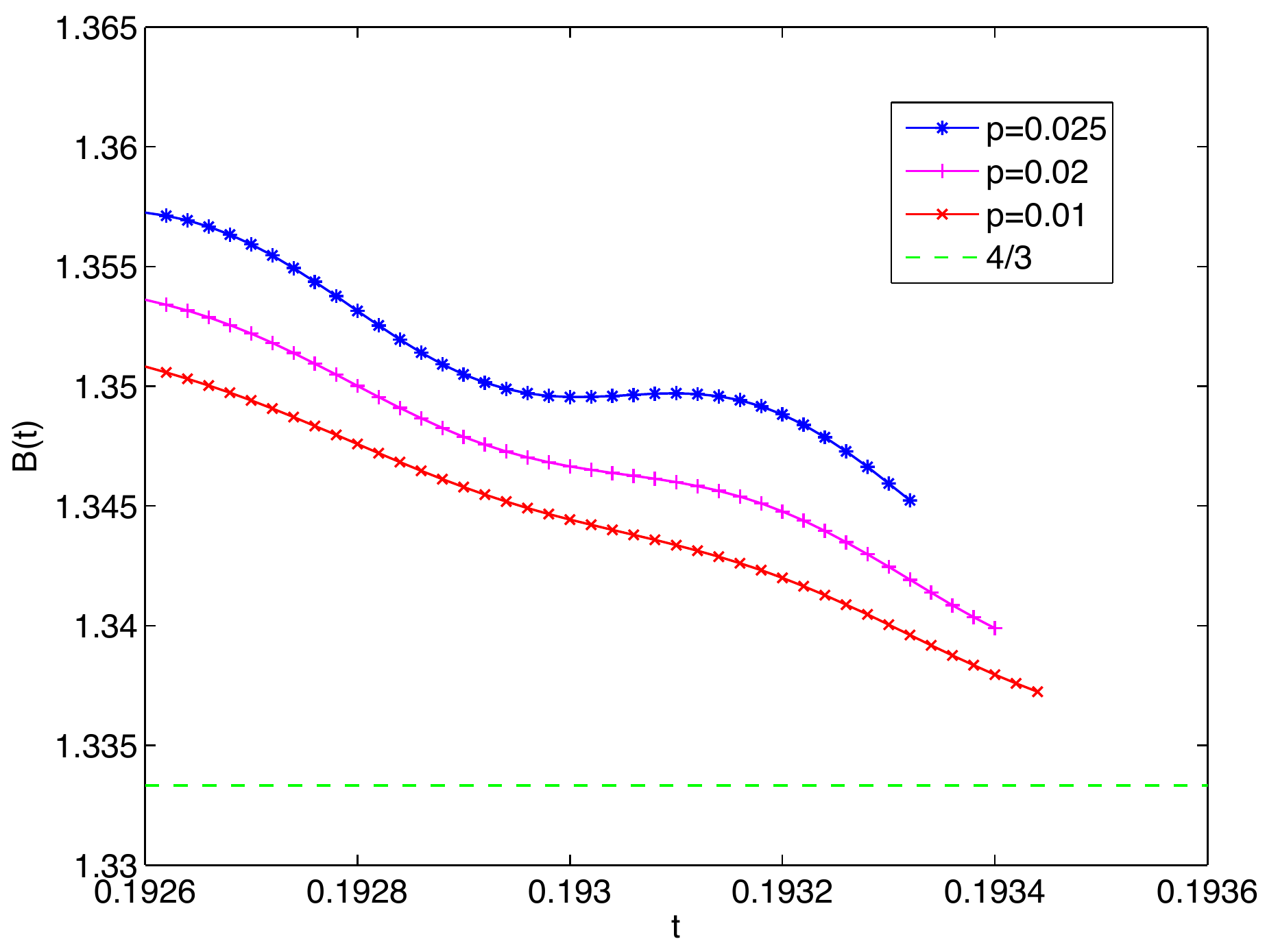}
 \caption{Time dependence of the fitting parameters $\delta$ (left), 
 and $B$ (right) close to $t_c\sim0.1934$ for the numerical solution 
 to the dKP I equation for the initial data \ref{uini2} for different values of $p$.}
 \label{u1fitp}
\end{figure}

To determine the time of the appearance of this shock, we 
apply the fitting procedure described in the previous section for 
the Fourier coefficients (again denoted by $v(k_x, k_y, t)$) in 
$x$-direction, i.e., $v(k_x, 0, t)$,  to formula (\ref{abd}).
As 
before (see Section 2.4) we use at least half of the Fourier coefficients with values 
above the rounding error, and an appropriately chosen interval 
$[k_{min},k_{max}]$.  
In this way, for $N_x=2^{14}$ we find the minimal possible precision to be $p=0.01$, for $k_{min}=5$ and $k_{max}=\max(k_x)/2$.
For this prescribed error,  we find that $t_c \sim 0.1934$, where $\delta$ vanishes, and $B(t_c)=1.337 \sim \frac{4}{3}$. We observe at this time the 
formation of a shock in the 
$x$-direction, which can be seen in Fig.~\ref{u1uts},
where we show the behavior of the solution at different times, plotted on the $x$-axis, and the corresponding Fourier coefficients, plotted on the $k_x$-axis.
The time dependence of the fitting parameters is shown in Fig.~\ref{u1f1}.
As expected, we observed a rapid decrease of $\delta(t)$, and the 
change of $B(t)$ 
from a value $\sim 1.5$ to the value $\sim 1.33$ when approaching the critical time, determined by the vanishing of $\delta(t)$.  
In Fig. \ref{u1fitp} we show the time dependence of the fitting parameters close to $t_c$, for different values of $p$.
By doing a fitting of $\Im\log(v)$, we find that at $t_c$, $x_c=\alpha(t_c)=2.4134$.

In this case, we  reach an error in the fitting of $p \sim 0.01$, 
which is the same as obtained  for the Hopf equation. Note also that the approximation of $B$ is of the same precision as in the latter. 
However, such a 
precision cannot always be achieved as will be seen below.

\begin{figure}[htb!]
\centering
  \includegraphics[width=0.45\textwidth]{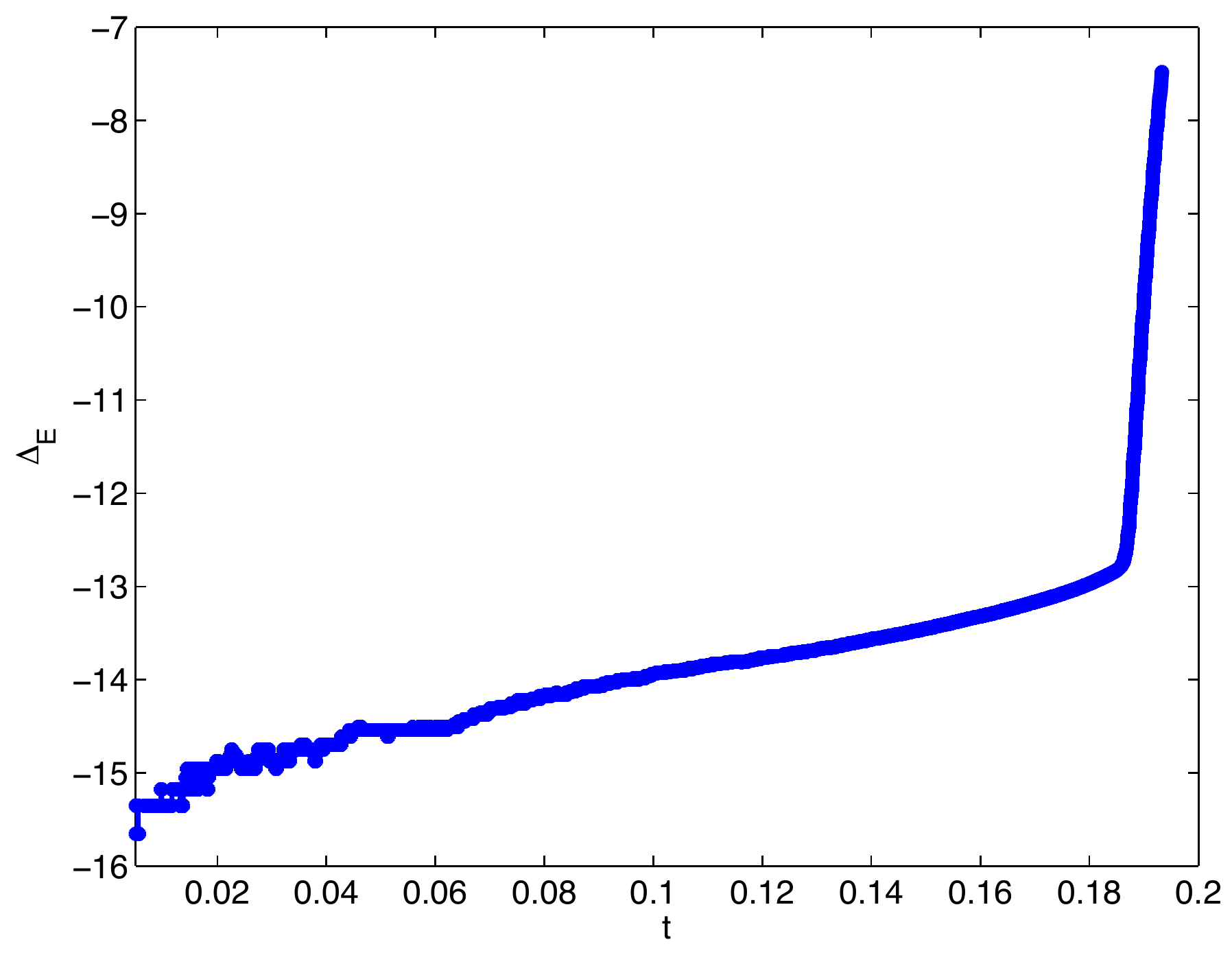} 
    \includegraphics[width=0.45\textwidth]{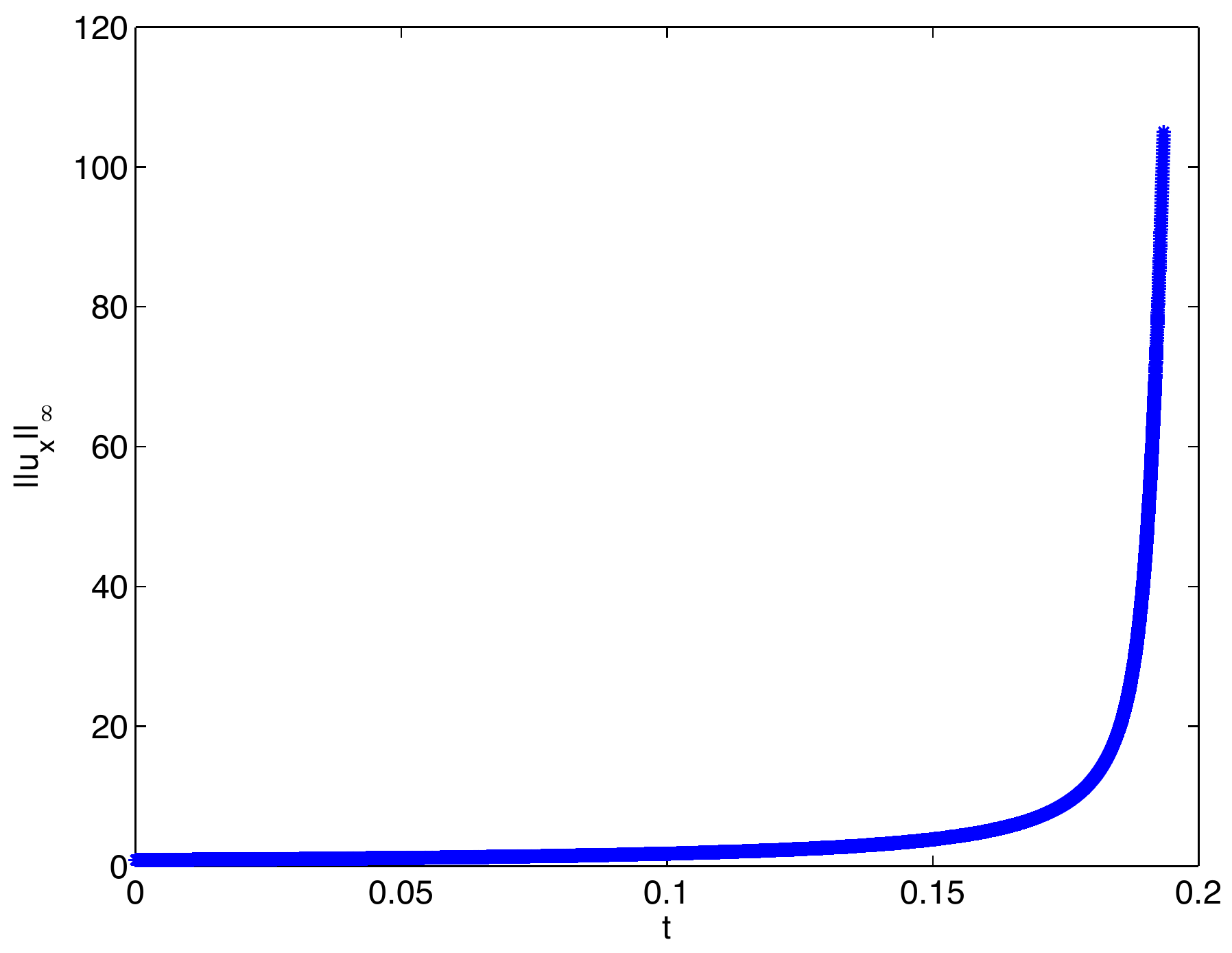}
 \caption{Time evolution of the quantity $\Delta_{E}$ as an indicator 
 of the numerical accuracy (left), and of 
 $\|\partial_x u \|_{\infty}$ (right) for the solution to the dKP I 
 equation for initial data (\ref{uini2}).}
 \label{massampl}
\end{figure}

The time evolution of the quantity $\Delta_E$ in 
(\ref{deltaE}), used as an indicator of the numerical accuracy, and of the $L_{\infty}$-norm of $\partial_x u$ 
until the critical time are shown in Fig.~\ref{massampl}. 
The quantity $\Delta_{E}$ indicating the numerical mass conservation 
increases close to $t_c$, but stays below $10^{-7}$, which implies 
that the system is still well resolved up to the shock formation.
Whereas the $L_{\infty}$ norm of $u$ stays finite, there are strong
indications of a blow-up of $u_{x}$ close to the critical time, which 
indicates that the fitting in fact identified correctly the gradient 
catastrophe. 

\begin{figure}[htb!]
\centering
  \includegraphics[width=0.45\textwidth]{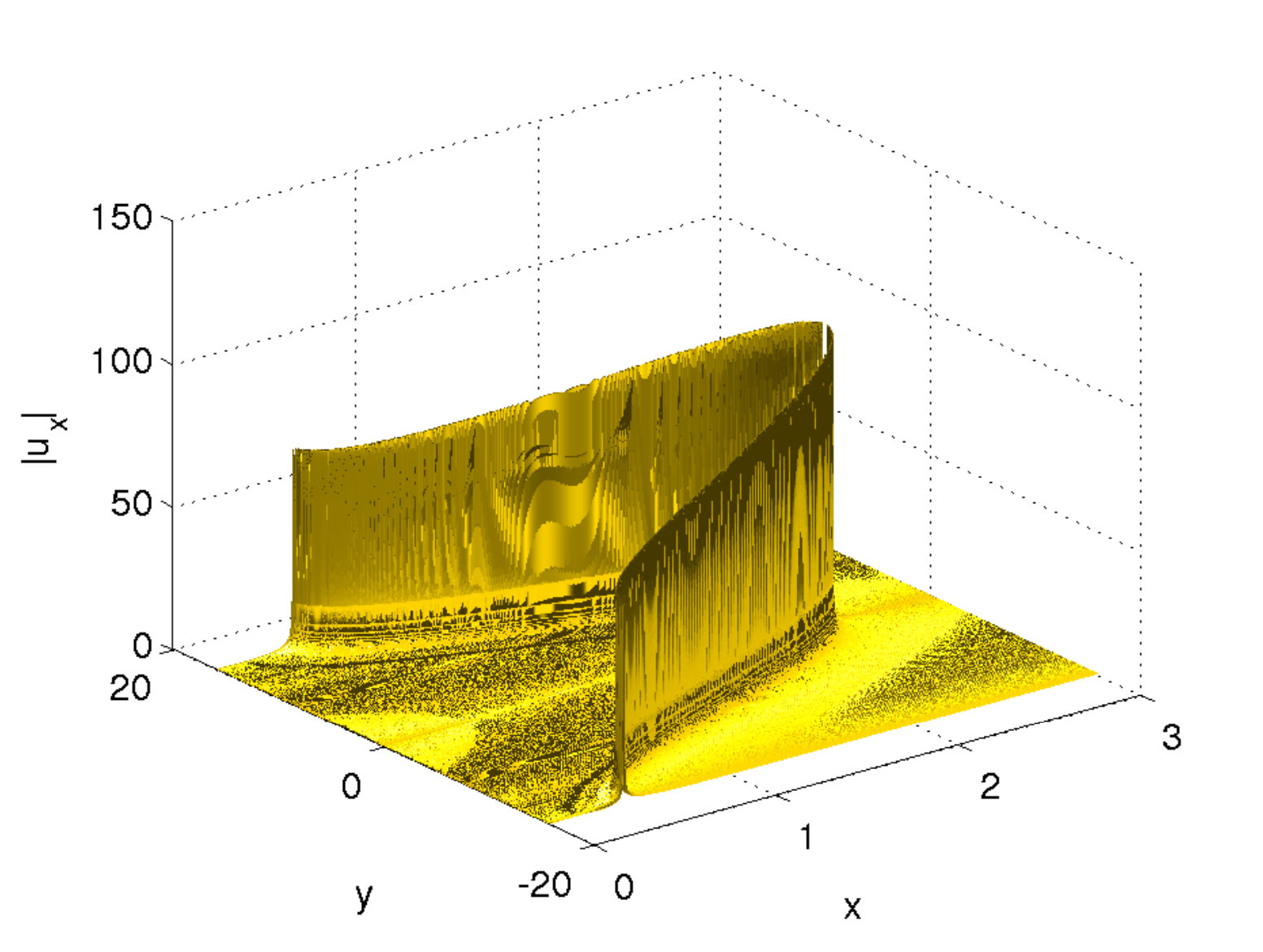} 
  \includegraphics[width=0.45\textwidth]{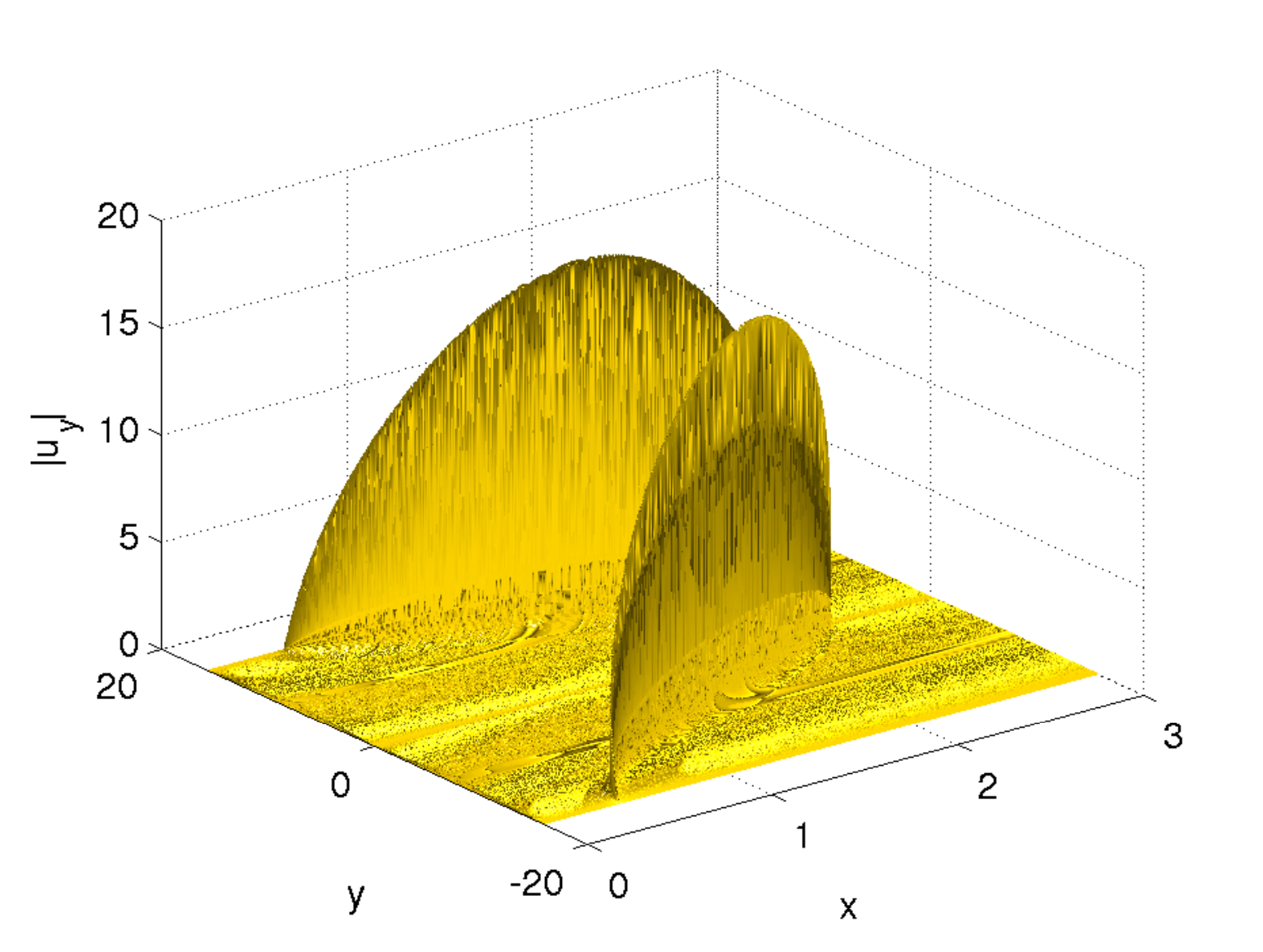}
 \caption{Derivatives  of the dKP I solution for the initial data 
 (\ref{uini2}) at the critical time $t_{c}=0.1934$, $|\partial_x u|$ (left) and of $|\partial_y u|$.}
 \label{derivativesofu}
\end{figure}

As one can see in Fig.~\ref{derivativesofu}, where  
$|\partial_x u|$  and $|\partial_y u|$  are plotted at 
$t=0.1934$, the gradient catastrophe does not appear in only one 
spatial point for the infinitely extended initial data.     

\begin{figure}[htb!]
\centering
  \includegraphics[width=0.45\textwidth]{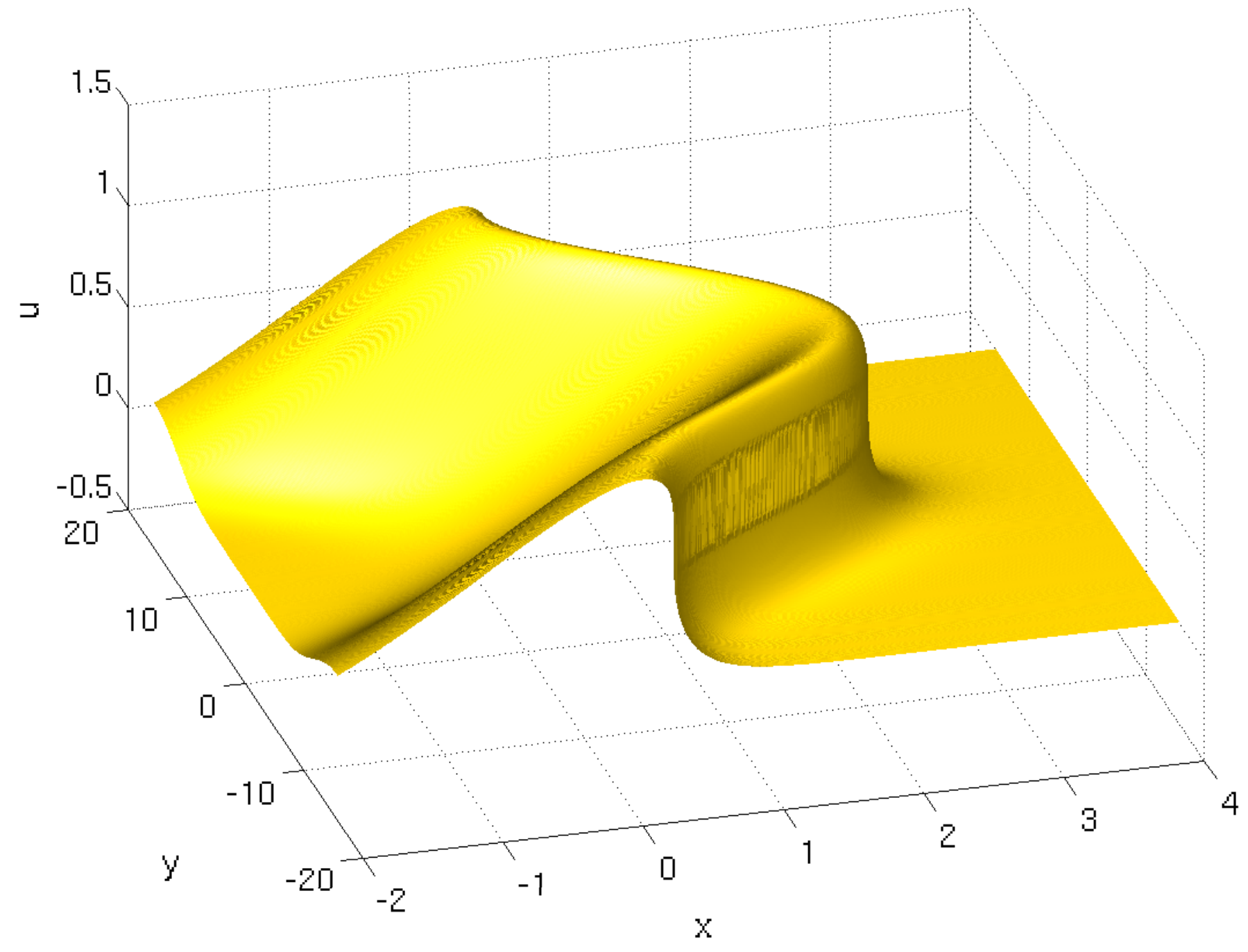} 
 \caption{Solution to the dKP II equation for the initial data (\ref{uini2}), at $t_c=0.1934$.}
 \label{udkp2exp}
\end{figure}

The situation is rather similar for the dKP II case, for the same initial data (\ref{uini2}).
By using the same fitting bounds as before, ($k_{min}=5$ and 
$k_{max}=\max(k)/2$), we find that the solution of the dKP II equation develops a shock at $t_c \sim 0.1934$, 
which is exactly the previously identified value for the critical 
time for dKP I.  The precision we could achieve here is also similar, $\Delta=0.01$, $B(t_c)=1.337$ 
and $\alpha(t_c)=2.41$. The solution at $t_c = 0.1934$ is shown in Fig. \ref{udkp2exp}. 
  
\begin{figure}[htb!]
\centering
  \includegraphics[width=0.45\textwidth]{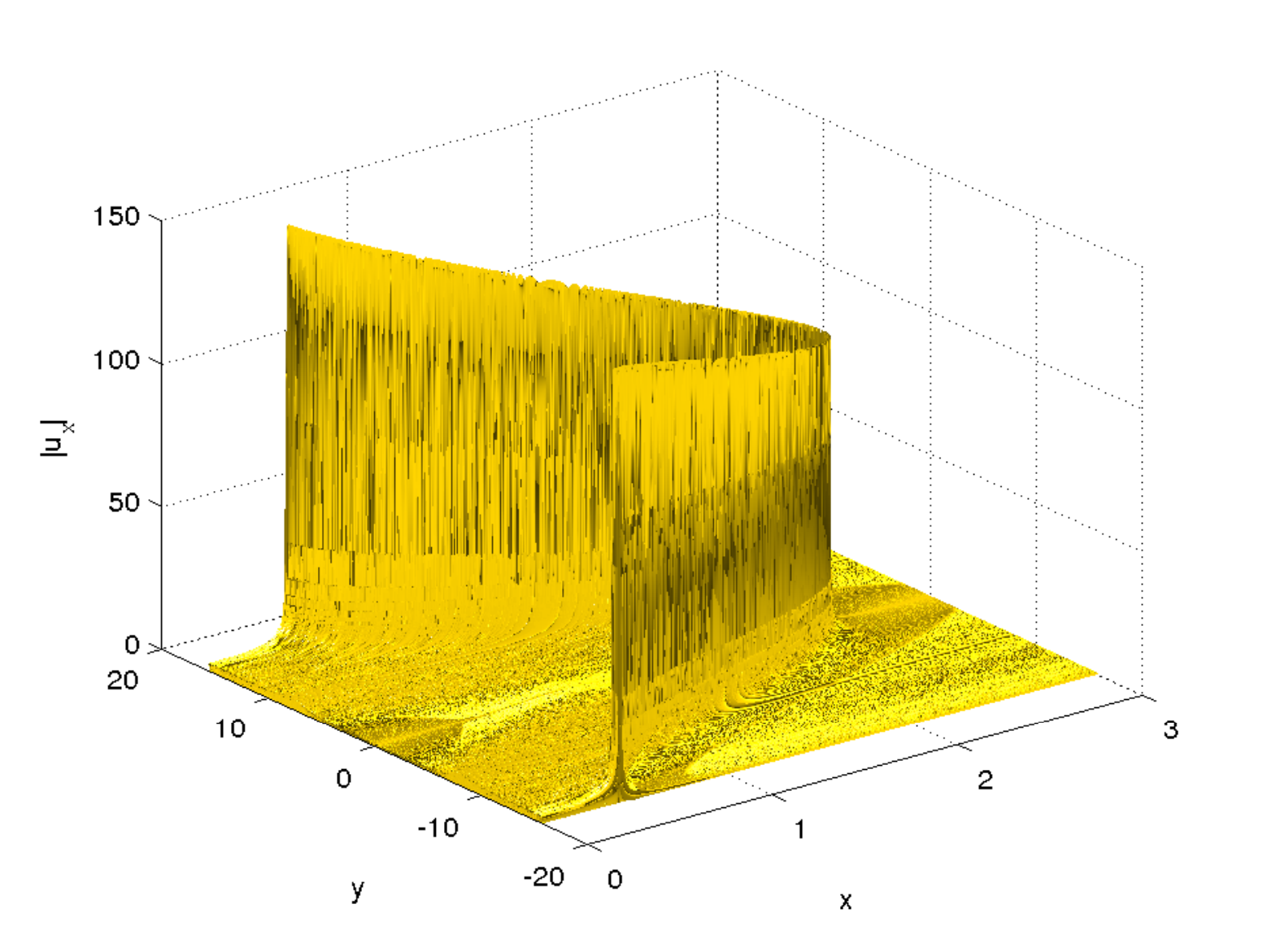} 
  \includegraphics[width=0.45\textwidth]{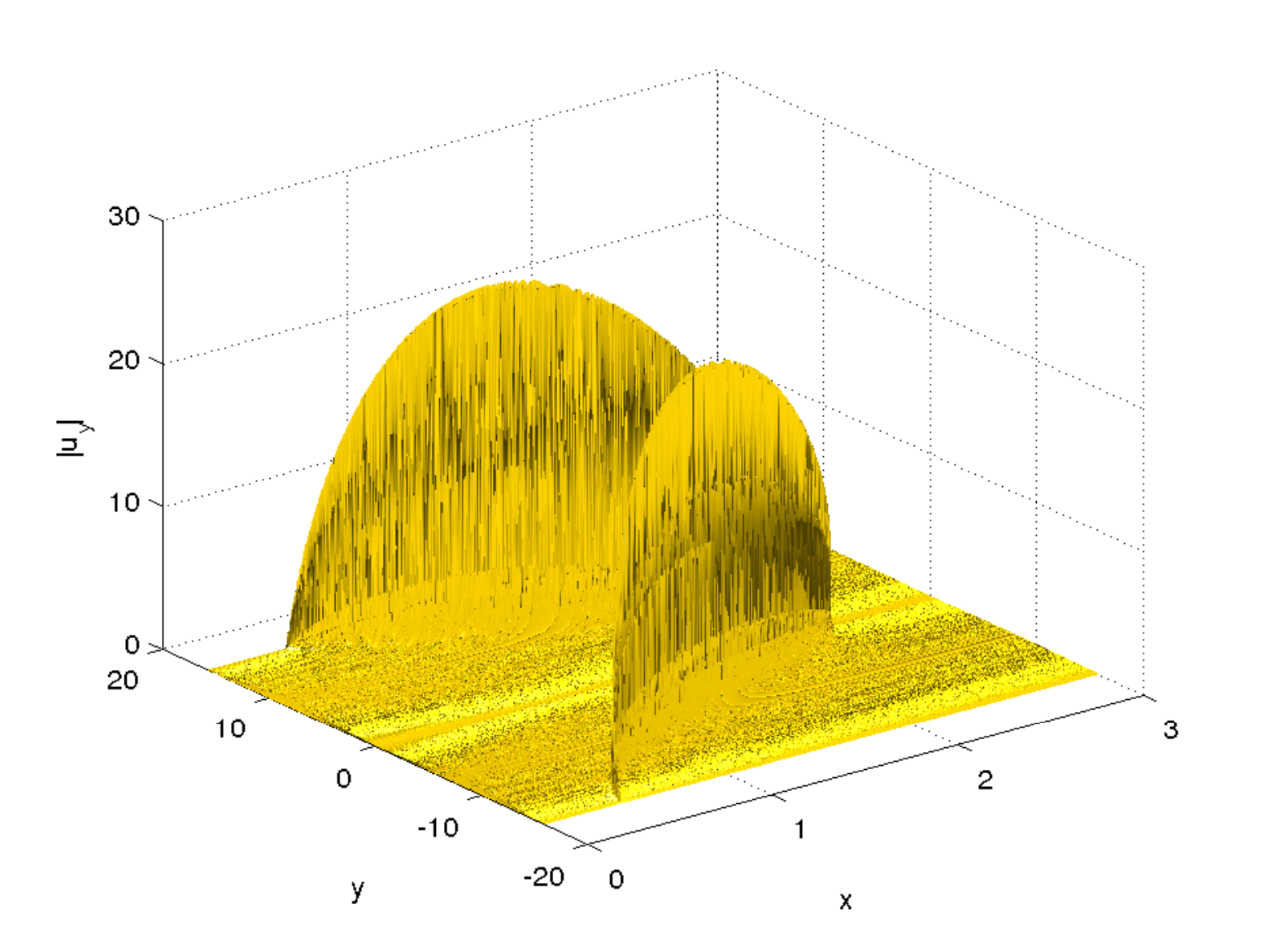}
 \caption{Derivatives  of the dKP II solution for the initial data 
 (\ref{uini2}) at the critical time $t_{c}=0.1934$, $|\partial_x u|$ (left) and of $|\partial_y u|$.}
 \label{uxuytcdkp2exp}
\end{figure}

A noticeable difference appears however in the profile of $u_x$. One 
can observe in Fig. \ref{uxuytcdkp2exp}, that in contrast to the dKP I case, the strongest value of the 
$x$-gradient here occur at the trailing part of the solution, whereas 
in the case of dKP I, it was occurring at roughly $y=0$. This is 
related to the focusing and defocusing character of dKP I and dKP II 
respectively. In the first case, the most advanced point in 
propagation direction gets focused, whereas in the latter case, this 
happens for the trailing points. 

\subsection{Shock formation in dKP solutions for 
initial data localized in both spatial dimensions}
In this subsection we will study shock formation in dKP solutions for 
initial data of the form (\ref{uini1}) which are localized in both 
spatial directions. It is known, see 
\cite{alinhac,MS08}, that localized data with a single 
maximum will lead to a break-up in a single point as expected. 

\begin{remark}\label{remark}
    In order to satisfy the constraint (\ref{const}), 
    the data studied here do not have a single 
maximum. Therefore there will 
be actually two points of gradient catastrophe.  
The second gradient catastrophe appears for dKP I for
negative $x$ for a slightly larger $t$ (for dKP II the role of the 
two critical points is interchanged). Since the code can only be 
reliably run until the first shock formation, we cannot obtain the 
second break-up with the method, and we will 
only discuss this time in the following. 
\end{remark}

\begin{figure}[htb!]
\centering
\includegraphics[width=\textwidth]{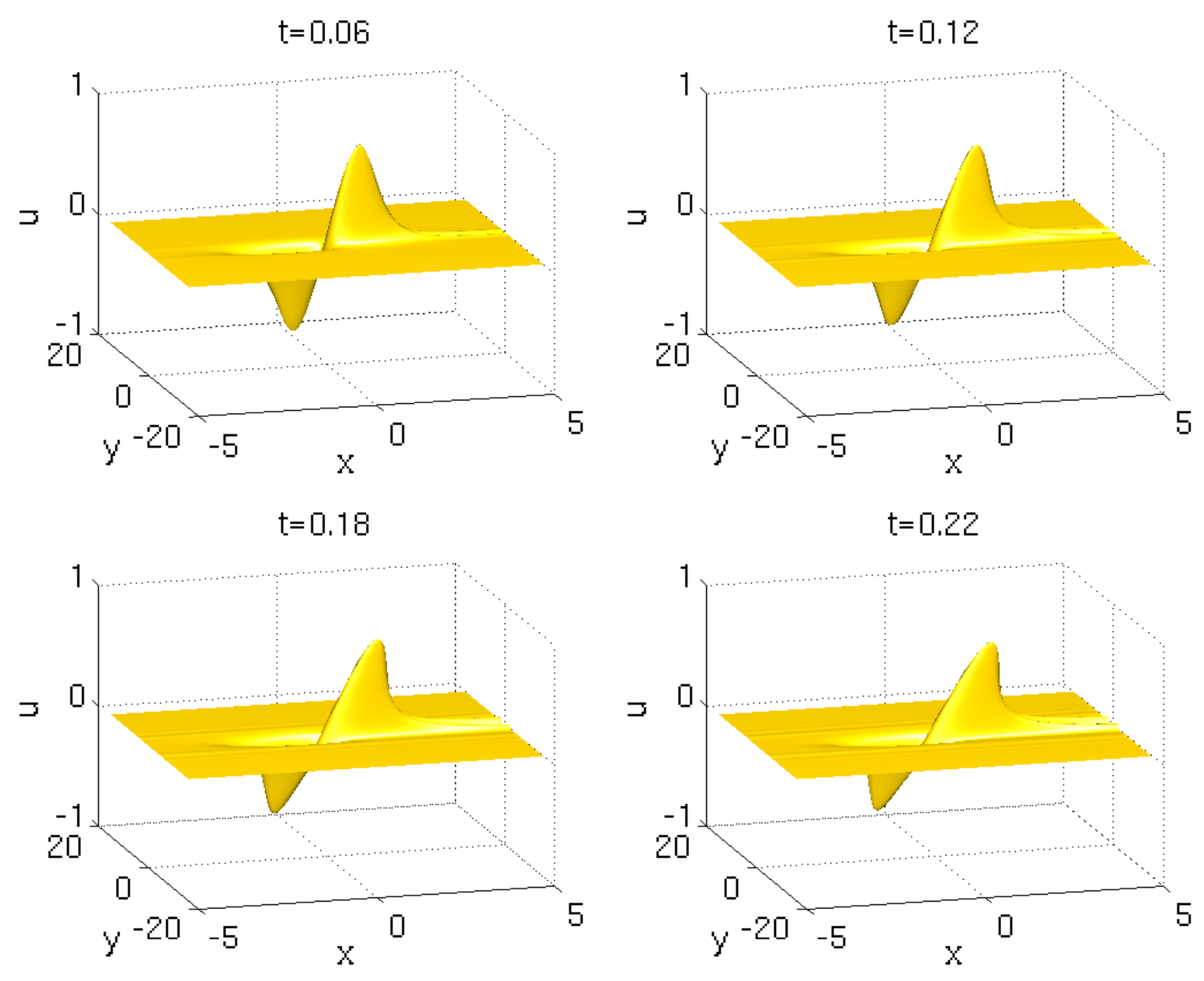} 
\caption{Solution to the dKP I equation for initial data 
(\ref{uini1}), for several values of t.}
\label{u2uts2d}
\end{figure}

\begin{figure}[htb!]
\centering
  \includegraphics[width=0.45\textwidth]{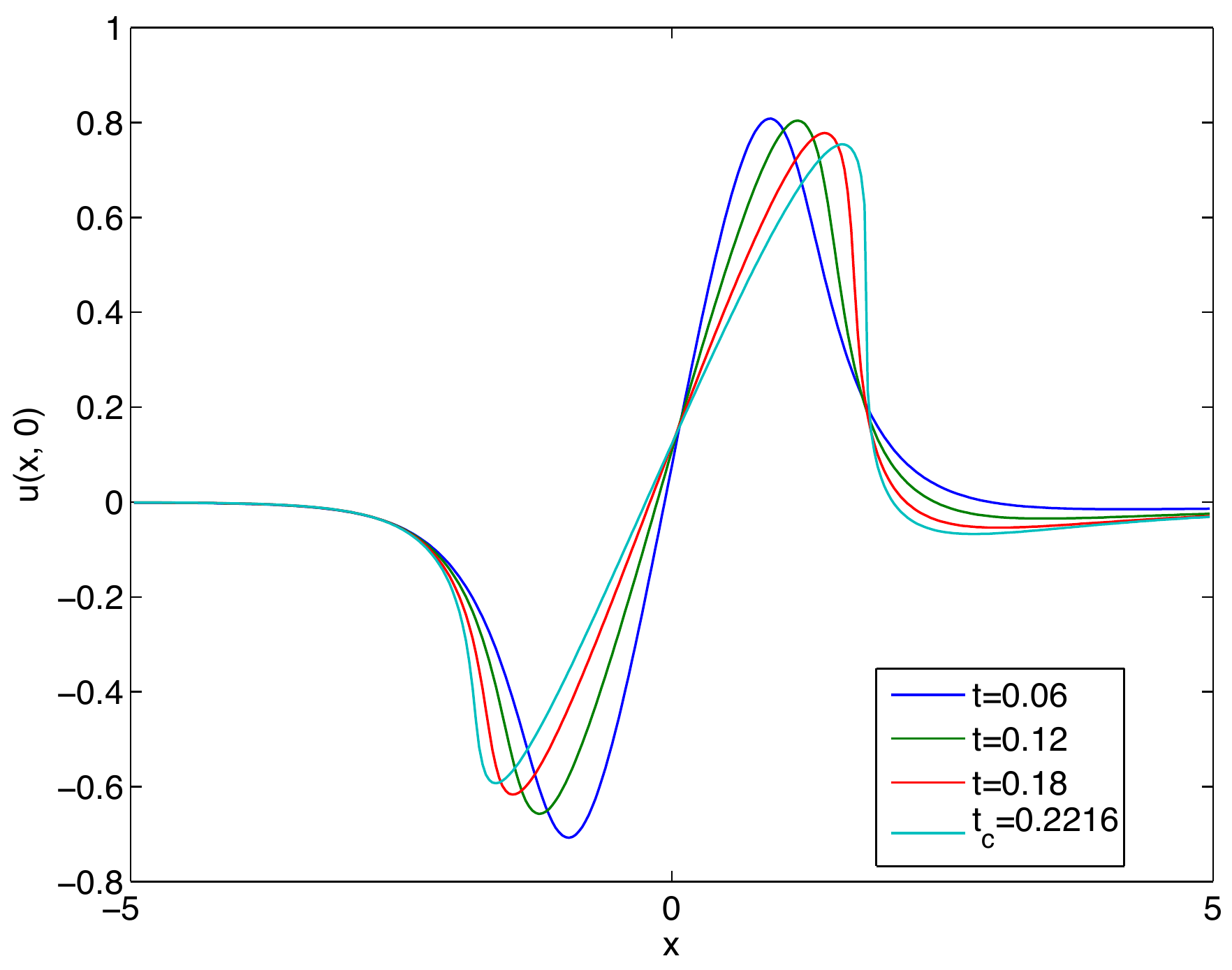}  
  \includegraphics[width=0.45\textwidth]{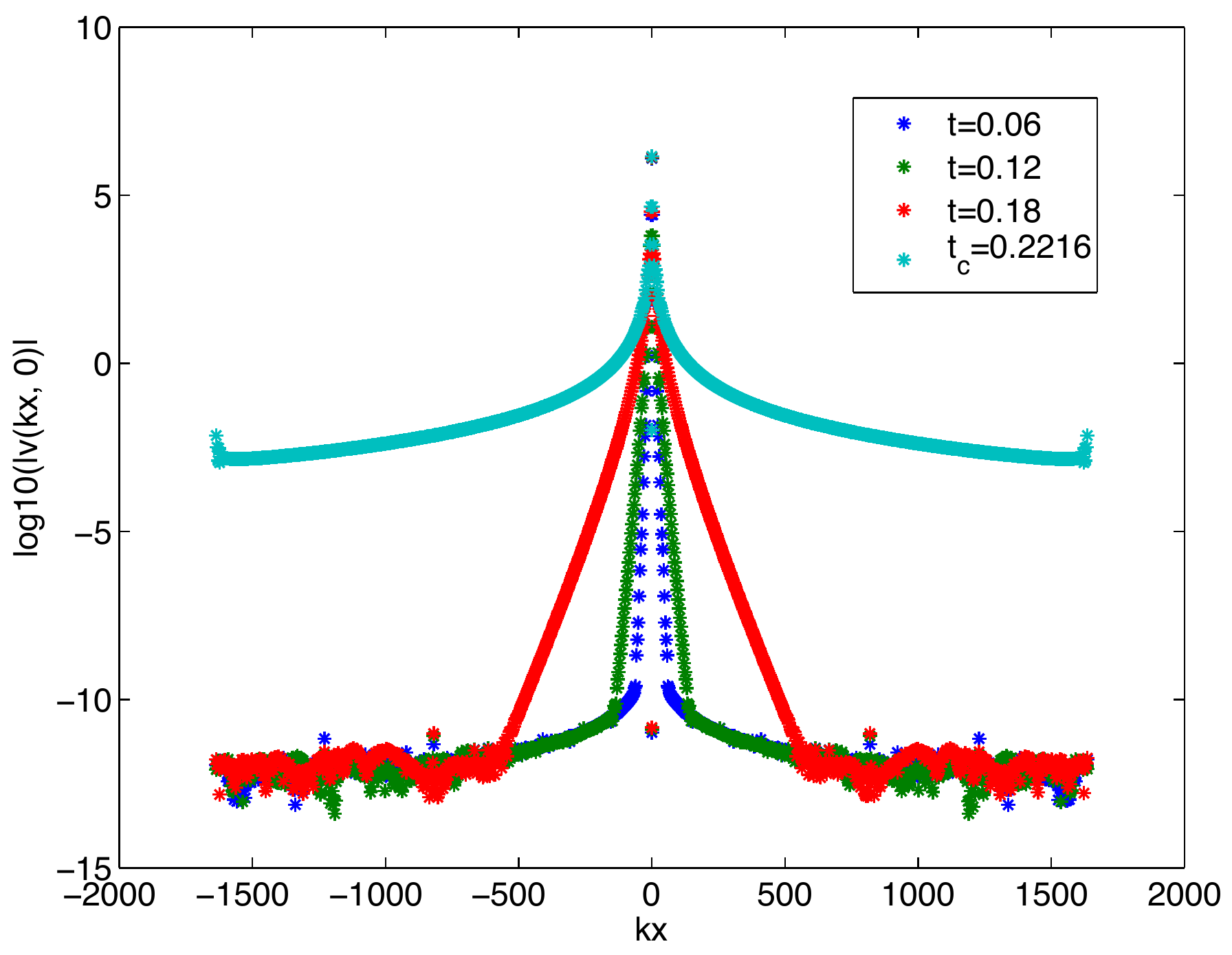}
 \caption{Solution to the dKP I equation for initial data 
 (\ref{uini1}), plotted on the $x$-axis for several values of $t$, 
 and the corresponding 
 Fourier coefficients, plotted on the $k_x$-axis.}
 \label{u2uts}
\end{figure}

For $N_x=N_y=2^{14}$ and $L_x=L_y=5$, we observe that the solution develops a shock in the $x$-direction at a time $t \sim 0.22$, as we can see in Fig. \ref{u2uts2d}.
This is in accordance with the experiments in \cite{KSM} for the same 
initial data, where the break-up was studied qualitatively for a 
dissipative regularization of the dKP equation. There it was found 
that the gradient catastrophe for $x > 0$ is reached roughly at $t 
\sim 0.23$. In Fig.~\ref{u2uts2d} it can also be seen that the 
initially rapidly decreasing function $u$ develops tails with an 
algebraic fall-off to infinity. Due to the imposed periodicity 
condition, this leads to a Gibbs phenomenon at the boundary of the 
computational domain. To reduce the effect of this phenomenon on the 
Fourier coefficients, we have used a considerably larger 
computational domain than shown in Fig.~\ref{u2uts2d}. We present 
in Fig. \ref{u2uts} the solution to the dKP I equation with initial data (\ref{uini1}), plotted on the $x$-axis for several values of $t$ 
 and the corresponding Fourier coefficients, plotted on the 
 $k_x$-axis. From the latter it can be inferred that the asymptotic 
 behavior of the Fourier coefficients is dominated by the break-up 
 singularity and not by the algebraic 
 fall off of the solution.

\begin{figure}[htb!]
\centering
  \includegraphics[width=0.45\textwidth]{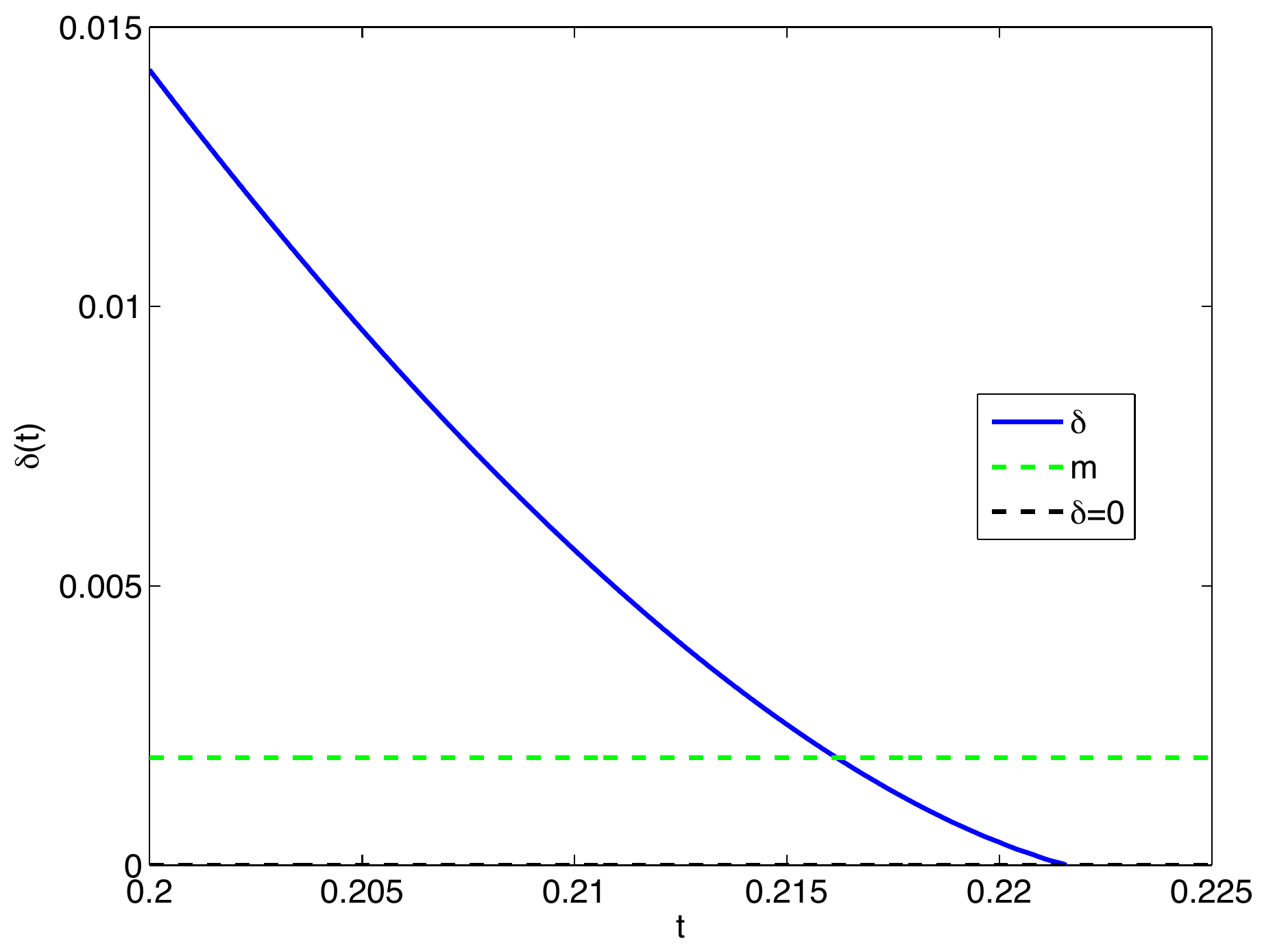}
  \includegraphics[width=0.45\textwidth]{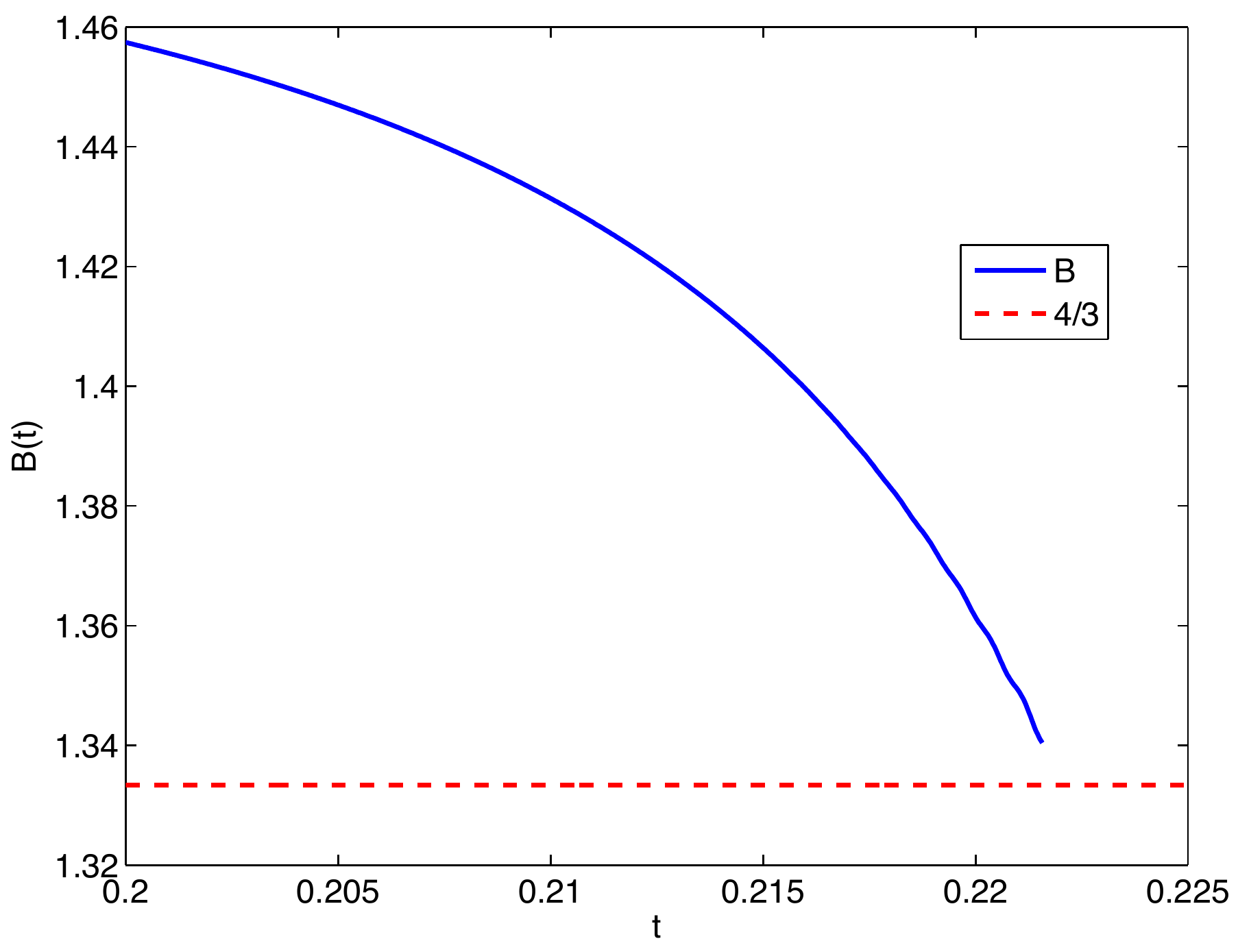} 
 \caption{Time dependence of the fitting parameters $\delta$ (left), 
 $B$ (right) 
 for the numerical solution to the dKP I 
 equation for  initial data \ref{uini1}.
The fitting  is done for  $10<k_x<\max(k_x)/2$. }
 \label{u2f1}
\end{figure}

In this case, we find that the minimal attainable fitting error (with 
the above described requirement to use at least half of the Fourier 
coefficients) is $p\sim0.5$ for $k_{min}=10$ and $k_{max}=max(k_x)/2$.
For such bounds, $\delta$ vanishes at 
$t_c=0.2216$, where $B(t_c)=1.34$ and  $\Delta= \|   \ln |v(k_x, 0)| 
- (A - B \ln k_x - k_x \delta)   \|_{\infty}=0.48$.  The considerably 
larger error in the fitting than for the initial data (\ref{uini2}) 
is related to the appearence of a second 
break-up singularity for the data (\ref{uini1}) as explained in 
Remark~\ref{remark}. Nonetheless the results of the fitting are 
reliable as can be seen in the following. 
A fitting of $\Im\log(v)$ yields 
$x_c=\alpha(t_c)= 1.7961$. The value of $B$ is compatible with the 
expected value 
$\frac{4}{3}$. 
The time dependence of the fitting parameters is shown in Fig. \ref{u2f1} 
for $t\sim t_c$. 

\begin{figure}[htb!]
\centering
  \includegraphics[width=0.45\textwidth]{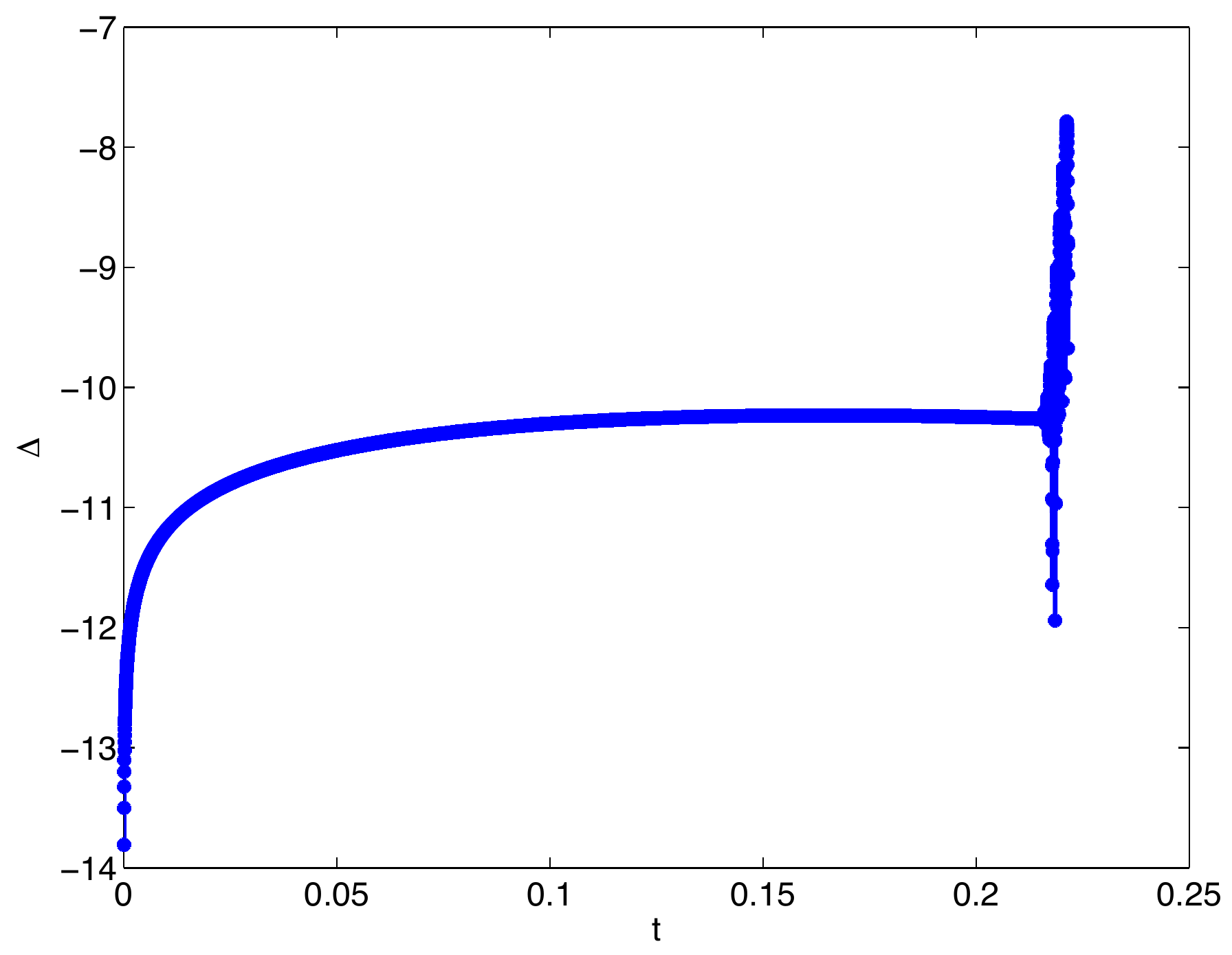} 
\includegraphics[width=0.45\textwidth]{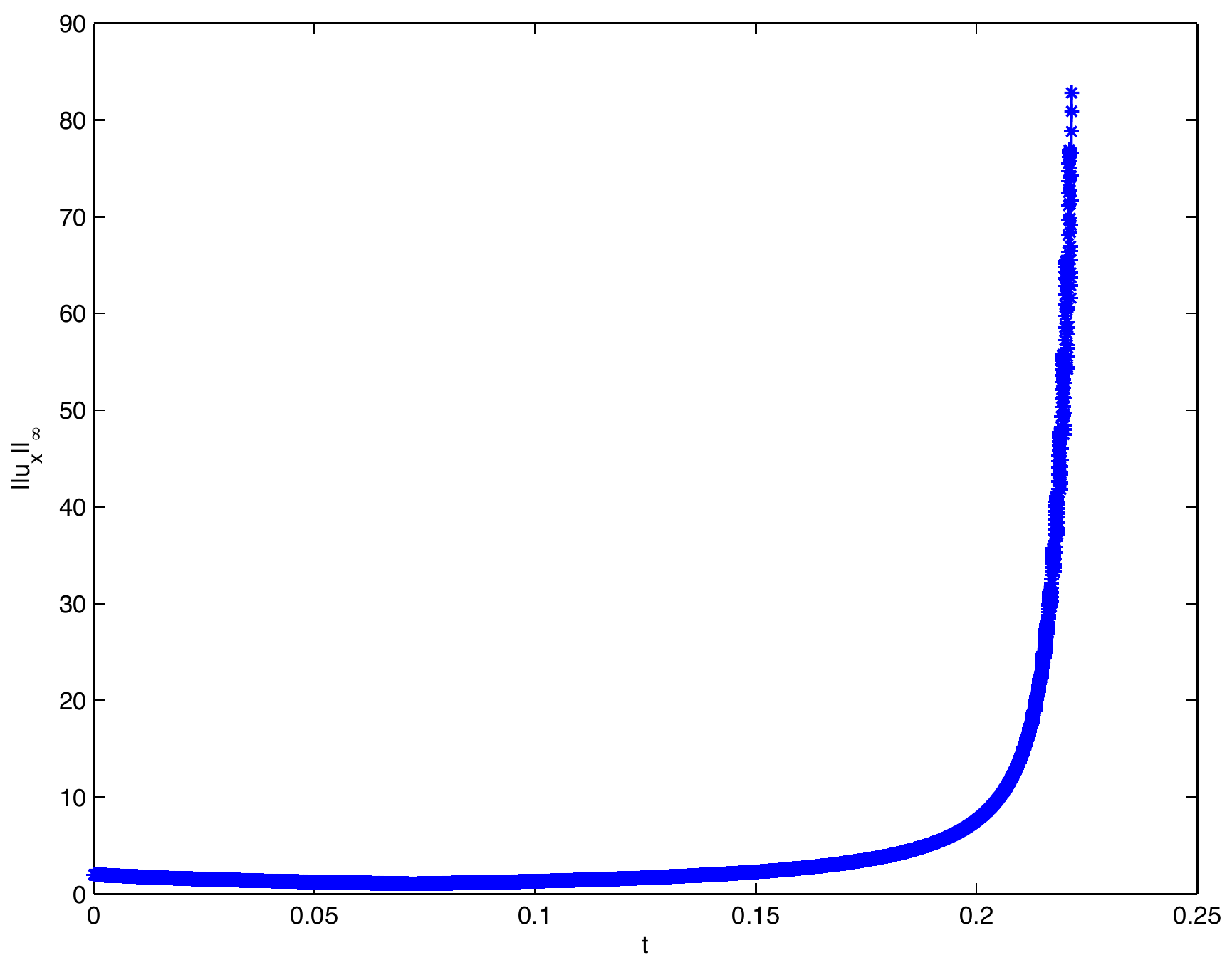} 
 \caption{Time evolution of the quantity $\Delta_{E}$(\ref{deltaE}) 
 (left) and of  $\|u_x\|_{\infty}$ (right) 
 for the solution to the dKPI equation for the initial data (\ref{uini1}).}
 \label{massamplu2}
\end{figure}

The time dependence of the quantity $\Delta_E$ 
(\ref{deltaE}) controlling the numerical error,  and of  $\|u_x\|_{\infty}$ until $t_c=0.2216$ are shown in Fig. \ref{massamplu2}.

\begin{figure}[htb!]
\centering
 \includegraphics[width=0.45\textwidth]{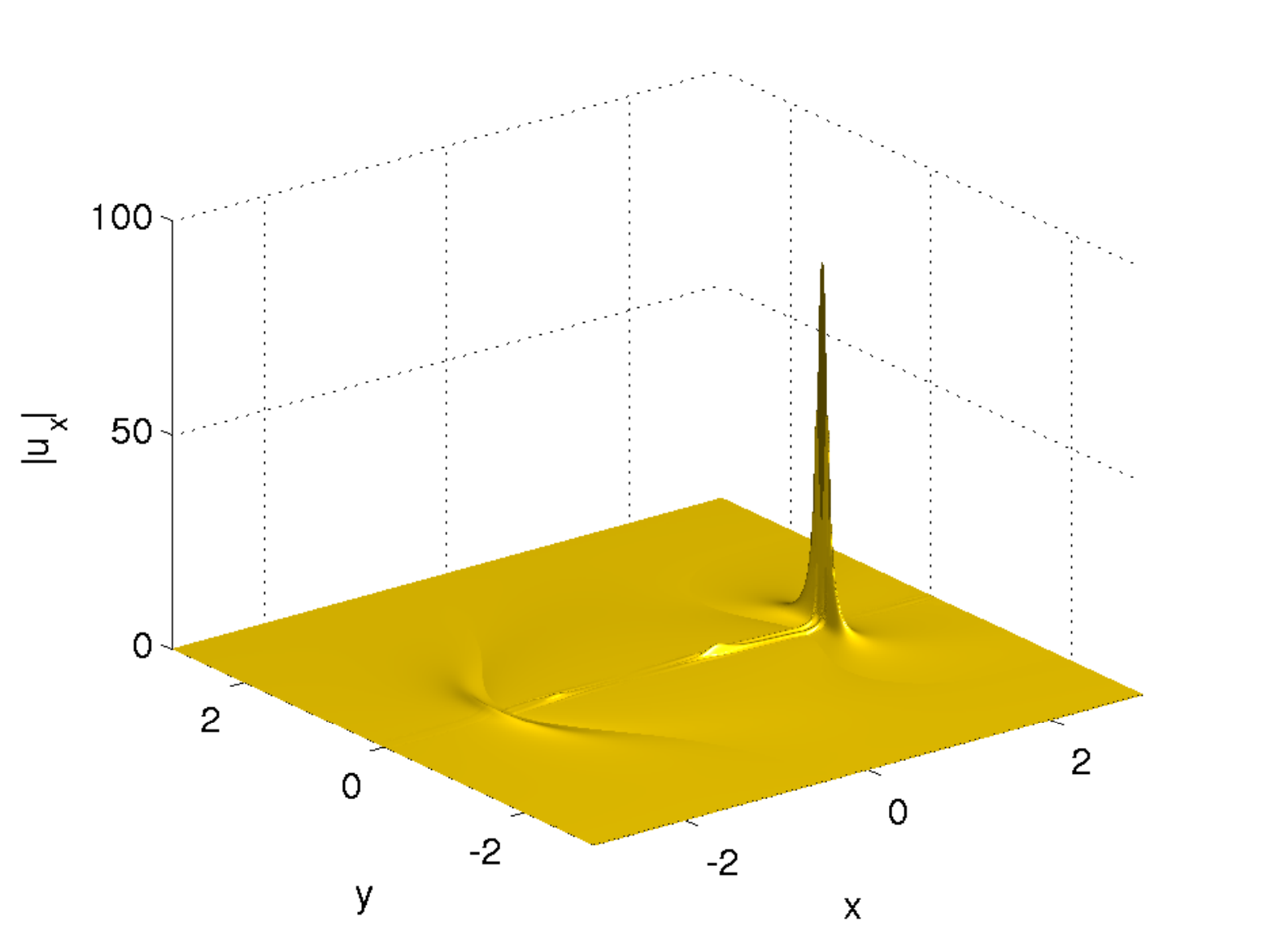} 
 \includegraphics[width=0.45\textwidth]{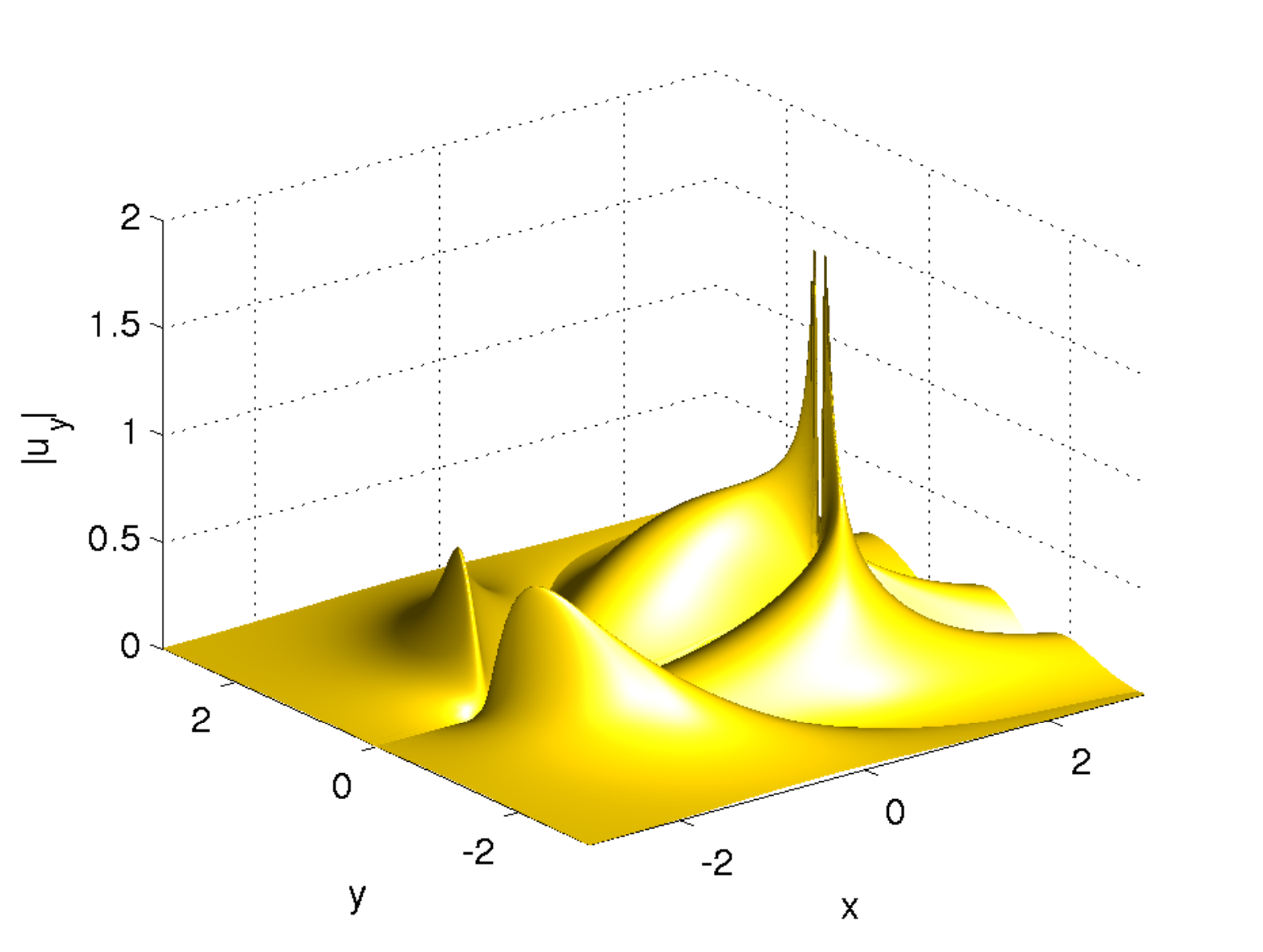}
 \caption{Derivatives  $|\partial_x u|$  and of $|\partial_y u|$ 
 at $t=0.2216$ for the solution to the dKP I equation for the initial data (\ref{uini1}).}
 \label{derivativesofu2}
\end{figure}

\begin{figure}[htb!]
\centering
 \includegraphics[width=0.45\textwidth]{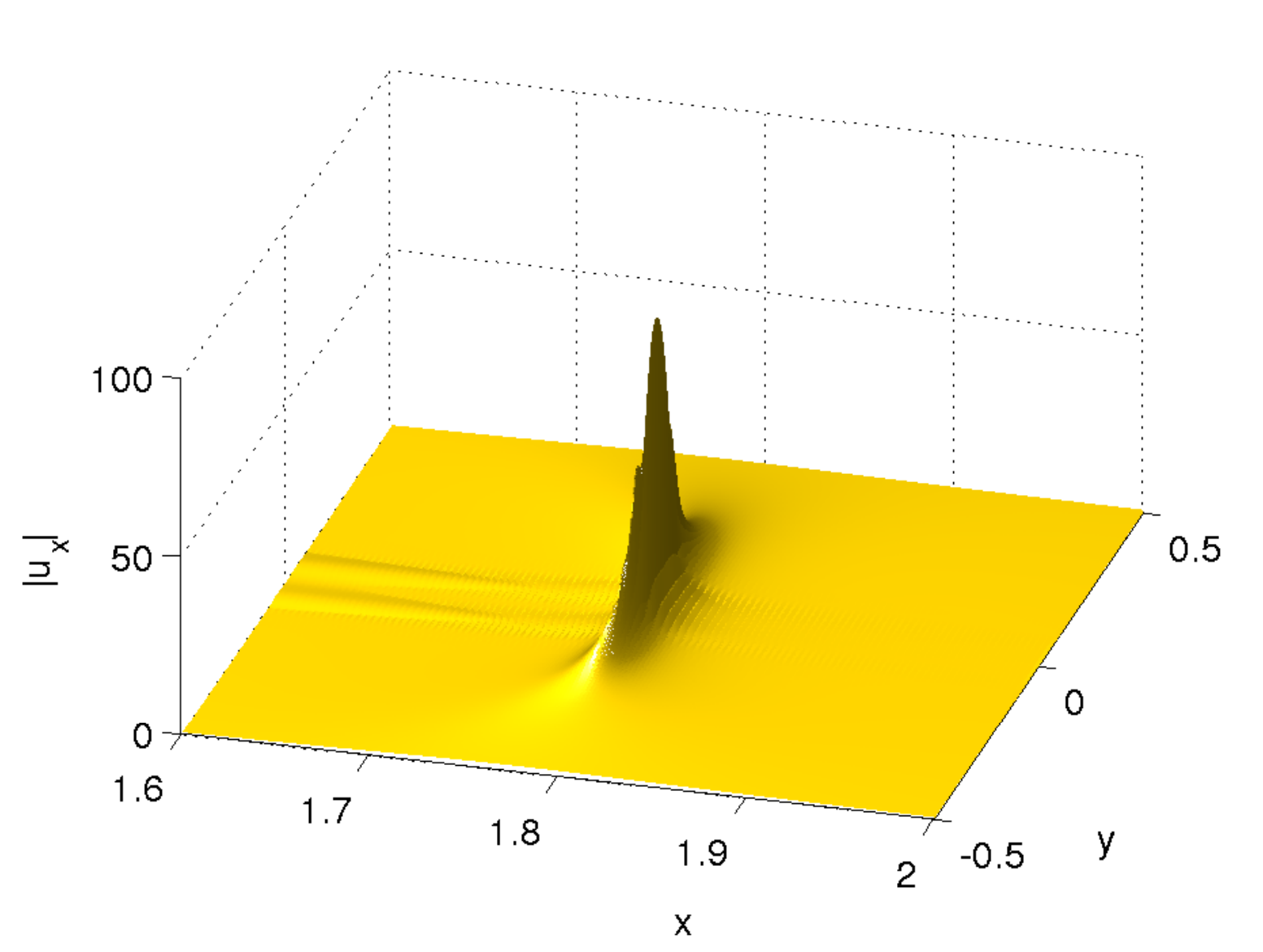} 
\includegraphics[width=0.45\textwidth]{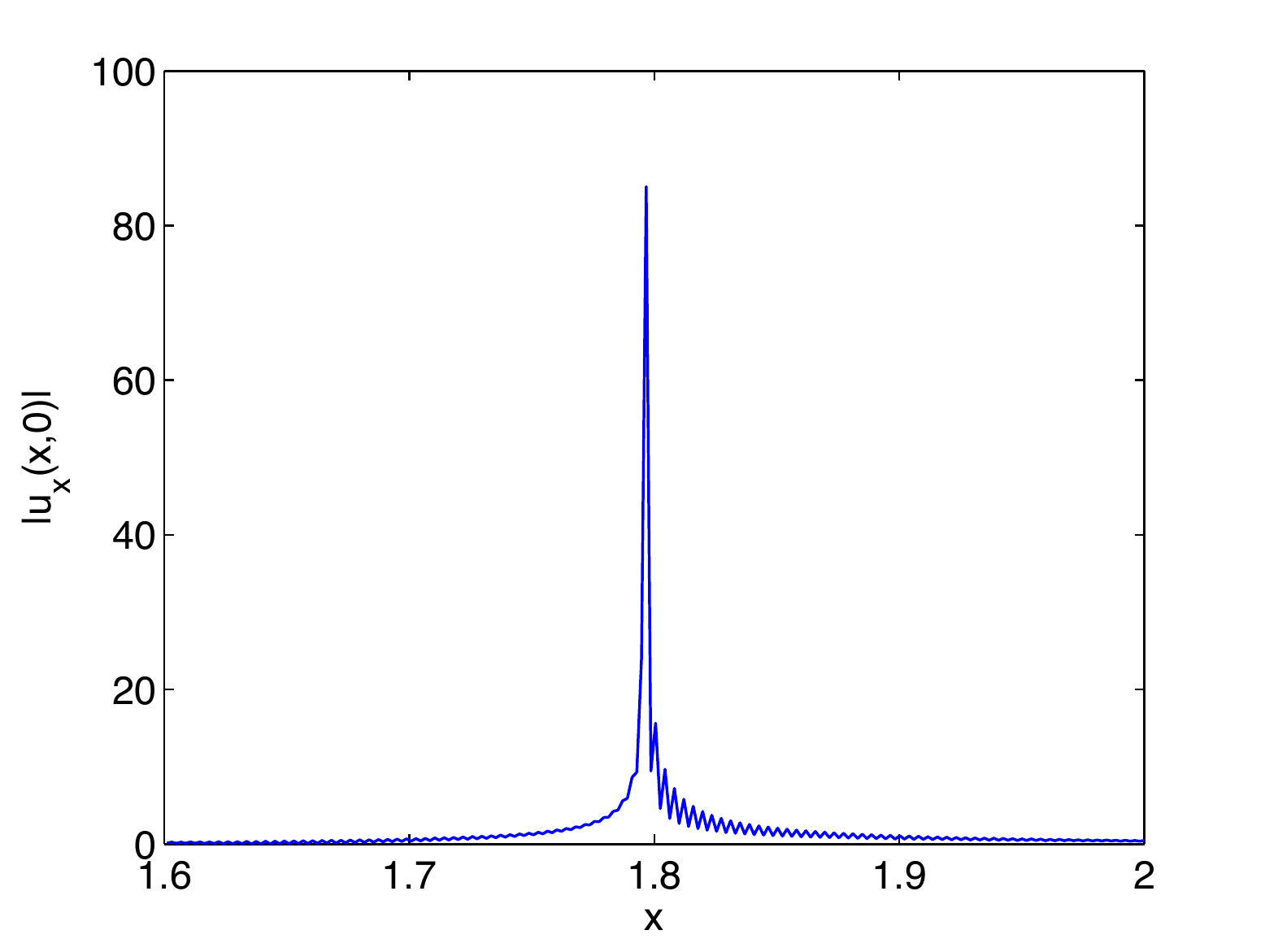}
 \caption{Close up of Fig.~\ref{derivativesofu2} of $|\partial_x 
 u|$ (left) and the profile of 
 $|\partial_x u(x,0)|$ (right) at $t=0.2216$.}
 \label{derzoom}
\end{figure}

In this case, the gradient catastrophe appears in only one point $(x_c, y_c)$, as we can see 
from Fig. \ref{derivativesofu2}, where we show the behavior of 
$|\partial_x u|$ (left) and of $|\partial_y u|$ (right) at 
$t=0.2216$. Note that the gradient in $y$-direction is much smaller 
than in $x$-directions which confirms the theoretical expectation 
that break-up only occurs in one direction whereas the solution remains 
regular in the other direction. 
This is even more obvious in a close up in Fig. \ref{derzoom} (top), 
where we show  $|\partial_x u|$ in the region where the gradient 
catastrophe occurs. We recover the previously computed value  
$\alpha(t_c)\sim 1.79$ in Fig. \ref{derzoom} (bottom), where 
$|\partial_x u(x,0)|$ is shown. It can be seen that  
$(x_c,y_c)=(1.79,0)$, and we find that $u_y(t_c,x_c,0)=4.8 * 10^{-13} \sim 0$, 
as also predicted by the theoretical results in \cite{MS12} for 
symmetry reasons.

\begin{figure}[htb!]
\centering
\includegraphics[width=0.5\textwidth]{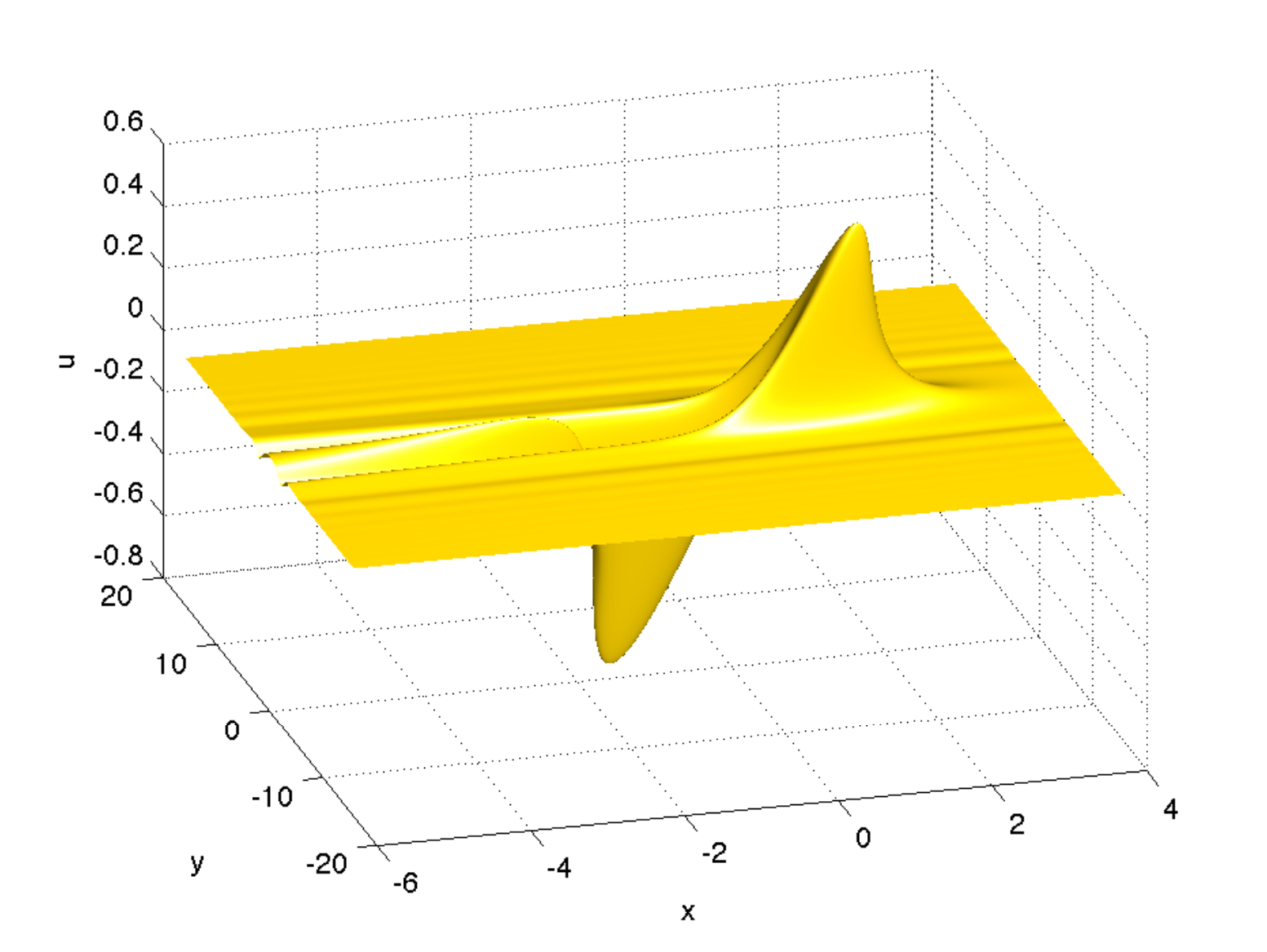} 
  \includegraphics[width=0.4\textwidth]{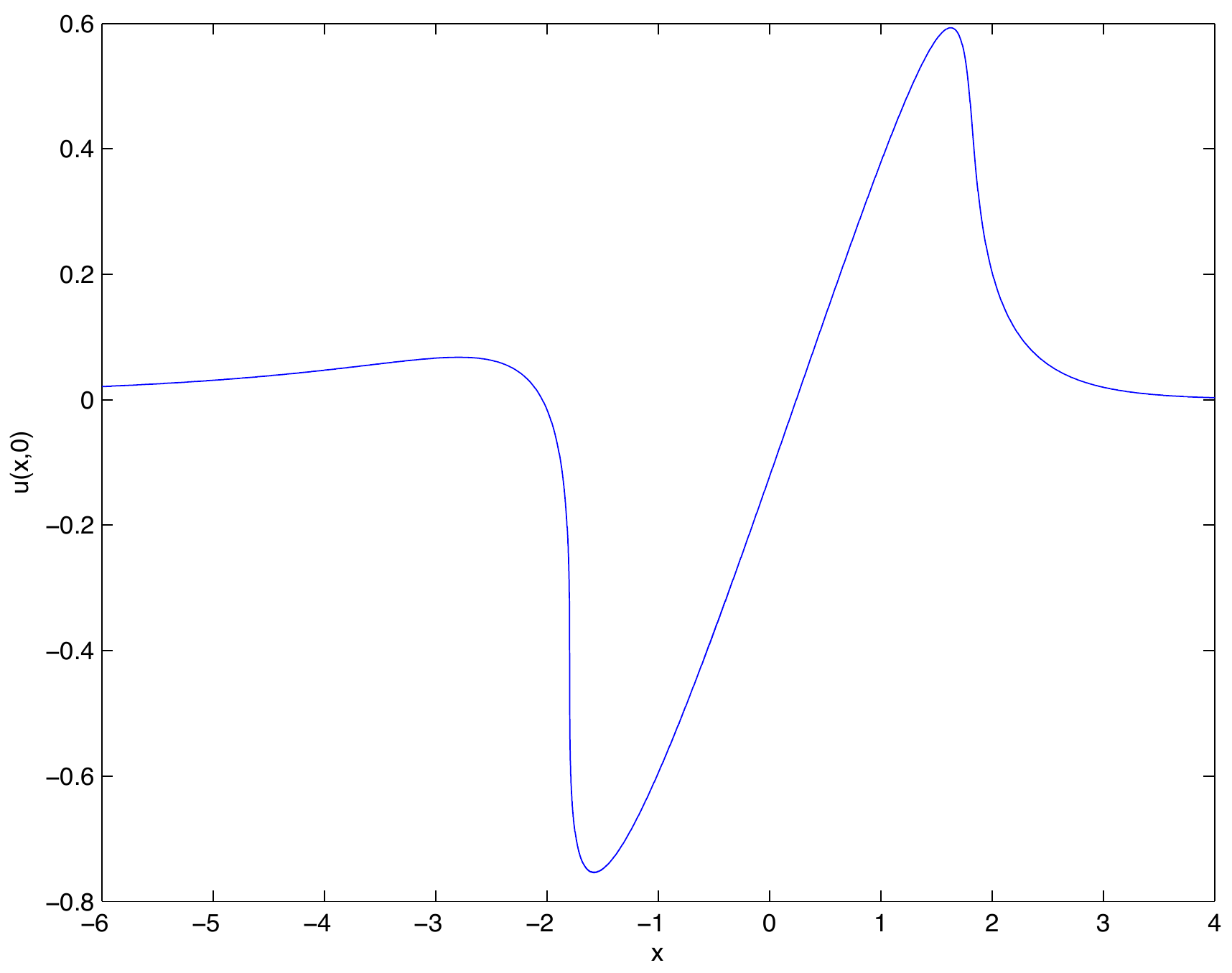}  
\caption{Solution to the dKP II equation with initial data 
(\ref{uini1}) at t=0.2216.}
\label{u2ddkp2}
\end{figure}

The situation is similar for the dKP II case for the initial data 
(\ref{uini1}), which we treat analogously. 
We find that the fitting used for dKP I is optimal also for dKP II
($k_{min}=10, k_{max}=\max(k)/2$), which yields the vanishing of $\delta$ at $t_c \sim 0.2216$, $B(t_c)=1.34$, and $\Delta= \|   \ln |v(k_x, 0)| - (A - B \ln k_x - k_x \delta)   \|_{\infty} = 0.48$.
We observe that the solution develops a shock in the $x$-direction 
($x<0$) at time $t_c$, see Fig.~\ref{u2ddkp2}, where we show the 
solution to the dKP II equation with initial data (\ref{uini1}), at 
$t=t_c$. Note that the tails with the algebraic fall off are now 
directed towards $-\infty$ in contrast to the dKP I case. Thus the 
gradient catastrophe appears in this case first for negative $x$ 
(again we only treat the point of gradient catastrophe appearing 
first). By doing a fitting of  $\Im\log(v)$, we find  $\alpha(t_c) 
\sim -1.79$. The good agreement of this value with the maximum of the 
gradient of $u$ shows the self consistency of the approach. 

\begin{figure}[htb!]
\centering
 \includegraphics[width=0.6\textwidth]{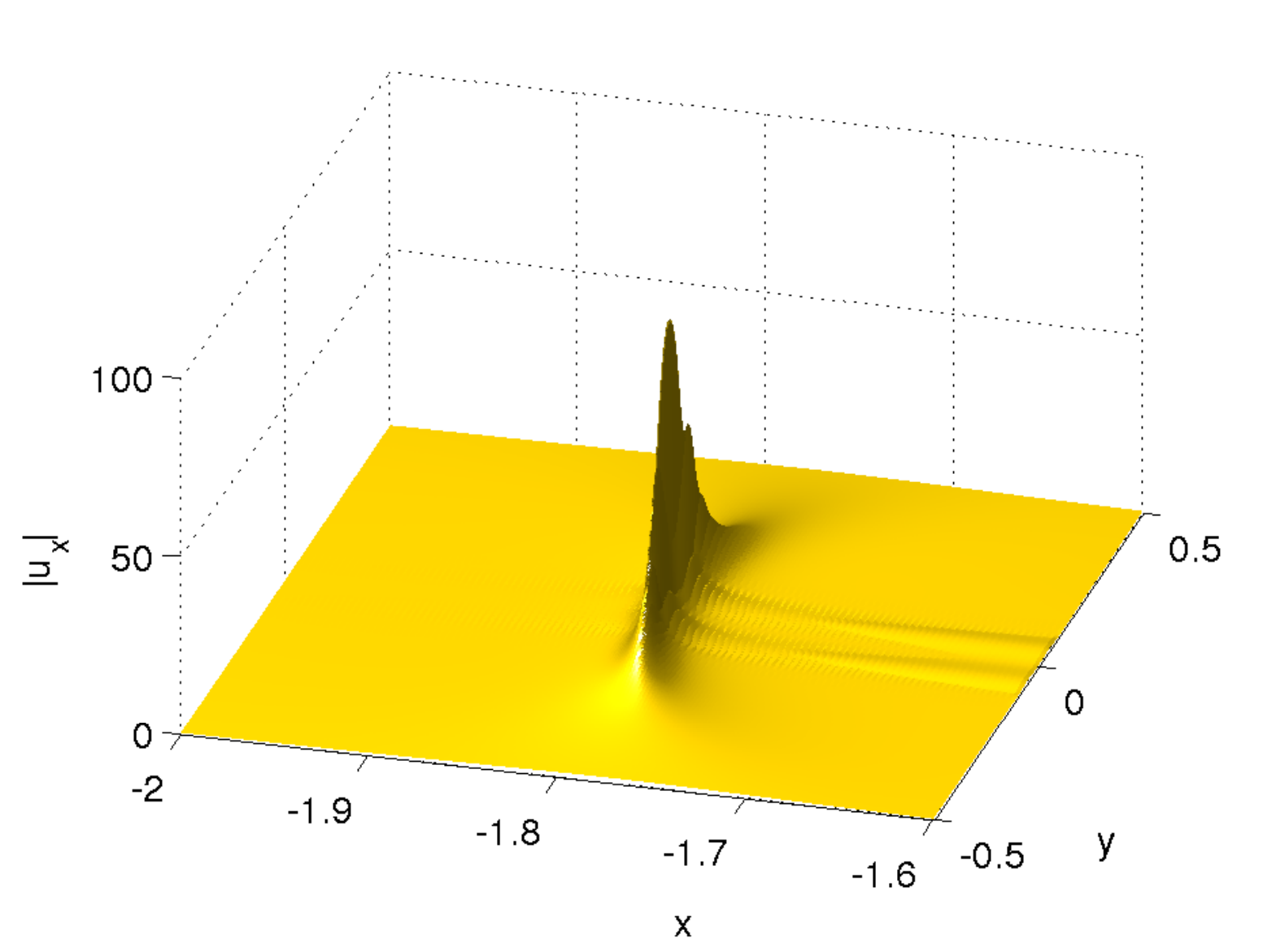} 
 \caption{Derivative  $|\partial_x u|$ at $t=0.2216$ of the 
 solution in Fig.~\ref{u2ddkp2}.}
 \label{derivativesofu2kp2}
\end{figure}

In this case, the gradient catastrophe appears again only in one point $(x_c, y_c) \sim(-1.79,0)$, as we can see 
from Fig.~\ref{derivativesofu2kp2}, where we show  $|\partial_x u|$ at $t=0.2216$, and $u_y(t_c,x_c,0)=3.3 *10^{-12}$.
The solutions to both dKP equations show thus a similar behavior, 
independent of the sign of $\lambda$ in (\ref{kpevol}), except that of 
course, in the case of dKP I, the gradient catastrophe appears at  
$(x_c,y_c)=(1.79,0)$, whereas in the case of dKP II, it appears at  
$(x_c,y_c)=(-1.79,0)$ (as actually expected). As already mentioned, 
this is due to the tails with an algebraic fall off towards infinity 
which are for dKP I directed towards $+\infty$, whereas they are 
directed towards $-\infty$ for dKP II.

\subsection{KP  solution in the small dispersion limit for initial 
data localized in one spatial dimension}
Solutions to the KP equation will develop \emph{dispersive shocks}, 
zones of rapid modulated oscillations, in the vicinity of a shock of 
solutions to the dKP equation for the same initial data. These were 
studied numerically for the first time in 
\cite{KSM}, see also  \cite{KR}. For KdV it was found that the 
difference between Hopf and KdV solution scales as $\epsilon^{2}$ for 
$t\ll t_{c}$ and as $\epsilon^{2/7}$ for $t\sim t_{c}$. 
In this subsection, we establish numerically such scaling laws for 
the initial data localized in only one spatial direction.

\begin{figure}[htb!]
\centering
  \includegraphics[width=0.8\textwidth]{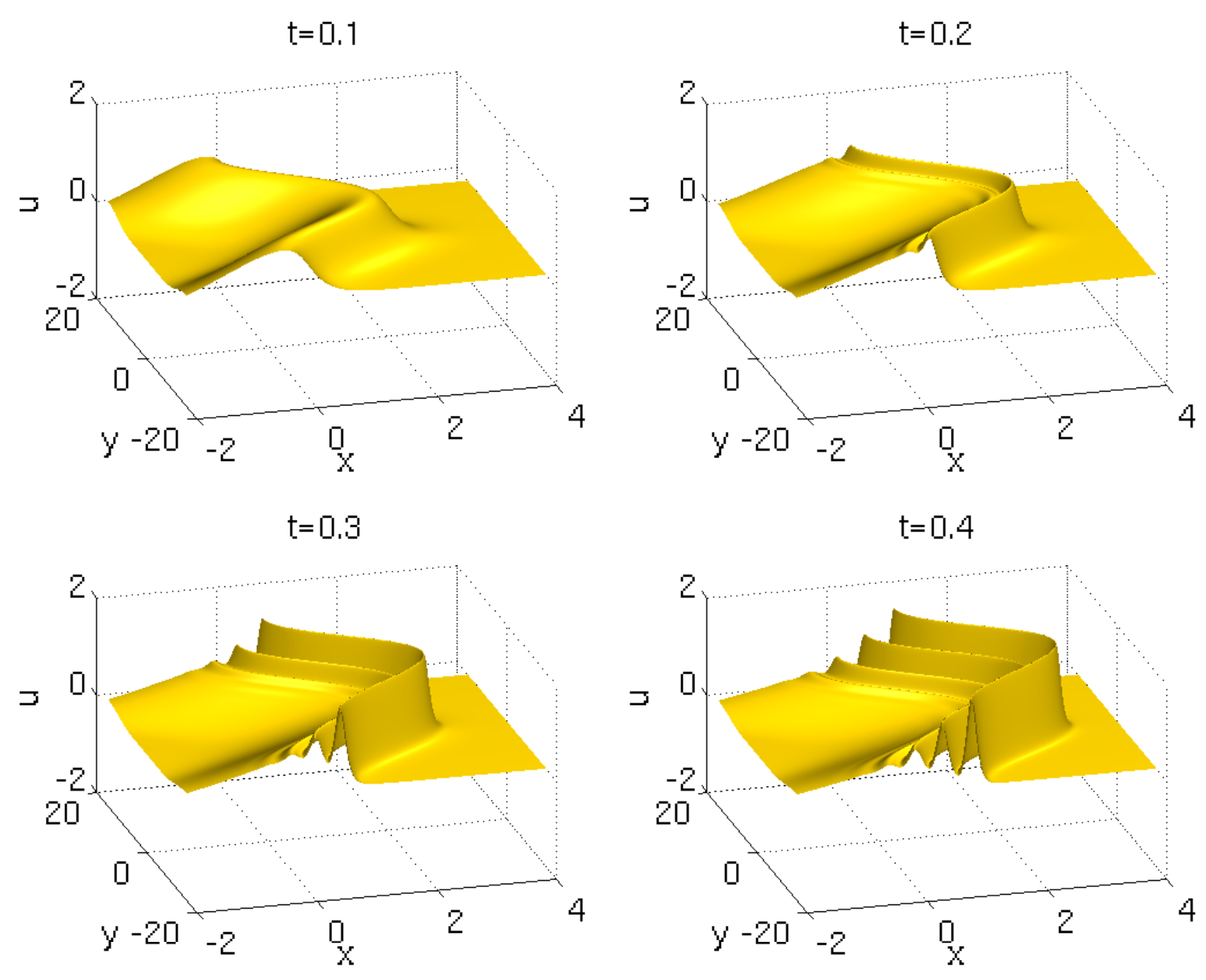} 
\caption{ Solution of the KP I equation in small dispersion limit for 
the initial data (\ref{uini2}) with $\epsilon=0.1$ for several values of $t$}
 \label{ut4s2d}
\end{figure}

For initial data of the form (\ref{uini2}), we show the 
solution to the KP I  equation in the small dispersion limit with $\epsilon=0.1$ in Fig. \ref{ut4s2d} for several times.
 The computations are carried out with $2^{14} \times 2^{14}$ points for $x \times y \in [-5\pi, 5\pi]$ and $\Delta_t =  4 * 10^{-5}$.
As expected, the dispersive regularization of the break-up singularity 
leads to rapid modulated oscillations in the region where the 
corresponding dKP I solution develops a shock.

\begin{figure}[htb!]
\centering
  \includegraphics[width=0.6\textwidth]{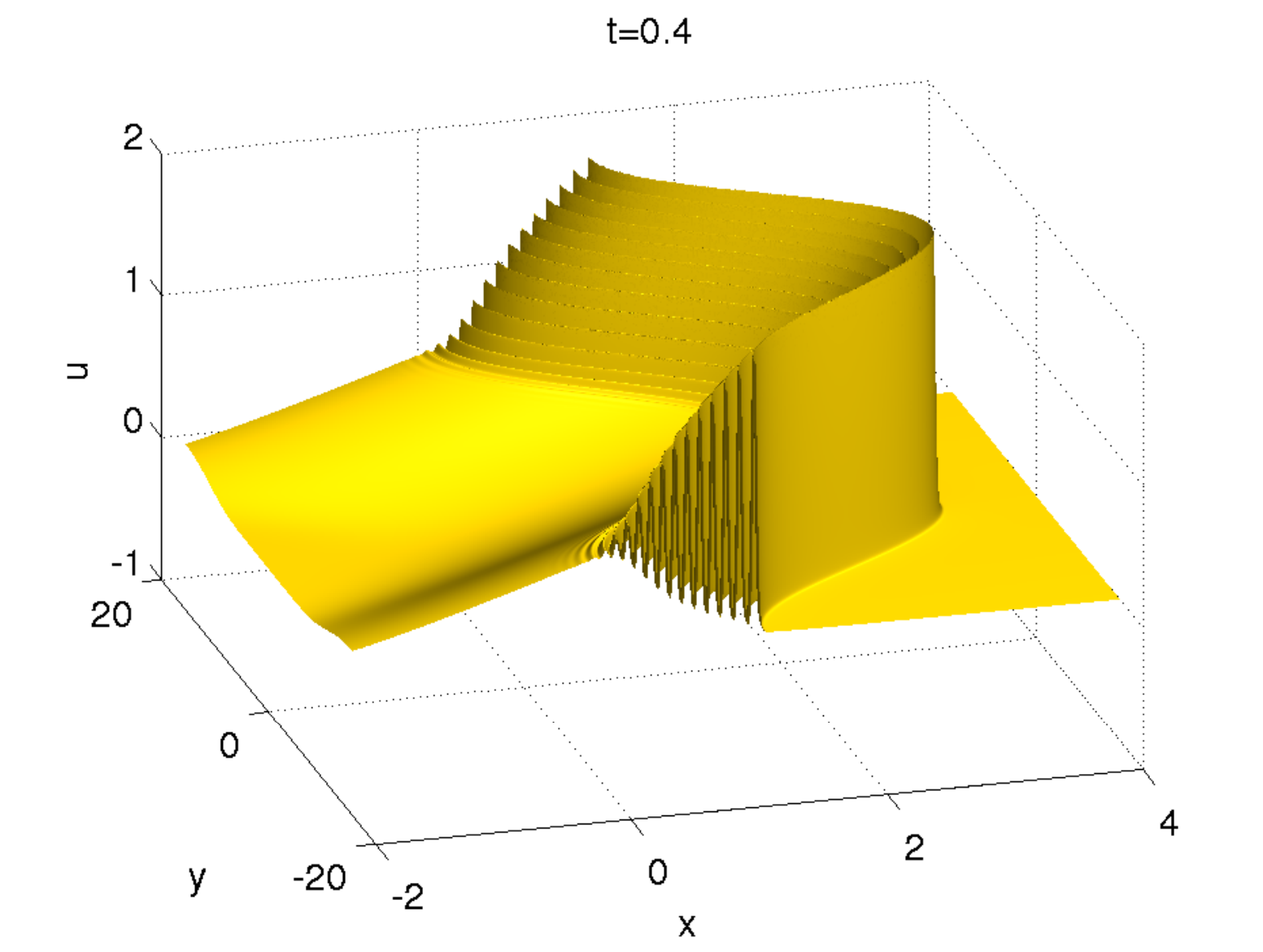} 
\caption{ Solution of the KP I equation in the small dispersion limit 
for initial data of the form (\ref{uini2}) with $\epsilon=0.02$ at t=0.4.}
 \label{ut4e0025}
\end{figure}

\begin{figure}[htb!]
\centering
  \includegraphics[width=0.8\textwidth]{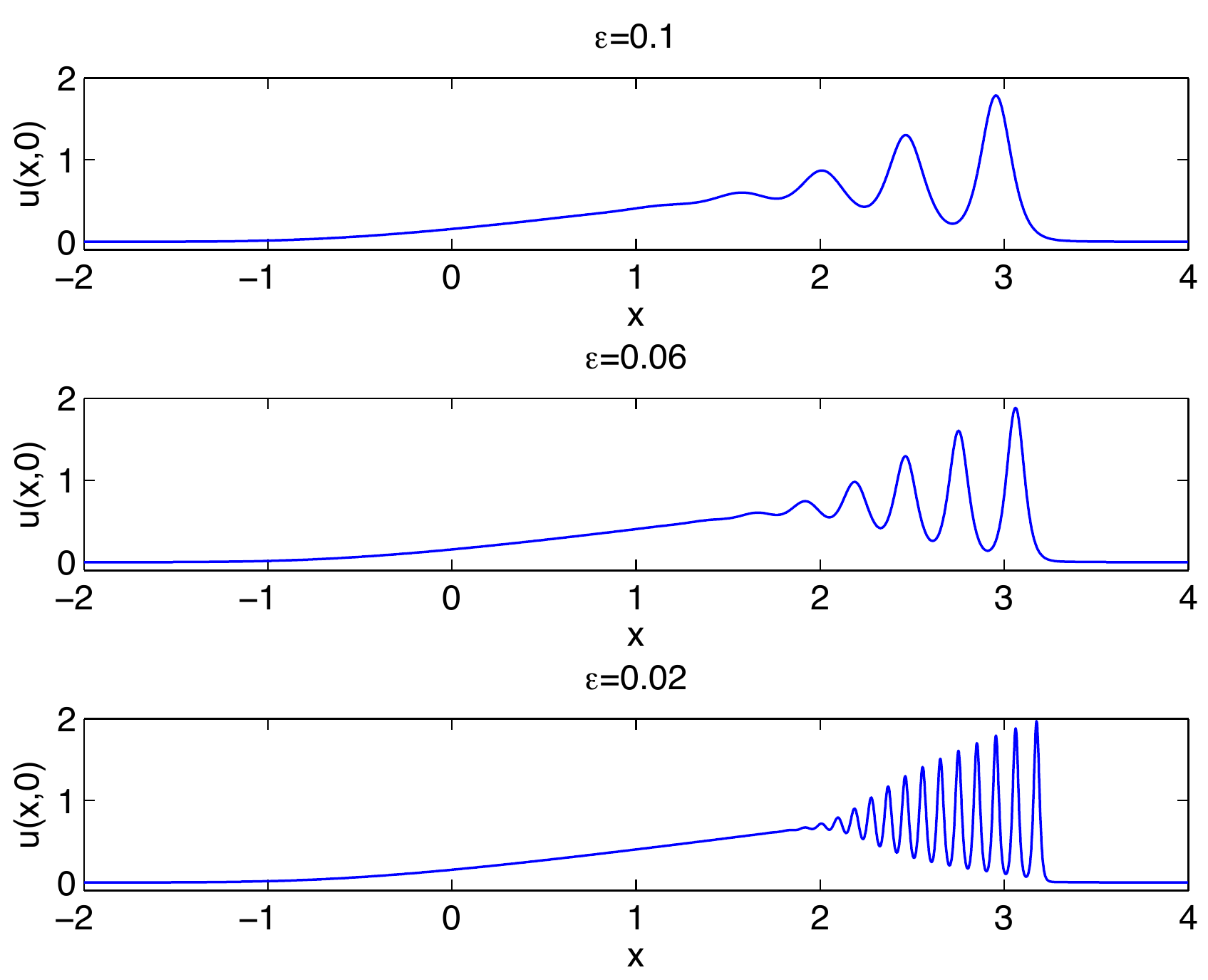} 
\caption{ Solutions to the KP I equation in the small dispersion limit on the $x$-axis at  $t=0.4$  for several values of $\epsilon$} 
 \label{ut4epss}
\end{figure}

\begin{figure}[htb!]
\centering
 \includegraphics[width=0.8\textwidth]{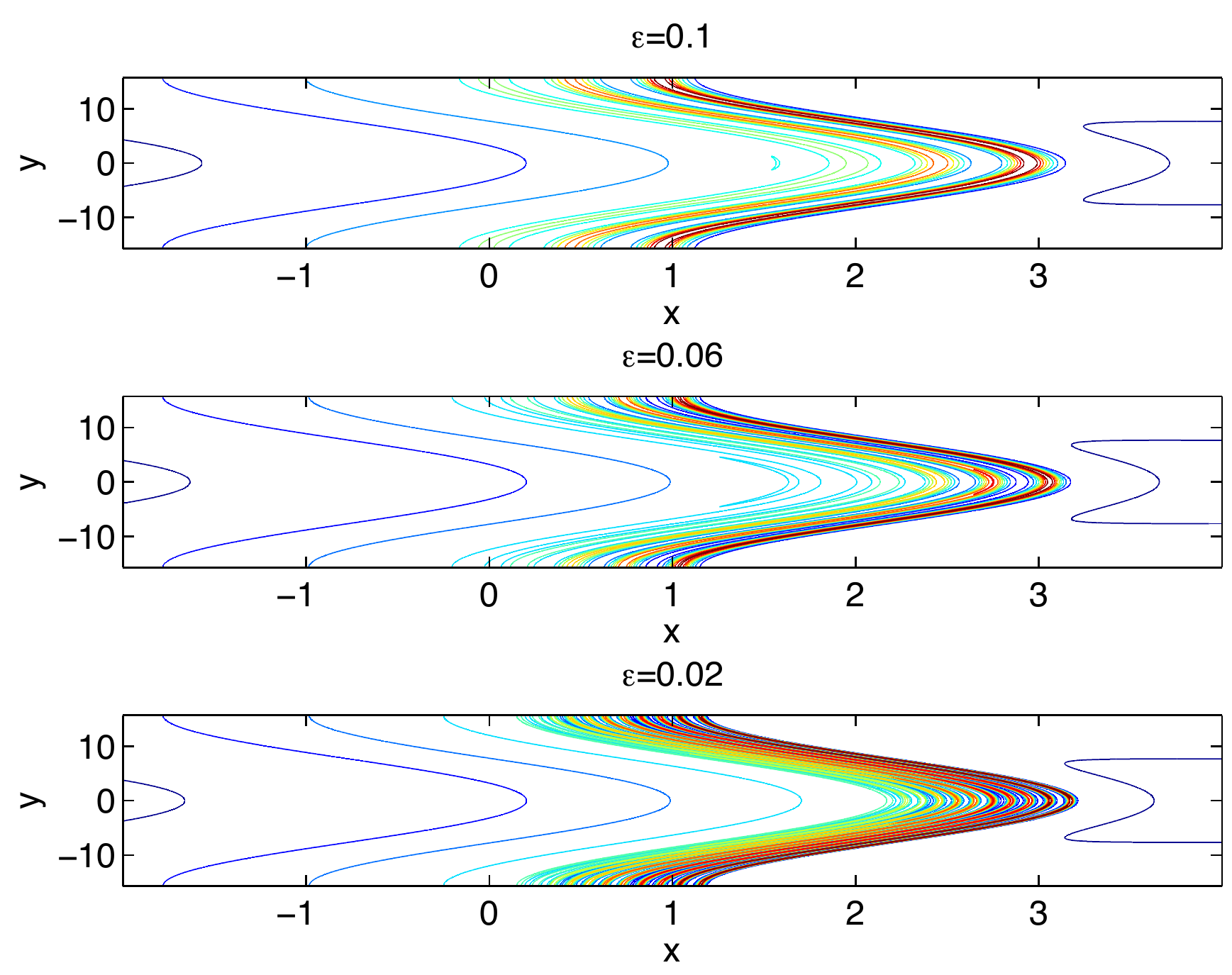}
\caption{Contour plots of the solutions to the KP I equation in the 
small dispersion limit at $t=0.4$  for several values of $\epsilon$.} 
 \label{cont1}
\end{figure}

The number of oscillations increases as $\epsilon$ tends to $0$, as 
already observed in \cite{KSM}. This can be seen in  
Fig.~\ref{ut4e0025} where 
the solution to the KP I equation in the small dispersion limit is 
shown for $\epsilon=0.02$ at $t=0.4$. 
In Fig.~\ref{ut4epss} we present 
the solutions to the KP I equation in the small dispersion limit for 
several values of $\epsilon$  on the $x$-axis for $t=0.4$.
\\
The contourplots  of the solutions at $t=0.4$ are shown in Fig.~\ref{cont1} for several values of $\epsilon$.

\begin{figure}[htb!]
\centering
  \includegraphics[width=0.4\textwidth]{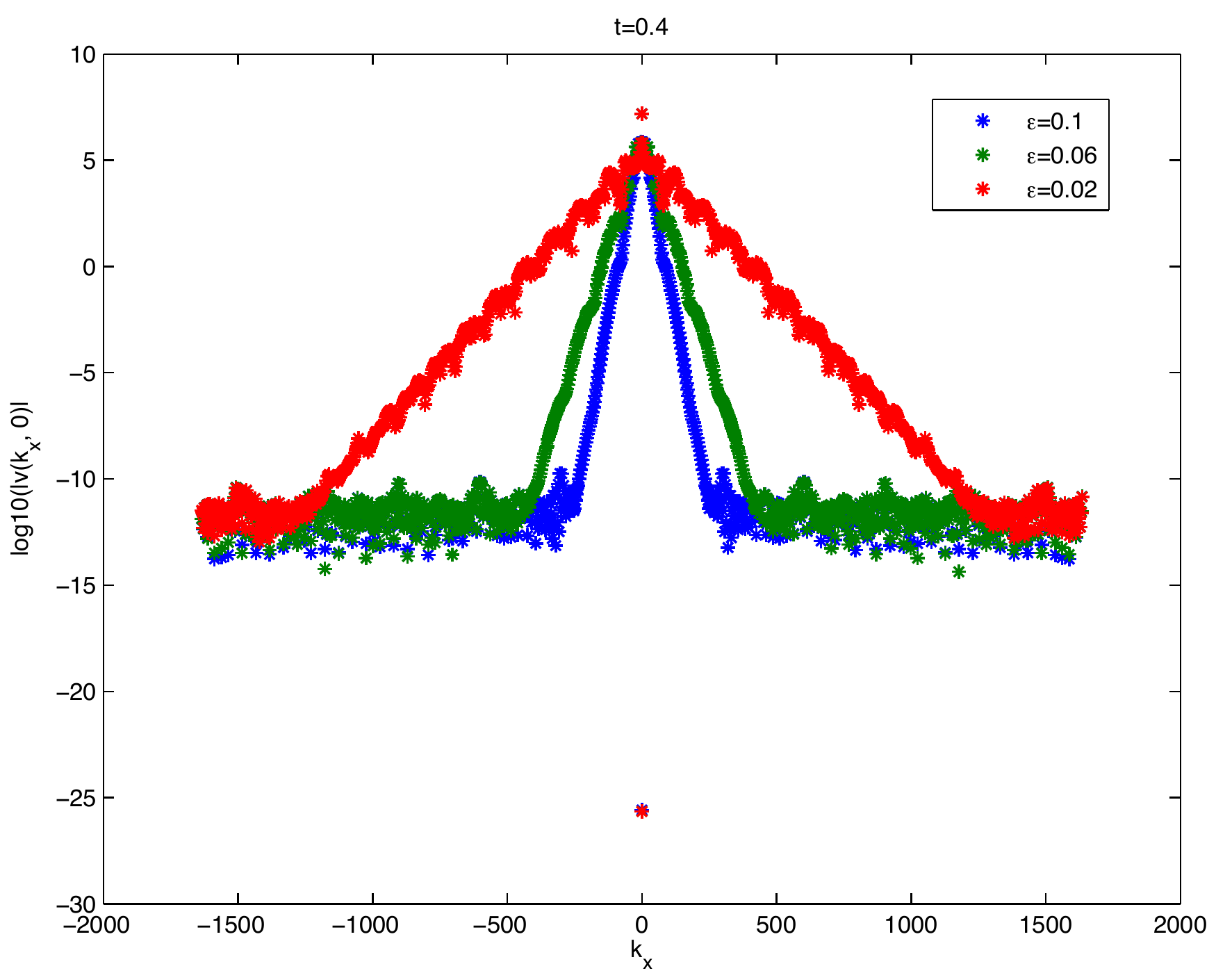}
  \includegraphics[width=0.4\textwidth]{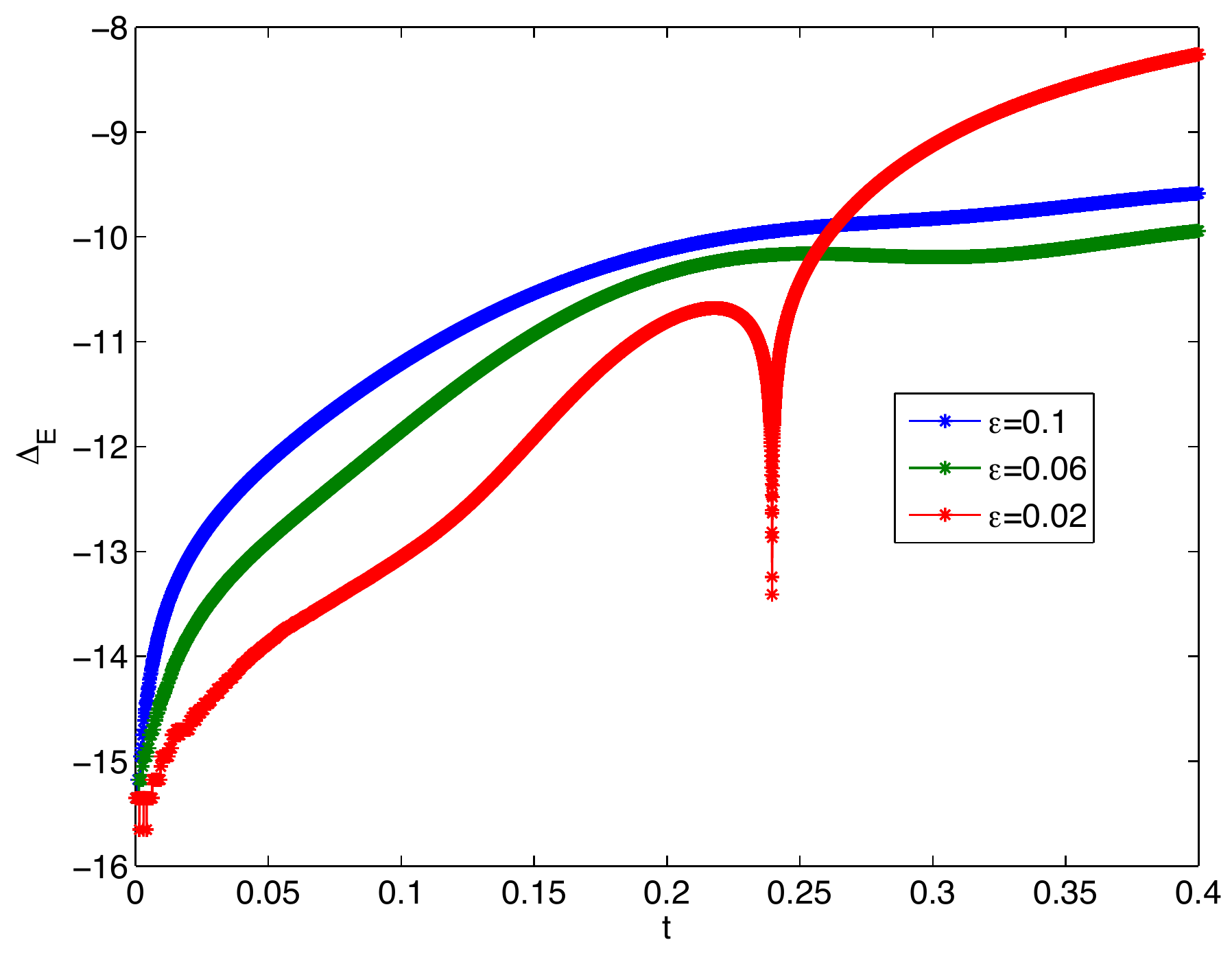}  
 \caption{Fourier coefficients to the solution of the KP I equation 
 in the small dispersion limit on the $k_x$-axis for several values 
 of $\epsilon$, 
 at $t=0.4$ (left), and the 
 time evolution of the quantity $\Delta_E$ (\ref{deltaE}) (right)}
 \label{coeftctmax}
\end{figure}

As usual, we ensure that the system is numerically well resolved by checking the 
decay of the Fourier coefficients, and the conservation of the 
numerically computed mass.
The Fourier coefficients of the solution of the KP I equation in the 
small dispersion limit are plotted on the $k_x$-axis for several 
values of $\epsilon$ in Fig.~\ref{coeftctmax} 
at $t=0.4$. They decrease to machine precision in all cases.
The quantity $\Delta_{E}$ is always smaller than $10^{-8}$, which 
indicates a numerical error well below plotting accuracy.

\begin{figure}[htb!]
\centering
  \includegraphics[width=0.45\textwidth]{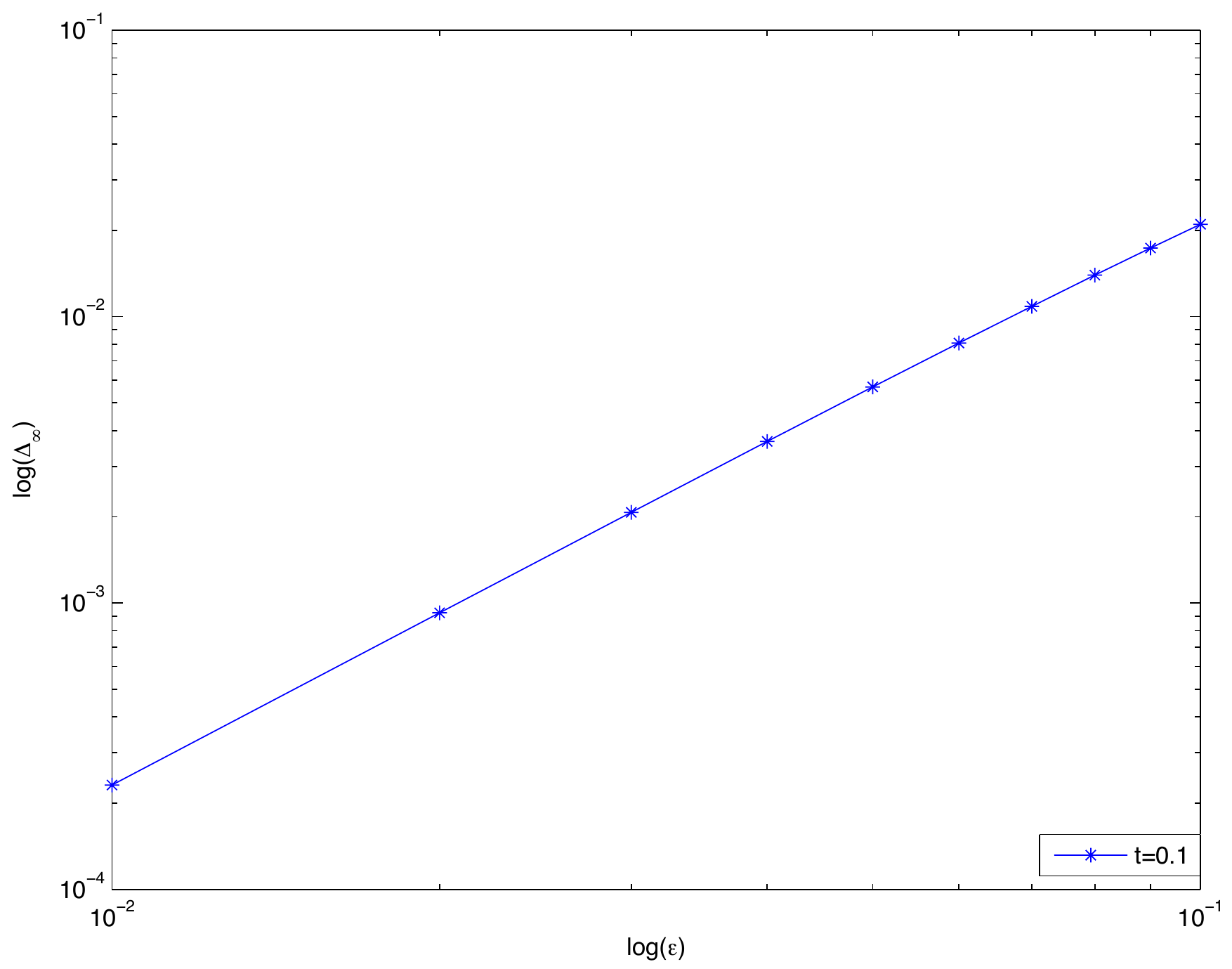} 
 \includegraphics[width=0.45\textwidth]{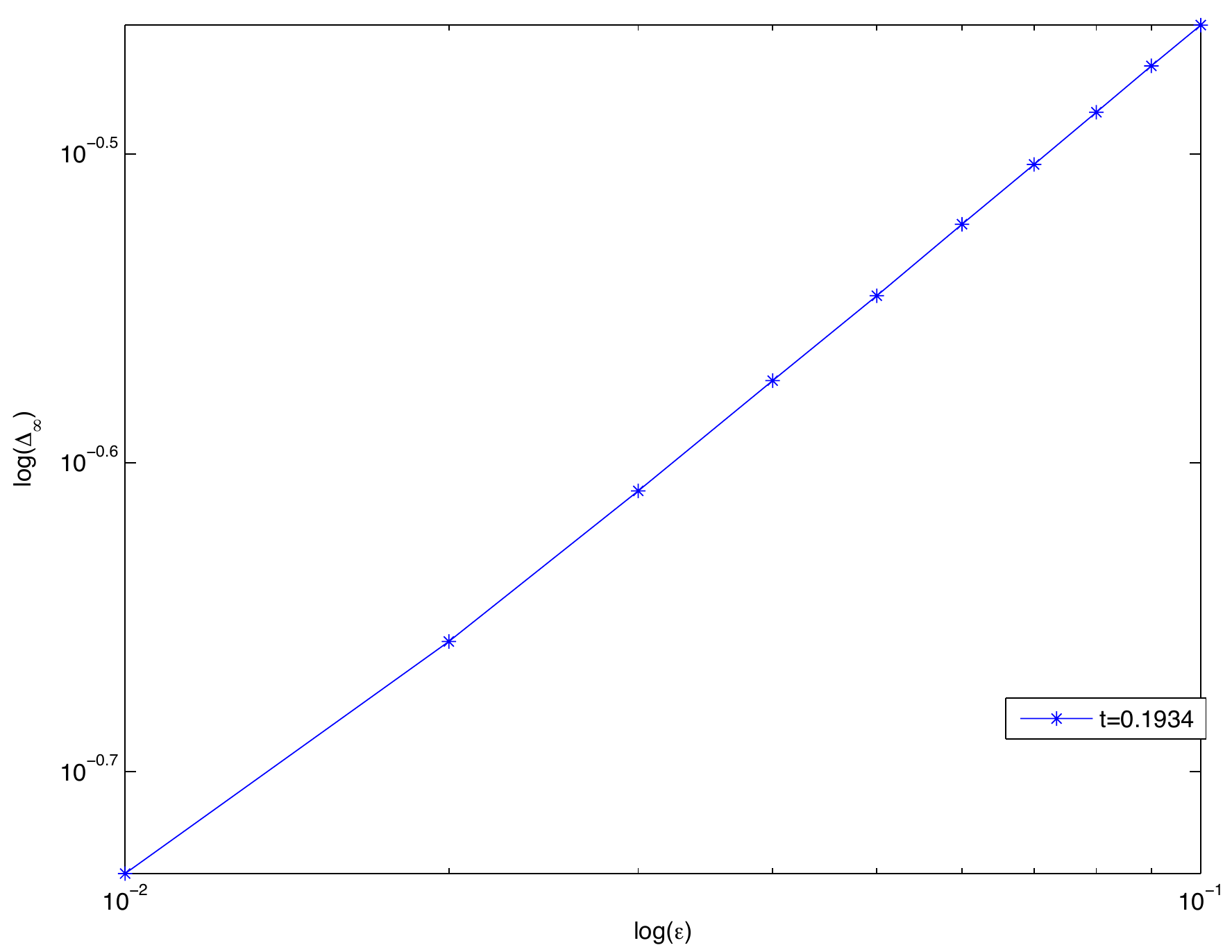} 
 \caption{$L_{\infty}$ norm  $\Delta_{\infty}$ of the difference 
 between KP and dKP solution for the initial data (\ref{uini2}) in dependence of $\epsilon$ at $t=0.1\sim t_c/2$ (left) and at $t_c=0.1934$ (right) for several values of $\epsilon$. }
 \label{scaltc}
\end{figure}

An important question is the scaling with $\epsilon$ of the 
$L_{\infty}$ norm of the difference between dKP 
and KP solutions for the same initial data. 
The $L_{\infty}$ norm of this difference is shown in 
Fig.~\ref{scaltc} at $t=0.1\sim t_c/2$ (left) and at $t_c=0.1934$ 
(right) in dependence of $\epsilon$ for  $0.01\leq  \epsilon \leq 0.1$. 
A linear regression analysis ($\log_{10} \Delta_{\infty} = a \log_{10} \epsilon + b$ ) shows that $\Delta_{\infty}$ decreases as 
\begin{align}
\mathcal{O} \left( \epsilon^{1.96} \right)  \sim \mathcal{O} \left( \epsilon^{2} \right) \,\, \mbox{at} \,\, t=0.1\sim t_c/2 , \,\, \mbox{with}\,\, a=1.961 \,\,\mbox{and} \,\, b= 0.6840  \\
\mathcal{O} \left( \epsilon^{0.28} \right)  \sim \mathcal{O} \left( \epsilon^{2/7} \right) \,\, \mbox{at} \,\, t=t_c=0.1934, \,\, \mbox{with}\,\, a=0.278  \,\,\mbox{and} \,\,  b= -0.4214.
\end{align}
In both cases, the correlation coefficient is $r = 0.999$. 

\begin{figure}[htb!]
\centering
  \includegraphics[width=0.8\textwidth]{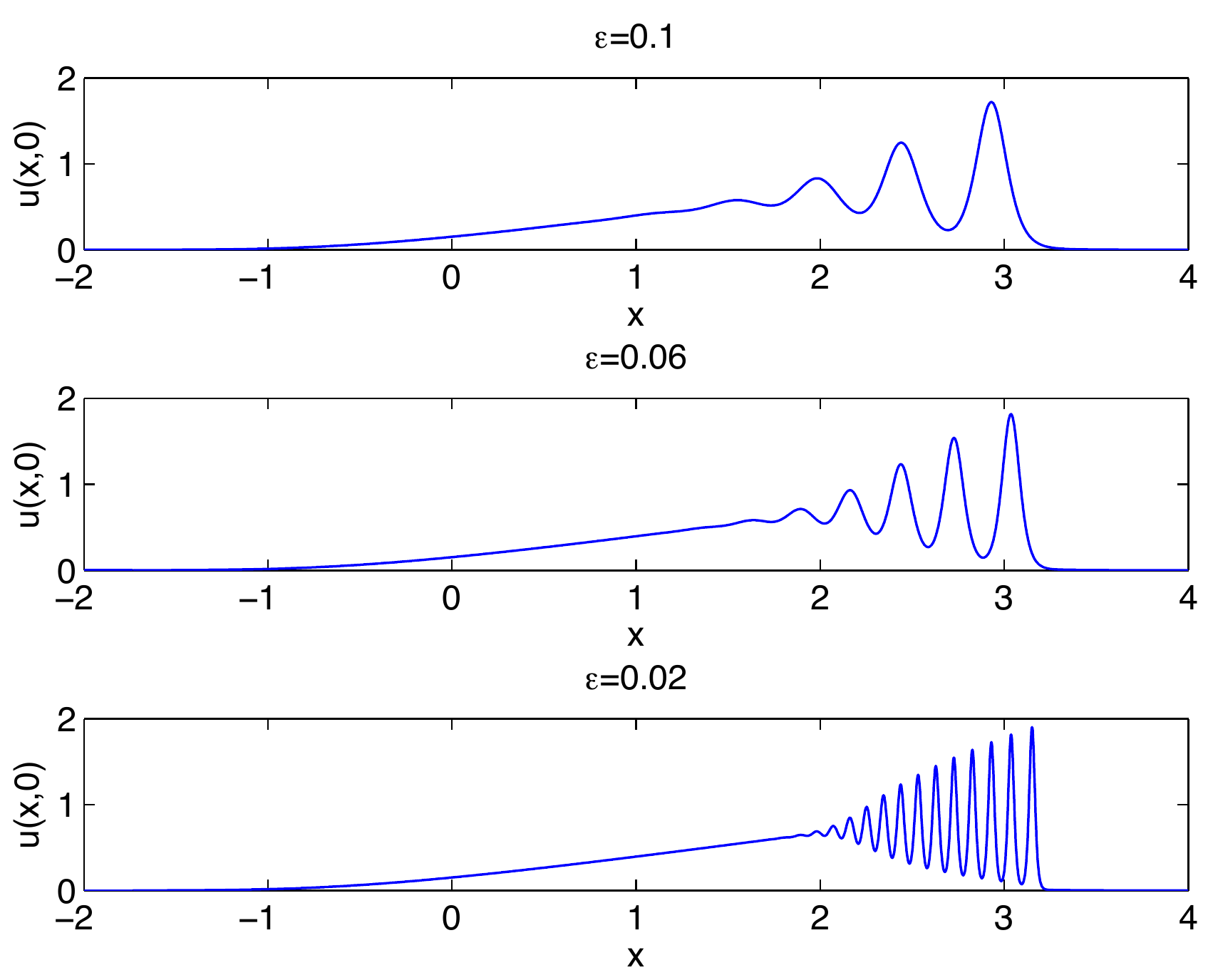} 
\caption{ Solutions to the KP II equation in the small dispersion limit on the $x$-axis at  $t=0.4$  for several values of $\epsilon$} 
 \label{ut4epsskp2uexp}
\end{figure}

The behavior of solutions of KP II in the small dispersion limit, for 
initial data of the form (\ref{uini2}) is completely similar to the KP I case, see for example 
Fig.~\ref{ut4epsskp2uexp}, where we show them plotted on the 
$x$-axis, for different values of $\epsilon$. Note that the 
oscillations are somewhat smaller than in the KP I case because of 
the discussed defocusing effect of KP II. 
We repeat the same study as before, and find that also in this case, the system is well resolved until $t=0.4$, and that 
the 
$L_{\infty}$ norm of the difference between dKP II
and KP II solutions for initial data (\ref{uini2}) decreases as in the case of dKP I/KP I:
We find that 
  $\Delta_{\infty}$ decreases as 
\begin{align}
\mathcal{O} \left( \epsilon^{1.96} \right)  \sim \mathcal{O} \left( \epsilon^{2} \right) \,\, \mbox{at} \,\, t=0.1\sim t_c/2 , \,\, \mbox{with}\,\, a=1.9608 \,\,\mbox{and} \,\, b= 0.6840  \\
\mathcal{O} \left( \epsilon^{0.28} \right)  \sim \mathcal{O} \left( \epsilon^{2/7} \right) \,\, \mbox{at} \,\, t=t_c=0.1934, \,\, \mbox{with}\,\, a=0.2751  \,\,\mbox{and} \,\,  b= -0.4122.
\end{align}
The correlation coefficient is $r = 0.999$ for both times.

\subsection{KP  solution in the small dispersion limit for initial 
data localized in both spatial dimension}
The same study as in the previous subsection will be carried our for 
the initial data (\ref{uini1}). Interestingly the results are very 
similar though the (first)  dKP shock appears here in just one point.

\begin{figure}[htb!]
\centering
\includegraphics[width=\textwidth]{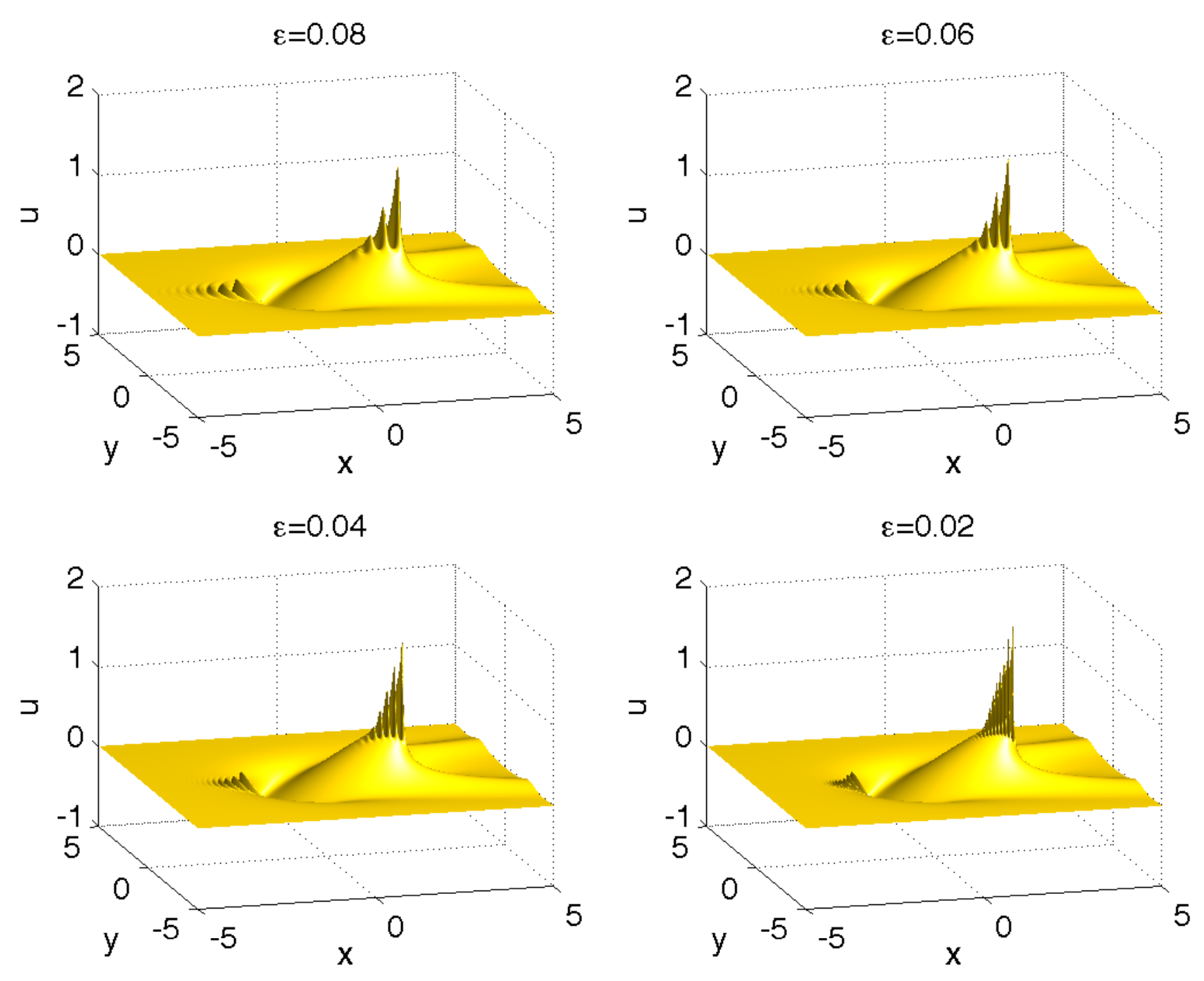}
\caption{ Solutions to the KP I equation in the small dispersion 
limit for $t=0.4$  for several values of $\epsilon$} 
 \label{ut4sechepss2d}
\end{figure}

\begin{figure}[htb!]
\centering
  \includegraphics[width=0.8\textwidth]{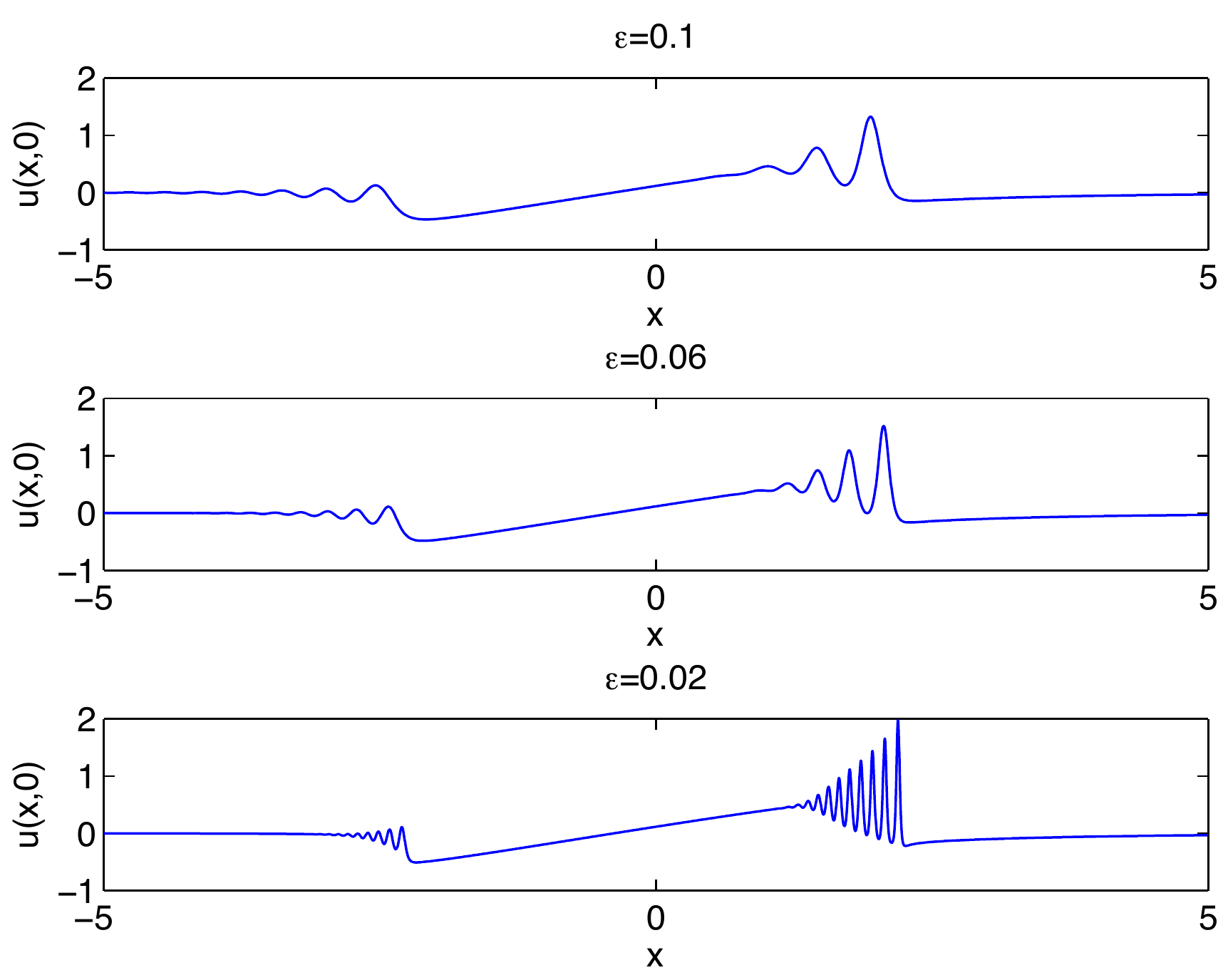} 
\caption{ Solutions to the KP I equation in the  small dispersion 
limit on the $x$-axis for $t=0.4$  for several values of $\epsilon$} 
 \label{ut4sechepss}
\end{figure}

\begin{figure}[htb!]
\centering
\includegraphics[width=0.8\textwidth]{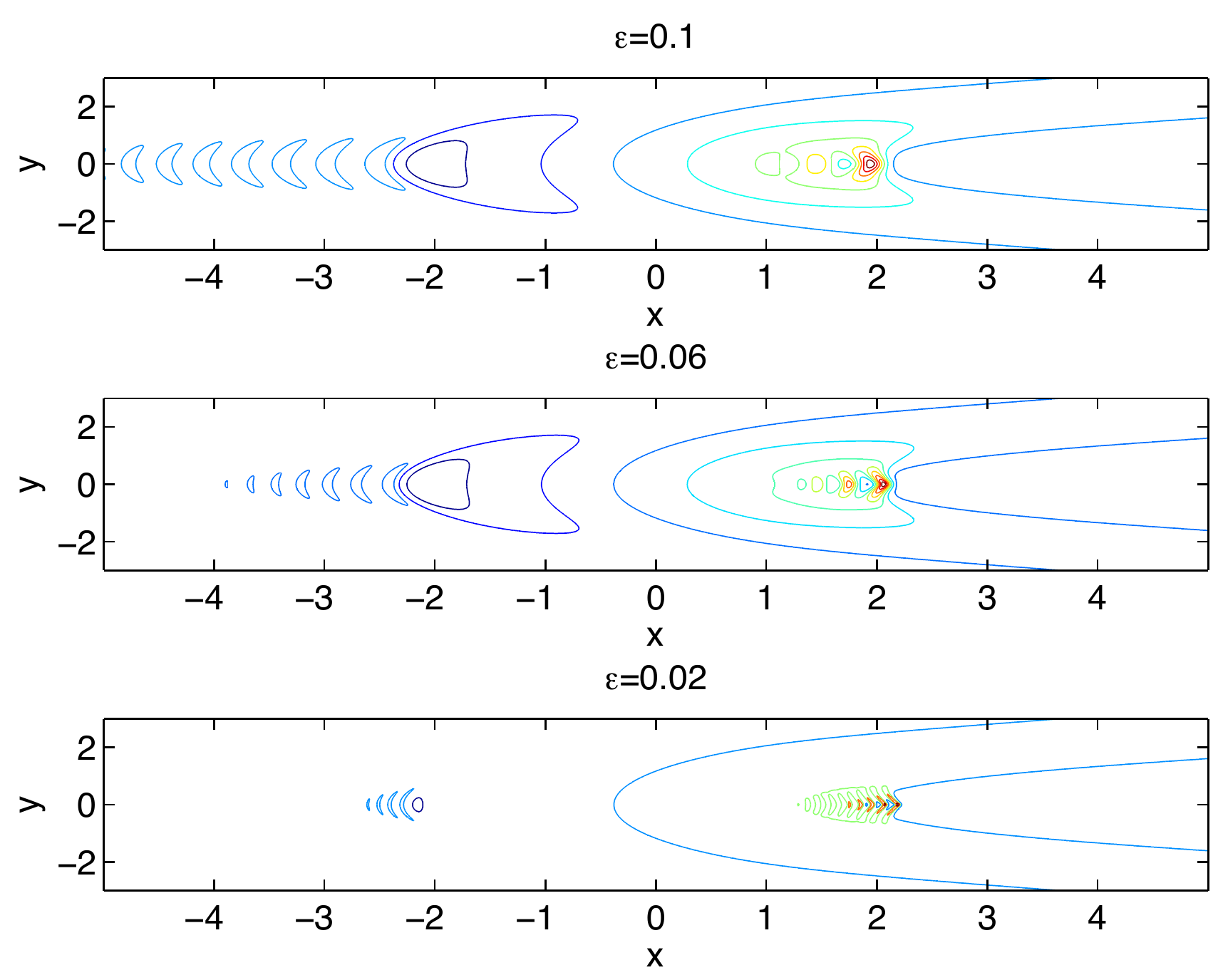} 
\caption{Contour plots of the solutions to the KP I equation in the 
small dispersion limit for $t=0.4$  for several values of $\epsilon$} 
 \label{cont1usech}
\end{figure}

For the initial data of the form (\ref{uini1}) one gets the 
dispersive shock already studied in \cite{KSM}. The solution at time 
$t=0.4>t_{c}$ can be seen in Fig.~\ref{ut4sechepss2d} for several 
values of $\epsilon$. It can be recognized that the number of oscillations 
increases with decreasing $\epsilon$, and that the oscillations are more 
and more confined to a well defined zone. 
This is even more visible on the $x$-axis as can be inferred from 
Fig.~\ref{ut4sechepss}.
The corresponding contour plots  are shown in Fig.~\ref{cont1usech}.

\begin{figure}[htb!]
\centering
  \includegraphics[width=0.45\textwidth]{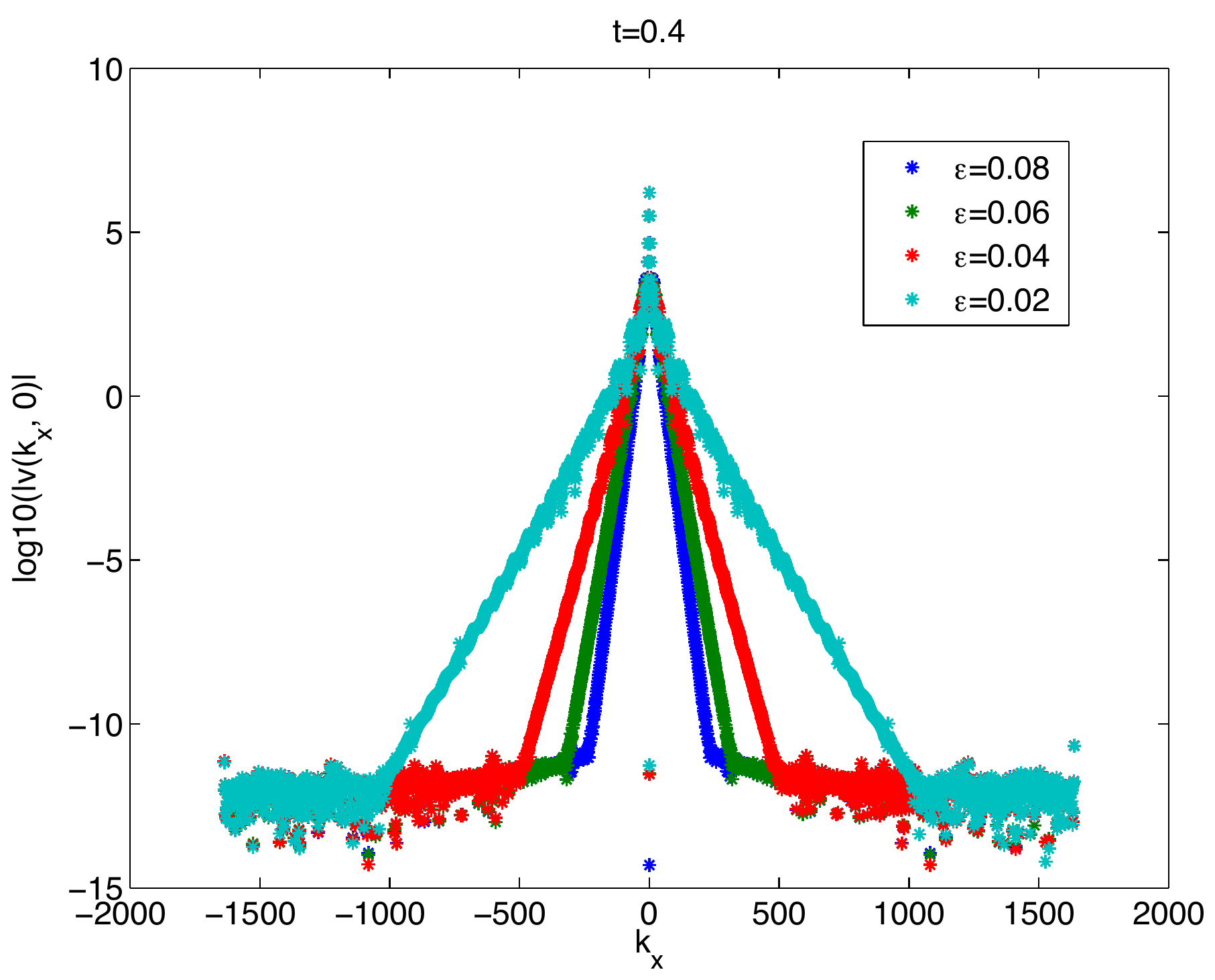} 
  \includegraphics[width=0.45\textwidth]{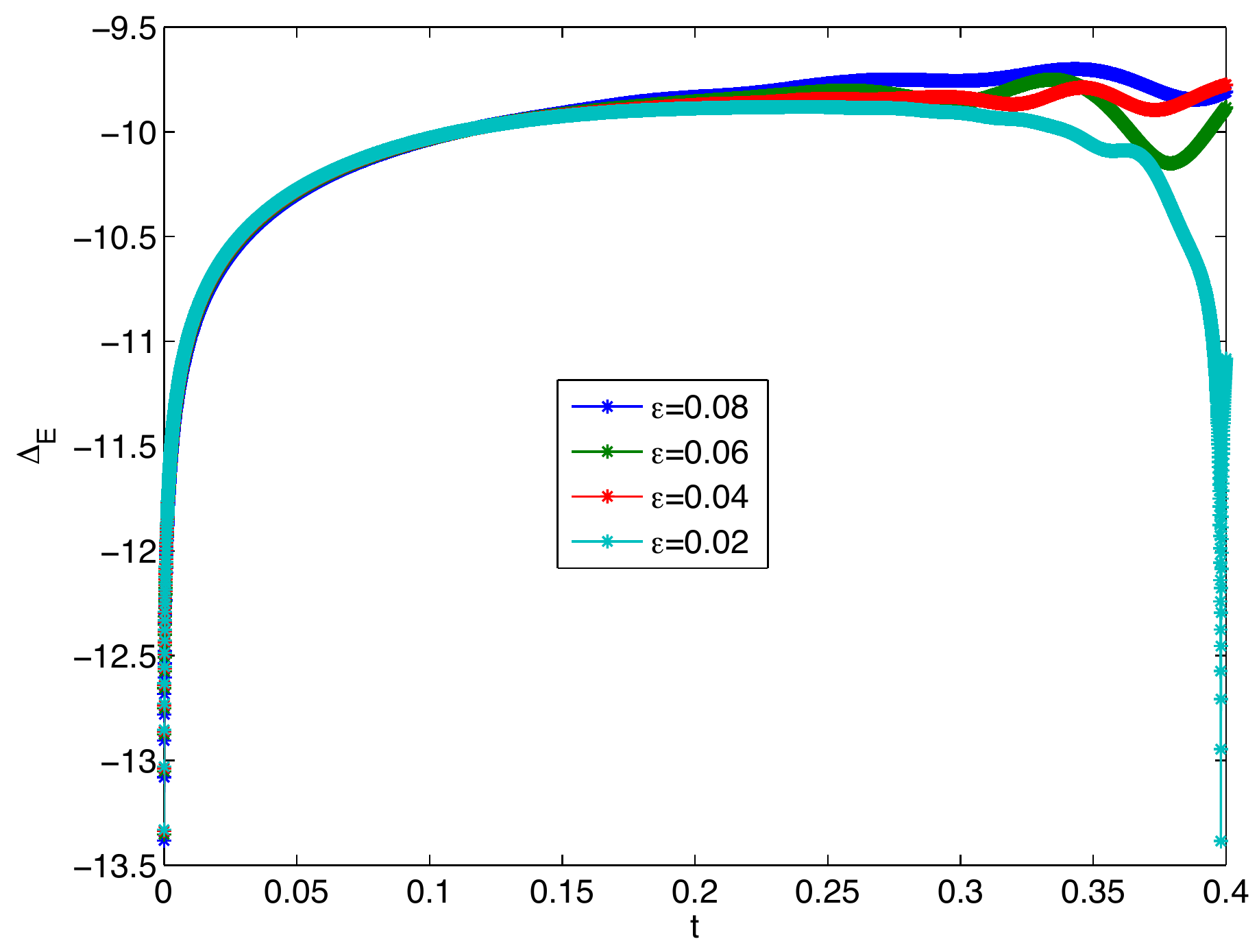} 
 \caption{Fourier coefficients to the solution of the KP I equation 
 in the small dispersion limit on the $k_x$-axis for several values 
 of $\epsilon$, at $t=0.4$ (left)  and the time evolution of the 
 quantity $\Delta_E$ (\ref{deltaE}) (right)}
 \label{coeftctmaxusech}
\end{figure}

The Fourier coefficients of the solution  on the $k_x$-axis can be seen  in Fig.~\ref{coeftctmaxusech} for several values of $\epsilon$ at $t=0.4$. They decrease to machine precision in all cases.
The quantity $\Delta_{E}$ used as an indicator for the numerical 
accuracy is always smaller than $10^{-9}$, see Fig. 
\ref{coeftctmaxusech}, which is more than satisfactory for our 
purposes.

For the initial data (\ref{uini1}), we find a similar scaling 
with $\epsilon$ as for the data (\ref{uini2}): 
We obtain via a linear regression analysis ($\log_{10} 
\Delta_{\infty} = a \log_{10} \epsilon + b$ ) that the $L_{\infty}$ 
norm  $\Delta_{\infty}$ of the difference between dKP and KP solution decreases as 
\begin{align}
\mathcal{O} \left( \epsilon^{1.94} \right)  \sim \mathcal{O} \left( \epsilon^{2} \right) \,\, \mbox{at} \,\, t=0.125 \ll t_c , \,\, \mbox{with}\,\, a=1.943 \,\,\mbox{and} \,\, b= 1.2239  \\
\mathcal{O} \left( \epsilon^{0.30} \right)  \sim \mathcal{O} \left( \epsilon^{2/7} \right) \,\, \mbox{at} \,\, t=t_c=0.2216, \,\, \mbox{with}\,\, a=0.30  \,\,\mbox{and} \,\,  b= -0.5122.
\end{align}
In both cases, the correlation coefficient is $r = 0.999$. 
\begin{figure}[htb!]
\centering
  \includegraphics[width=0.45\textwidth]{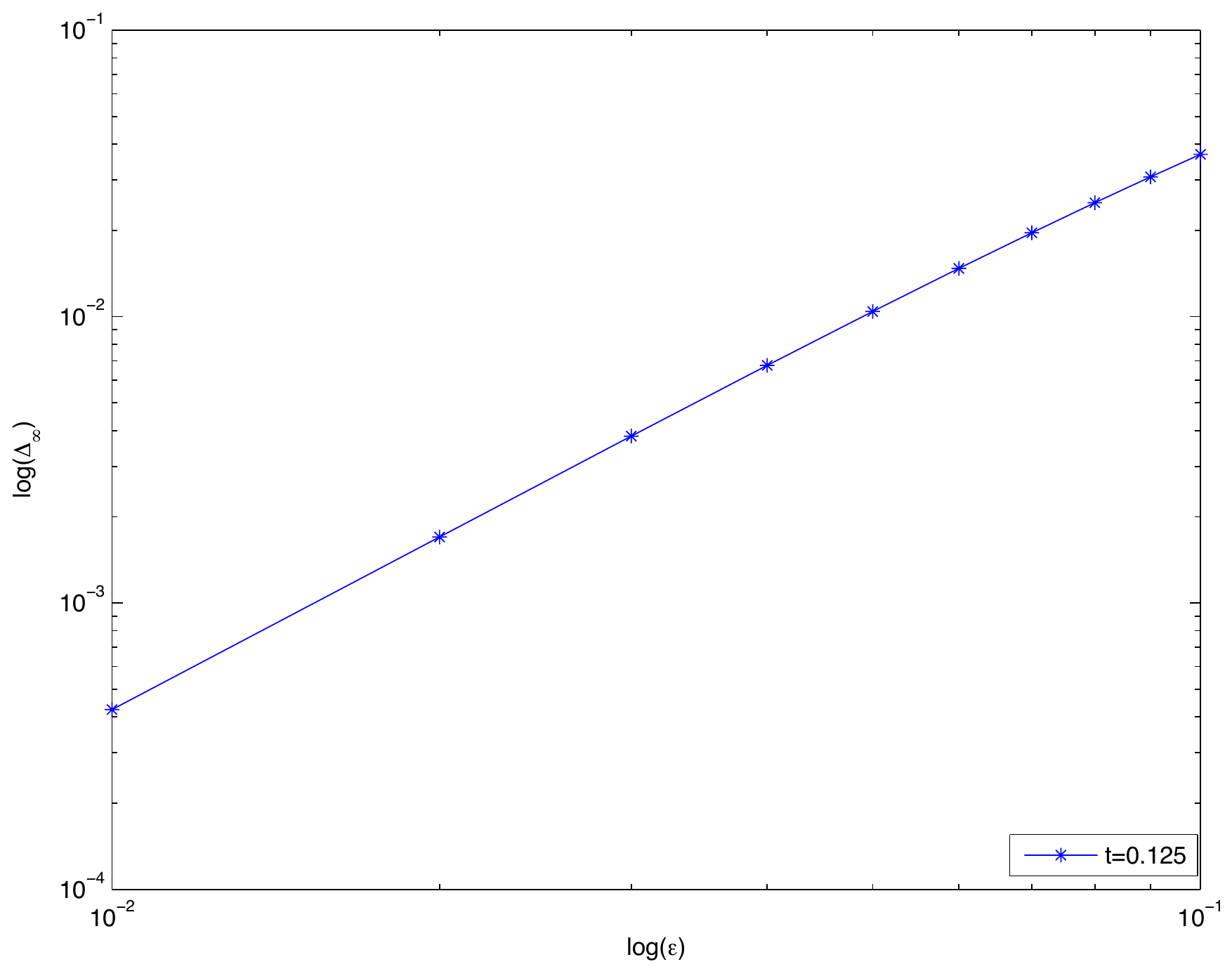} 
 \includegraphics[width=0.45\textwidth]{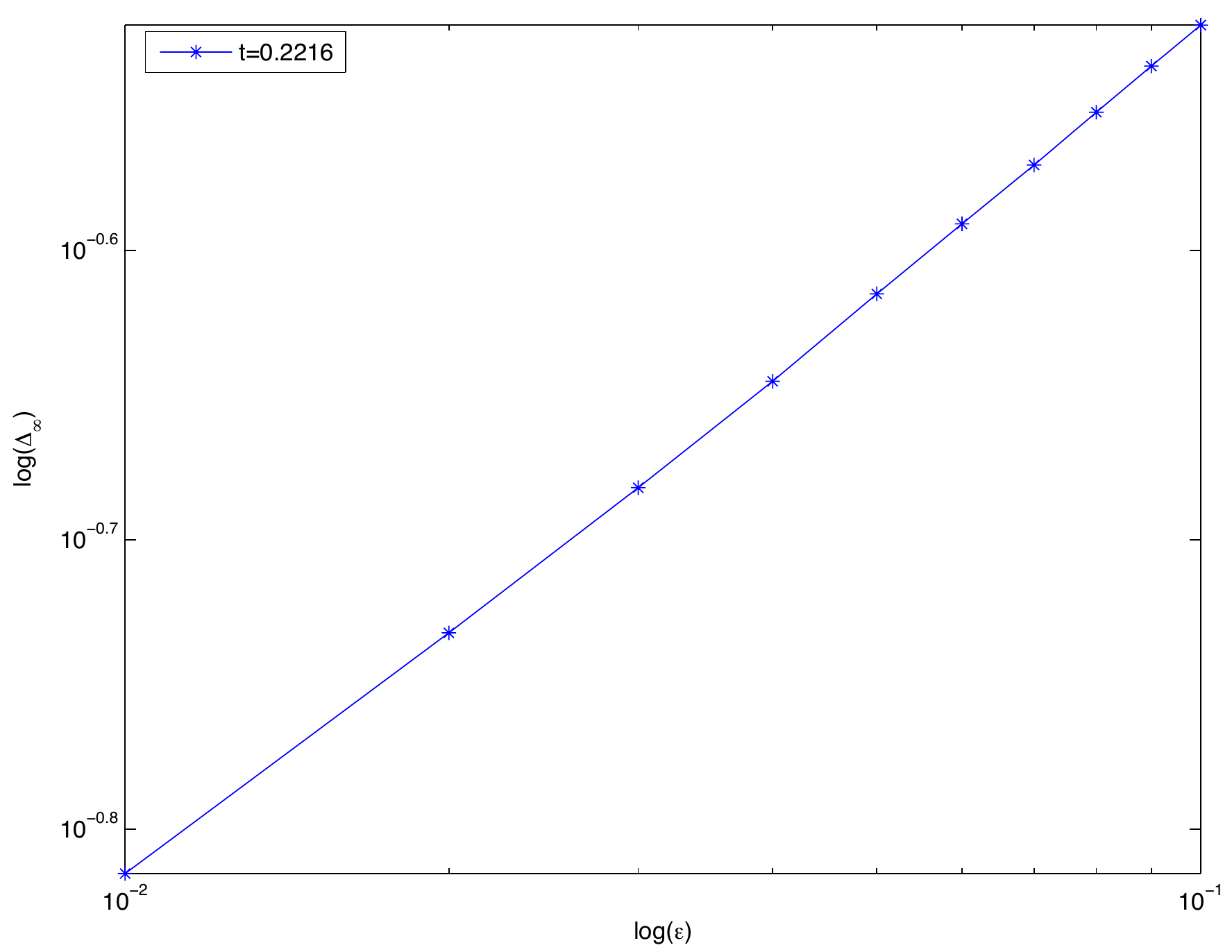} 
 \caption{Error $\Delta_{\infty}$ in dependence of $\epsilon$ at $t=0.125\sim t_c/2$ (left) and at $t_c=0.2216$ (right) for several values of $\epsilon$. }
 \label{scaltcu2}
\end{figure}

\begin{figure}[htb!]
\centering
  \includegraphics[width=0.8\textwidth]{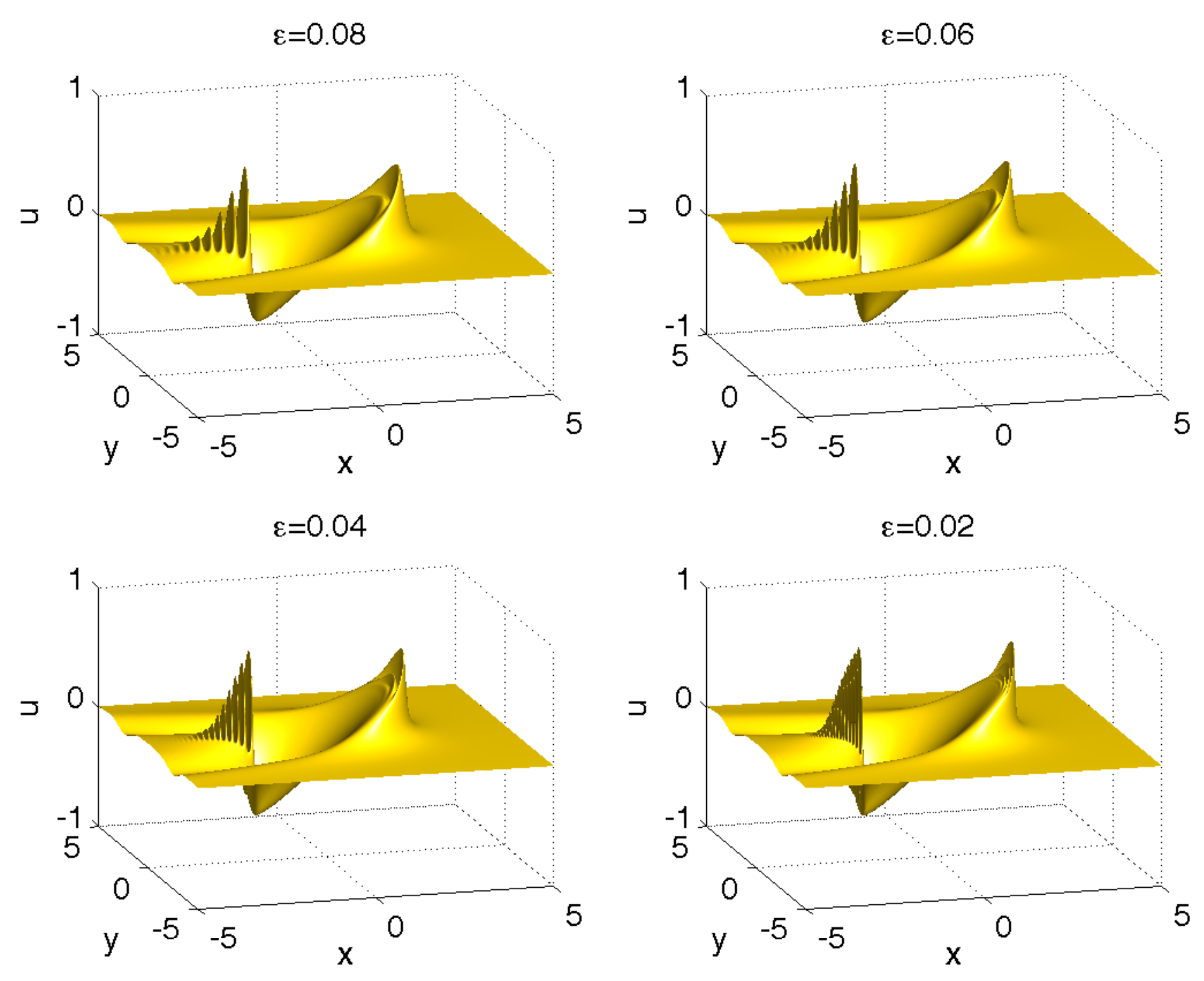}
\caption{ Solutions to the KP II equation in the small dispersion limit at the maximal time of computation, $t=0.4$  for several values of $\epsilon$}
 \label{ut4sech2dkp2}
\end{figure}

\begin{figure}[htb!]
\centering
  \includegraphics[width=0.8\textwidth]{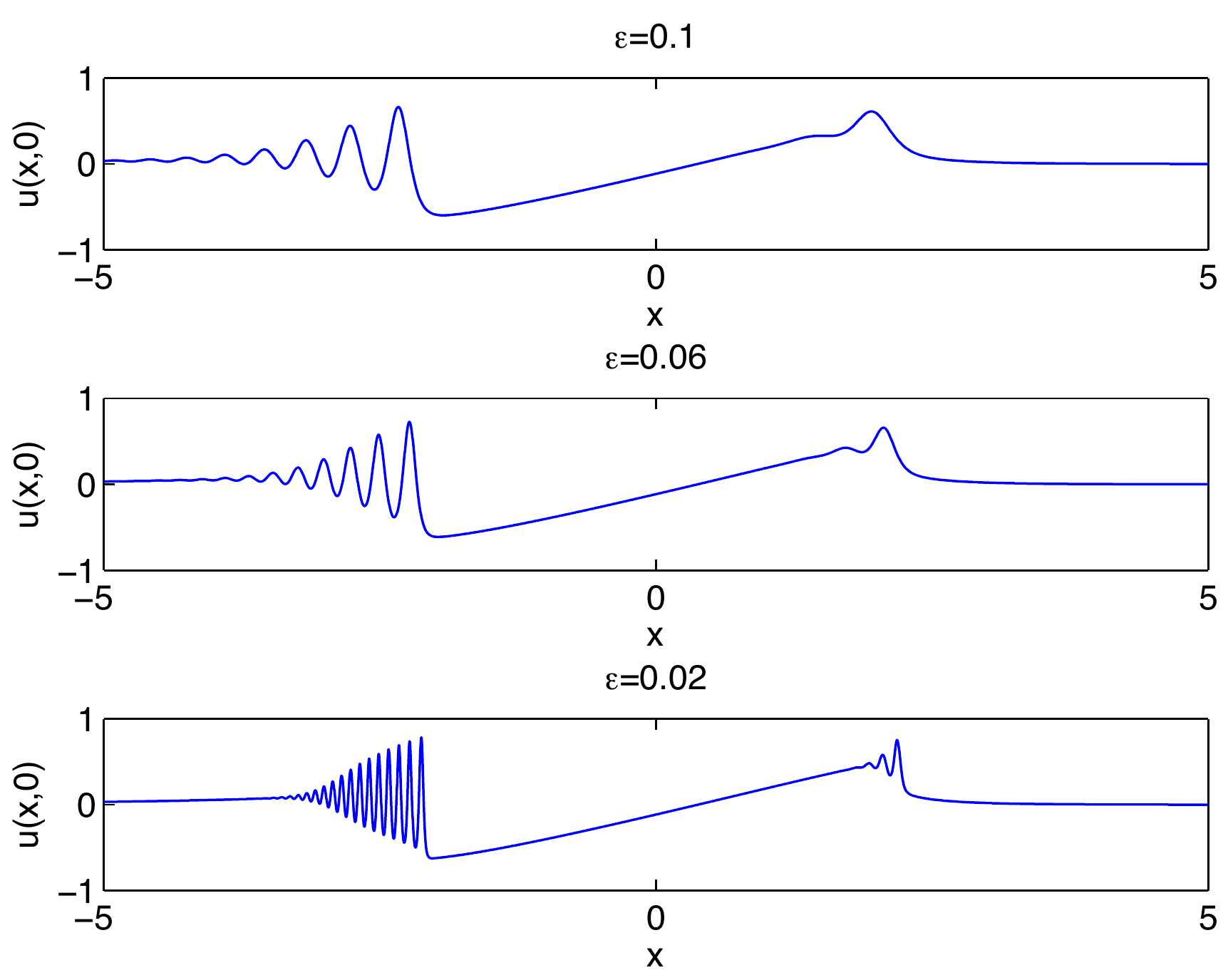} 
\caption{ Solutions to the KP II equation in the small dispersion limit plotted on the $x$-axis at the maximal time of computation, $t=0.4$  for several values of $\epsilon$}
 \label{ut4sechkp2}
\end{figure}

\begin{figure}[htb!]
\centering
\includegraphics[width=0.8\textwidth]{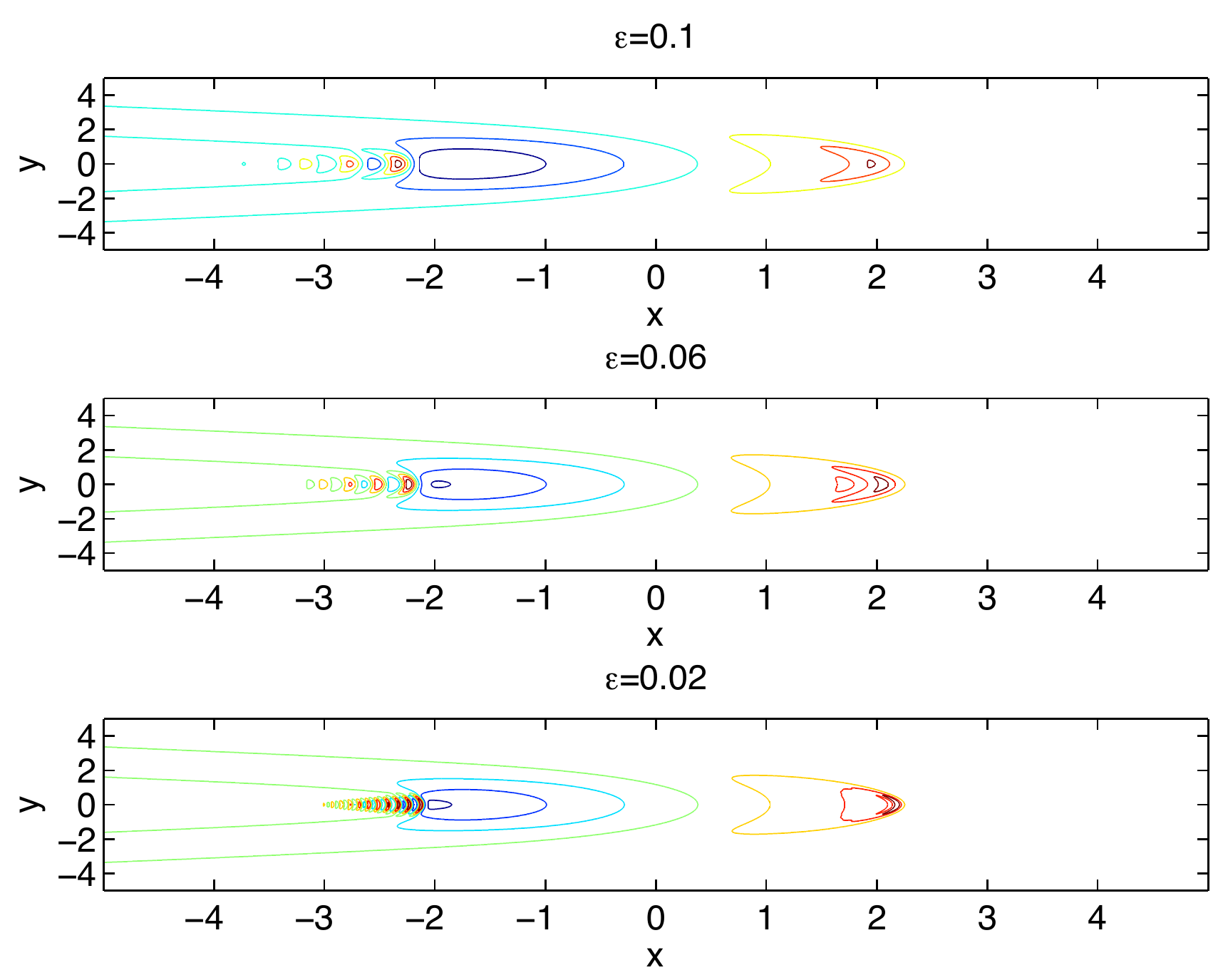} 
\caption{Contour of the solutions to the KP II equation in the small dispersion limit at the maximal time of computation, $t=0.4$  for several values of $\epsilon$}
 \label{cont1usechkp2}
\end{figure}

\begin{figure}[htb!]
\centering
  \includegraphics[width=0.45\textwidth]{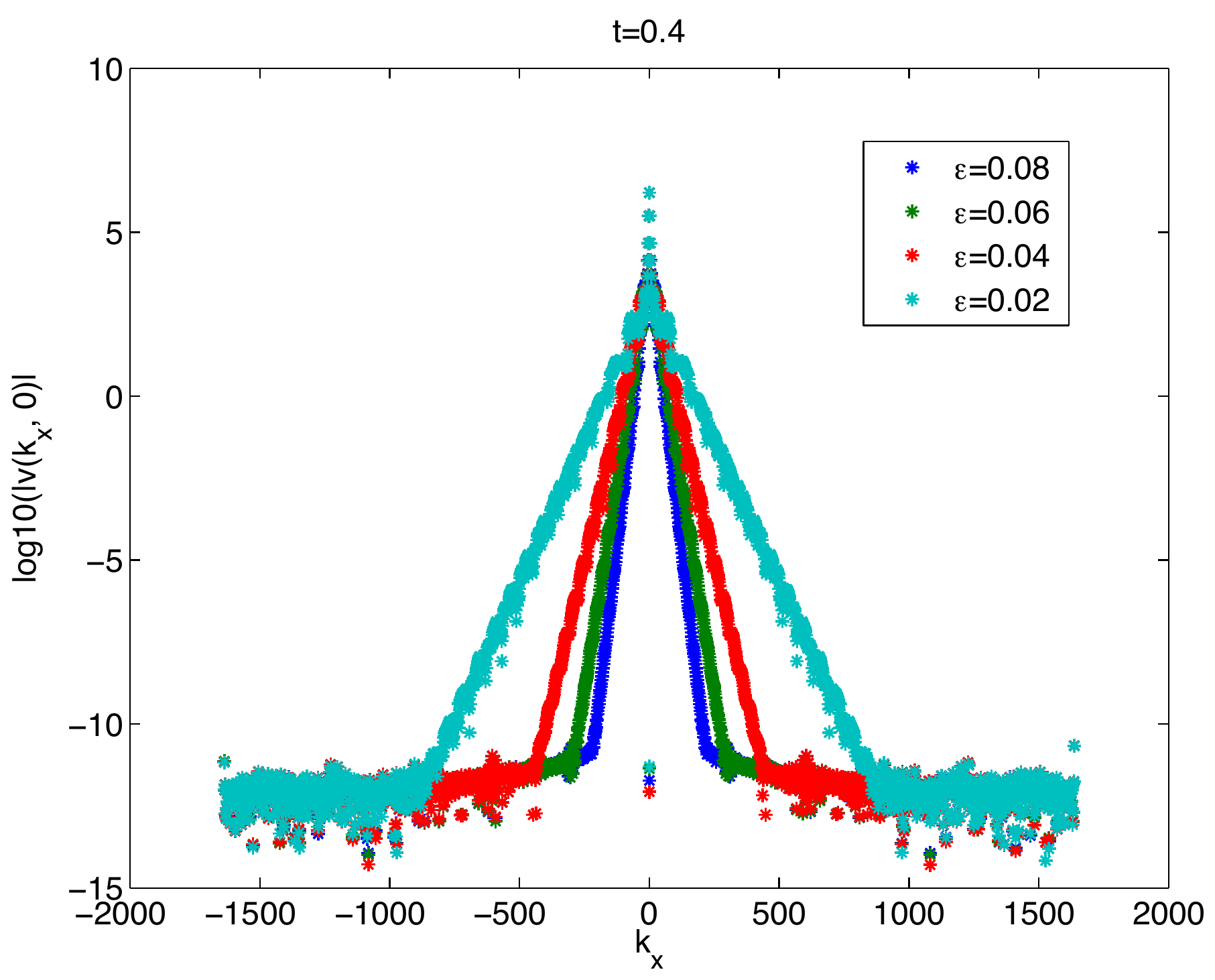}
    \includegraphics[width=0.45\textwidth]{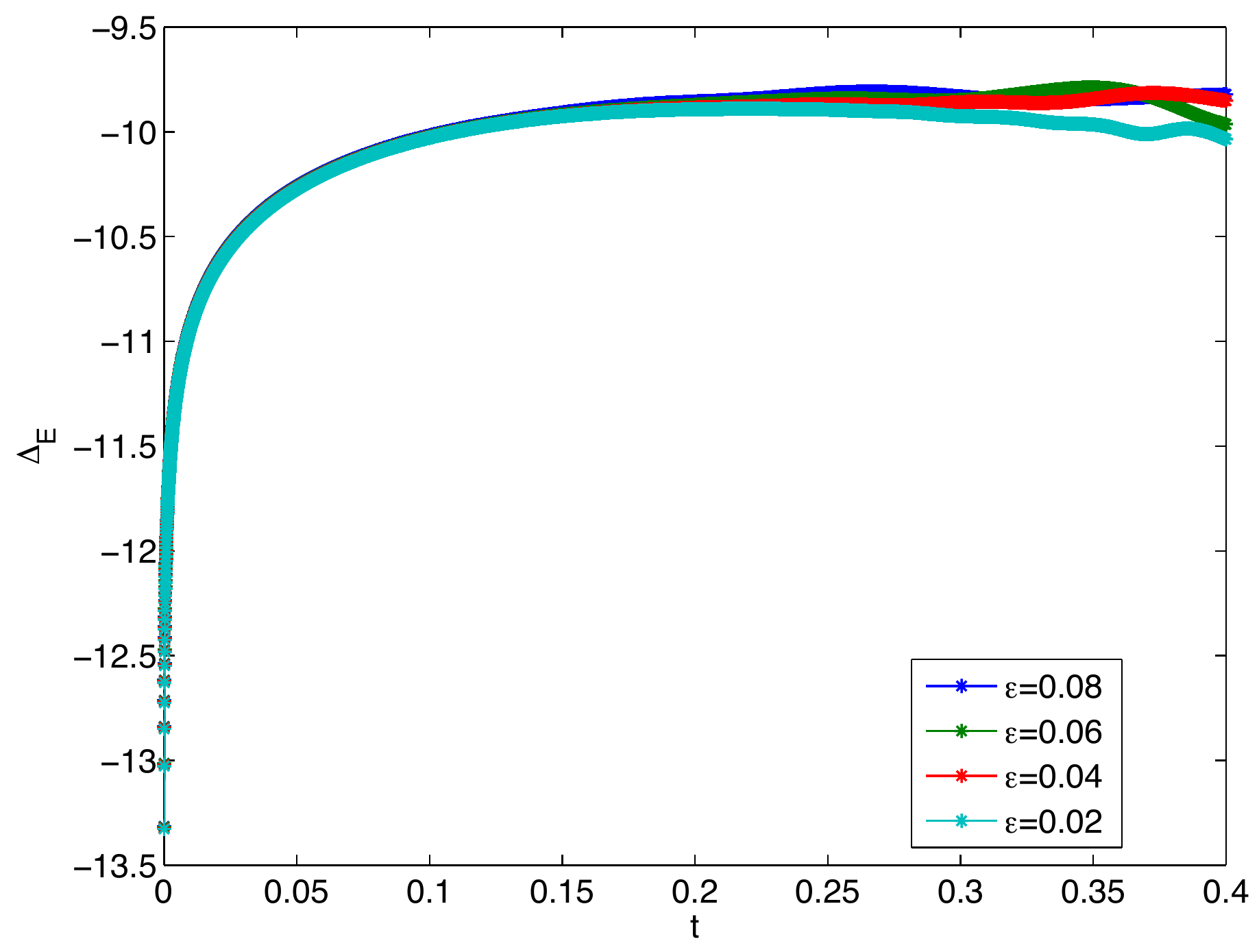}
 \caption{Fourier coefficients to the solution of the KP II equation 
 in the small dispersion limit, plotted on the $k_x$-axis for several values of $\epsilon$, at $t=0.4$ (left)  and time evolution of the mass conservation $\Delta_E$ (right)}
 \label{coeftctmaxusechkp2}
\end{figure}

The results for the KP II case are very similar to the KP I. 
We obtain again a dispersive shock in the small dispersion limit, 
see Fig.~\ref{ut4sechkp2} and Fig.~\ref{ut4sech2dkp2}, already 
observed in \cite{KSM} for the same initial data. Since the tails 
with the algebraic fall off are now directed towards $-\infty$, the 
first shock for dKP as well as the stronger oscillation appear for 
negative $x$. The mass transfer towards $-\infty$ simply implies that 
there is more mass in this region which triggers a dispersive shock. 
The corresponding solutions on the $x$-axis can be seen in 
Fig.~\ref{ut4sechkp2}.
The related contour plots are shown in Fig. \ref{cont1usechkp2}.
The Fourier coefficients of the solution to the KP II equation in the 
small dispersion limit are plotted on the $k_x$-axis for several 
values of $\epsilon$ in Fig.~\ref{coeftctmaxusechkp2} at $t=0.4$. They decrease to machine precision
 in all cases, which ensures that the system is well resolved until $t=0.4$.
The conservation of the numerically computed mass always reaches a 
precision of better than $10^{-9}$,  see Fig. \ref{coeftctmaxusechkp2}.

The scaling laws are also almost identical to the previous case: 
We find with a linear regression analysis ($\log_{10} \Delta_{\infty} 
= a \log_{10} \epsilon + b$) that  $\Delta_{\infty}$ decreases as 
\begin{align}
\mathcal{O} \left( \epsilon^{1.98} \right)  \sim \mathcal{O} \left( \epsilon^{2} \right) \,\, \mbox{at} \,\, t=0.11\sim t_c/2 , \,\, \mbox{with}\,\, a=1.9815 \,\,\mbox{and} \,\, b= 0.8998  \\
\mathcal{O} \left( \epsilon^{0.29} \right)  \sim \mathcal{O} \left( \epsilon^{2/7} \right) \,\, \mbox{at} \,\, t=t_c=0.2216, \,\, \mbox{with}\,\, a=0.295  \,\,\mbox{and} \,\,  b= -0.5288.
\end{align}
In both cases, the correlation coefficient is $r = 0.999$. 

Thus in all cases, we find approximately the same scaling laws as in the KdV case.

 \section{Outlook}
In this paper we have shown that it is possible to use asymptotic 
Fourier analysis to determine the time of the singularity formation 
in nonlinear evolution equations. A precondition for the 
applicability of the method is obviously that the PDE is solved with 
sufficient accuracy and resolution. The former can be controlled via 
a conserved quantity, typically the $L_{2}$ norm or the energy, the 
latter via the decrease of the Fourier coefficients for large wave 
numbers. 
The order of magnitude of the Fourier coefficients for the highest 
wave numbers is in some sense a measure for the numerical resolution. The 
conserved quantities, which will numerically depend on time because 
of unavoidable errors,  cannot indicate a higher precision than this 
resolution, and typically overestimate the accuracy by 2-3 orders of 
magnitude. If sufficient numerical resolution is provided, normally 
more than $2^{14}$ Fourier modes, a minimal wave number $k_{min}$ for the least 
square fitting to the asymptotic formula (\ref{fourierasym}) has to be 
chosen. This has to be clearly larger than 1, but should include 
sufficient Fourier coefficients with larger modulus than the rounding 
error for the considered times. Values between 10 and 100 appear to 
be convenient here. Then one concludes an aimed at fitting error in dependence 
of the resolution (several runs might be necessary to make this self 
consistent, but generally it can be chosen of the same order as the 
estimated accuracy of the solution). The upper limit for the wave numbers for which the 
fitting is performed is then chosen as the maximal $|k|$ for which
the difference 
between fitted curve and Fourier coefficients on the whole considered 
interval is smaller than the numerical accuracy of the solution. If 
this cannot be achieved with at least half of the Fourier 
coefficients with a modulus larger than the rounding error, the number of Fourier modes 
has to be increased. In this way it is possible to trace the 
singularity with acceptable precision up to the point where it hits 
the axis, which defines the critical time. This approach gives 
the critical quantities with the necessary precision to study the scaling of 
solutions to a dispersive regularization as KdV for Hopf. 

In $2+1$ dimensions and higher, the needed resolution can only be 
achieved on parallel computers. As an example we studied shock 
formation in the dKP equation which could be treated with a parallel 
version of the code with the same precision as for Hopf for initial data 
which are infinitely extended in one direction. For data being 
localized in all spatial directions, the achievable accuracy was 
slightly lower. The reason for this is the appearance of a second 
singularity shortly after the first we studied. This is due to the 
fact that the initial data for KP have to satisfy the constraint 
(\ref{const}) in order to allow a smooth solution in time for $t=0$. 
The presence of a second singularity in the complex plane limits 
somewhat the accuracy of the fitting procedure for the Fourier 
coefficients. Nonetheless we could identify the critical point and 
the critical singularity with 
the needed precision. This allowed to study the scaling of the 
difference between KP solutions in the small dispersion limit and the 
dKP solution at break-up. It was shown that the same scaling (within 
numerical tolerance) is 
observed at the critical time as in the KdV case, $\epsilon^{2/7}$. 
The task is now to give an asymptotic description of the KP solutions 
in the small dispersion limit close to the break-up of dKP as in 
\cite{Dub06}.  

The numerical analysis presented in this paper for dKP will also be 
applied to other $2+1$- and higherdimensional PDEs, for instance the 
Davey-Stewartson equations in the semiclassical limit. One problem in 
this context is that the singularity might in these cases be 
genuinely multi-dimensional, not as for dKP where the gradient blows 
up only in one spatial direction. The essentially one-dimensional 
approach can of course still be applied in this case.  
Alternatively  instead of the Fourier transform $\hat{u}(k,t)$, one 
can consider as in \cite{SSP}
the angle averaged energy spectrum defined by
\[
 E(K,t) = \underset{K<|k'|<K+1}{\sum} |\hat{u}(k',t)|^{2},
\]
where 
$|k'|=\sqrt{ k_{x}^2 + k_{y}^2 }$.
Slightly weaker estimates hold for $E(K,t)$ for an analytic
function $u$ in $S_{\alpha}$, and thus, to apply the asymptotic 
fitting to the Fourier coefficients, one assumes
that 
$ E(K,t)=C(t)K^{-\alpha(t)} e^{-\delta(t)K}.
$
Similary to the one-dimensional case, the 
appearance of a real singularity implies that $\delta(t)$ vanishes at a finite time $t_c$.
This implies that 
sufficient numerical resolution has to be provided which can as in 
the present paper in 
general only be achieved on parallel computers.  With this approach it 
should be possible to study singularity formation in general 
nonlinear evolution equations without dissipation and dispersion. 

\bibliographystyle{siam}
\bibliography{bibliof}{}

\def\cprime{$'$}
\begin{thebibliography}{10}

\bibitem{AC}
{\sc M.~Ablowitz and P.~Clarkson}, {\em Solitons, nonlinear {E}volution
  {E}quations and inverse {S}cattering}, London Mathematical Society Lecture
  Note Series, 149,  (1991).

\bibitem{alinhac}
{\sc S.~Alinhac}, {\em Blowup of small data solutions for a quasilinear wave
  equation in two space dimensions}, Ann. Maths., 149 (1999), pp.~97--127.

\bibitem{arnold}
{\sc V.~I. Arnol{\cprime}d, V.~V. Kozlov, and A.~I. Ne{\u\i}shtadt}, {\em
  Dynamical {S}ystems. {III}}, vol.~3 of Encyclopaedia of Mathematical
  Sciences, Springer-Verlag, Berlin, 1988.
\newblock Translated from the Russian by A. Iacob.

\bibitem{CR}
{\sc R.~E. Caflisch}, {\em Singularity formation for complex solutions of the
  {$3$}{D} incompressible {E}uler equations}, Phys. D, 67 (1993), pp.~1--18.

\bibitem{Can}
{\sc C.~Canuto, M.~Y. Hussaini, A.~Quarteroni, and T.~A. Zang}, {\em Spectral
  methods}, Scientific Computation, Springer-Verlag, Berlin, 2006.
\newblock Fundamentals in single domains.

\bibitem{asymbook}
{\sc G.~Carrier and M.~K.~C. Pearson}, {\em Functions of a {C}omplex
  {V}ariable, {T}heory and {T}echnique}, Society for Industrial and Applied
  Mathematics (SIAM), Philadelphia, PA, 2005.

\bibitem{CM}
{\sc S.~Cox and P.~Matthews}, {\em Exponential time differencing for stiff
  systems}, Journal of Computational Physics, 176 (2002), pp.~430--455.

\bibitem{CG}
{\sc P.~E. Crouch and R.~Grossman}, {\em Numerical {I}ntegration of ordinary
  differential {E}quations on {M}anifolds}, J. Nonlinear Sci., 3 (1993),
  pp.~1--33.

\bibitem{Dub06}
{\sc B.~Dubrovin}, {\em On hamiltonian perturbations of hyperbolic systems of
  conservation laws, ii: universality of critical behaviour}, Comm. Math.
  Phys., 267 (2006), pp.~117 -- 139.

\bibitem{DE13}
{\sc B.~Dubrovin and M.~Elaeva}, {\em On critical behavior in nonlinear
  evolutionary pdes with small viscosity}, preprint, arXiv:1301.7216 (2013).

\bibitem{DGK11}
{\sc B.~Dubrovin, T.~Grava, and C.~Klein}, {\em Numerical {S}tudy of breakup in
  generalized {K}orteweg-de {V}ries and {K}awahara equations}, SIAM J. Appl.
  Math., 71 (2011), pp.~963--1008.

\bibitem{DMT01}
{\sc M.~Dunajski, L.~Mason, and K.~Tod}, {\em Einstein--weyl geometry, the dkp
  equation and twistor theory}, J. Geom. Phys., 37 (2001), pp.~63--93,.

\bibitem{FK04}
{\sc E.~Ferapontov and K.~Khusnutdinova}, {\em On the integrability of (2 + 1)-
  dimensional quasilinear systems.}, Comm. Math. Phys., 248(1) (2004),
  pp.~187--206.

\bibitem{FS}
{\sc A.~Fokas and L.~Sung}, {\em The inverse spectral method for the {KP I}
  equation without the zero mass constraint}, Math. Proc. Camb. Phil. Soc., 125
  (1999), pp.~113--138.

\bibitem{Fornb}
{\sc B.~Fornberg}, {\em A practical {G}uide to pseudospectral {M}ethods},
  vol.~1 of Cambridge Monographs on Applied and Computational Mathematics,
  Cambridge University Press, Cambridge, 1996.

\bibitem{FJ}
{\sc M.~Frigo and S.~G. Johnson}, {\em {FFTW} for version 3.2.2}, July 2009.

\bibitem{FMB}
{\sc U.~Frisch, T.~Matsumoto, and J.~Bec}, {\em Singularities of {E}uler flow?
  {N}ot out of the blue!}, J. Statist. Phys., 113 (2003), pp.~761--781.
\newblock Progress in statistical hydrodynamics (Santa Fe, NM, 2002).

\bibitem{GK}
{\sc T.~Grava and C.~Klein}, {\em Numerical solution of the small dispersion
  limit of {K}orteweg de {V}ries and {W}hitham equations}, Comm. Pure Appl.
  Math., 60 (2007), pp.~1623--1664.

\bibitem{GK08}
\leavevmode\vrule height 2pt depth -1.6pt width 23pt, {\em Numerical study of a
  multiscale expansion of {K}d{V} and {C}amassa-{H}olm equation}, Comm. Pure
  Appl. Math., 60 (2007), pp.~1623--1664.

\bibitem{GK12}
\leavevmode\vrule height 2pt depth -1.6pt width 23pt, {\em Numerical study of
  the small dispersion limit of the {K}orteweg-de {V}ries equation and
  asymptotic solutions}, Physica D, 60 (2012), pp.~1623--1664.

\bibitem{KP}
{\sc B.~B. Kadomtsev and V.~I. Petviashvili}, {\em On the stability of solitary
  waves in weakly dispersing media}, Sov. Phys. Dokl., 15 (1970), pp.~539--541.

\bibitem{KassT}
{\sc A.-K. Kassam and L.~Trefethen}, {\em Fourth-{O}rder {T}ime-{S}tepping for
  stiff {PDE}s}, SIAM J. Sci. Comput, 26 (2005), pp.~1214--1233.

\bibitem{ckkdvnls}
{\sc C.~Klein}, {\em Fourth order time-stepping for low dispersion
  {K}orteweg-de {V}ries and nonlinear {S}chr{\"o}dinger {E}quation}, Electronic
  Transactions on Numerical Analysis, 39 (2008), pp.~116--135.

\bibitem{KR}
{\sc C.~Klein and K.~Roidot}, {\em Fourth order time-stepping for
  {K}adomtsev-{P}etviashvili and {D}avey-{S}tewartson equations}, SIAM J. Sci.
  Comp.,  (2011).

\bibitem{KSM}
{\sc C.~Klein, C.~Sparber, and P.~Markowich}, {\em Numerical {S}tudy of
  oscillatory {R}egimes in the {K}adomtsev-{P}etviashvili {E}quation}, J. Nonl.
  Sci., 17 (2007), pp.~429--470.

\bibitem{Kod1}
{\sc Y.~Kodama}, {\em A method for solving the dispersionless kp equation and
  its exact solutions.}, Physics Letters A, 129(4) (1988), pp.~223--226.

\bibitem{Kod2}
{\sc Y.~Kodama and J.~Gibbons}, {\em A method for solving the dispersionless kp
  hierarchy and its exact solutions. ii}, Physics Letters A, 135(3) (1989),
  pp.~167--170.

\bibitem{MS08}
{\sc S.~Manakov and P.~Santini}, {\em On the solutions of the dkp equation: the
  nonlinear riemann hilbert problem, longtime behaviour, implicit solutions and
  wave breaking}, J. Phys. A, 41 (2008), p.~055204.

\bibitem{MS12}
\leavevmode\vrule height 2pt depth -1.6pt width 23pt, {\em Wave breaking in
  solutions of the dispersionless kadomtsev-petviashvili equation at finite
  time}, Theor. Math. Phys., 172 (2012), pp.~1118--1126.

\bibitem{MS06}
{\sc S.~V. Manakov and P.~M. Santini}, {\em Cauchy problem on the plane for the
  dispersionless kadomtsev-petviashvili equation}, JETP Letters, 83(10) (2006),
  pp.~462--466.

\bibitem{MS07}
\leavevmode\vrule height 2pt depth -1.6pt width 23pt, {\em A hierarchy of
  integrable partial differential equations in dimension 2 + 1, associated with
  one-parameter families of vector fields}, Teoret. Mat. Fiz., 152(1) (2007),
  pp.~147--156.

\bibitem{MBF}
{\sc T.~Matsumoto, J.~Bec, and U.~Frisch}, {\em The analytic structure of 2{D}
  {E}uler flow at short times}, Fluid Dynam. Res., 36 (2005), pp.~221--237.

\bibitem{MST}
{\sc L.~Molinet, J.-C. Saut, and N.~Tzvetkov}, {\em Remarks on the mass
  constraint for {KP} type equations}, SIAM J. Math. Anal., 39 (2007),
  pp.~627--641.

\bibitem{NMPZ}
{\sc P.~L. Novikov~S., Manakov~S.V. and Z.~V. E.}, {\em Theory of Solitons: The
  Inverse Scattering Method}, 1984.

\bibitem{Pen76}
{\sc R.~Penrose}, {\em Nonlinear gravitons and curved twistor theory}, Gen.
  Rel. Grav., 7 (1976), pp.~31--52.

\bibitem{PS98}
{\sc M.~Pugh and M.~Shelley}, {\em Singularity {F}ormation in thin {J} with
  {S}urface {T}ension}, Comm. Pure Appl. Math., 51 (1998), pp.~733--795.

\bibitem{Rai12}
{\sc A.~Raimondo}, {\em Frobenius manifold for the dispersionless
  kadomtsev-petviashvili equation}, Comm. Math. Phys., 311 (2012),
  pp.~557--594,.

\bibitem{RLSS}
{\sc G.~D. Rocca, M.~C. Lombardo, M.~Sammartino, and V.~Sciacca}, {\em
  Singularity tracking for {C}amassa-{H}olm and {P}randtl's equations}, Appl.
  Numer. Math., 56 (2006), pp.~1108--1122.

\bibitem{Schme}
{\sc T.~Schmelzer}, {\em The fast evaluation of matrix functions for
  exponential integrators}, PhD thesis, Oxford University, 2007.

\bibitem{SCE96}
{\sc D.~Senouf, R.~Caflisch, and N.~Ercolani}, {\em Pole dynamics and
  oscillations for the complex {B}urgers equation in the small-dispersion
  limit}, Nonlinearity, 9 (1996), pp.~1671--1702.

\bibitem{SSF}
{\sc C.~Sulem, P.~Sulem, and H.~Frisch}, {\em Tracing complex singularities
  with spectral methods}, J. Comp. Phys., 50 (1983), pp.~138--161.

\bibitem{SSP}
{\sc P.-L. Sulem, C.~Sulem, and A.~Patera}, {\em Numerical simulation of
  singular solutions to the two-dimensional cubic {S}chr\"odinger equation},
  Comm. Pure Appl. Math., 37 (1984), pp.~755--778.

\bibitem{tref}
{\sc L.~Trefethen}, {\em Spectral {M}ethods in {MATLAB}}, vol.~10 of Software,
  Environments, and Tools, Society for Industrial and Applied Mathematics
  (SIAM), Philadelphia, PA, 2000.

\bibitem{WRS77}
{\sc R.~Ward}, {\em On self-dual gauge fields}, Phys. Lett. 61A, 81-2 (1977).

\bibitem{Wei}
{\sc J.~Weideman}, {\em Computing the {D}ynamics of {C}omplex {S}ingularities
  of {N}onlinear {PDE}s}, SIAM J Applied Dynamical Systems, 2 (2003),
  pp.~171--186.

\bibitem{ZK69}
{\sc E.~Zabolotskaya and R.~Khokhlov}, {\em Quasiplanar waves in nonlinear
  acoustics of bounded beams}, Sov. Phys. Acoust., 15 (1969), p.~35Ð40.

\bibitem{Zak2}
{\sc V.~E. Zakharov}, {\em Dispersionless limit of integrable systems in 2 + 1
  dimensions}, In Singular limits of dispersive waves (Lyon, 1991), NATO Adv.
  Sci. Inst. Ser. B Phys., 320 (1994), pp.~165--174.

\end{thebibliography}

\end{document}